\documentclass[letterpaper,12pt,leqno,oneside]{article}
\usepackage{color}
\usepackage{bm}
\usepackage{exscale}
\usepackage{amsmath}
\usepackage{amsfonts}
\usepackage{pagecolor,lipsum}
\usepackage{xcolor}
\usepackage{textcomp}
\usepackage{wasysym}
\usepackage{stmaryrd}
\usepackage{amscd}
\usepackage{graphicx}
\usepackage{amsxtra}
\usepackage{amssymb}
\usepackage{theorem}
\usepackage[final]{epsfig}
\usepackage{eqnarray}
\newtheorem{proposition}{Proposition}[subsection]
\newtheorem{definition}[proposition]{Definition}

\newtheorem{lemma}[proposition]{Lemma}
{\theorembodyfont{\rmfamily}\newtheorem{remark}[proposition]{Remark}}

\newtheorem{corollary}[proposition]{Corollary}

{\theorembodyfont{\rmfamily}}

\newfont{\abc}{cmtt10 scaled 1200}

\def\R{\mathbb{R}}

\def\B{\mathbf{A}}
\def\L{\mathbb{L}}
\def\T{\mathbb{T}}
\def\Z{\mathbb{Z}}
\def\Q{\mathbb{Q}}

\def\P{\mathbb{P}}

\def\P{\mathbb{P}}

\def\E{\mathbb{E}}

\def\F{\mathbb{F}}

\def\s{\mathbb{S}}
\def\sp{\mathbb{S}^{\bigoplus}}
\def\st{\mathbb{S}^{\bigotimes}}

\def\I{\mathbb{I}}
\def\ve{\varepsilon}

\def\ra{\rightarrow}
\def\cs{\symbol{35}}
\def\p{\partial}
\def\qed{\hfill $\Box$ \\}
\def\mm{\mbox}
\def\v{= \emptyset}
\def\n{\neq \emptyset}
\def\D{\mathbf{ID}}
\def\M{\mathbb{A}}
\def\RH{reg^{|A|}_H}

\def\bp{\langle A \rangle}

\def\CC{\mm{\textcolonmonetary\,}}
\def\CCC{{\tiny \mm{\textcolonmonetary\,}}}
\def\db{d^{\,\flat}}

\def\llb{\llbracket}
\def\rrb{\rrbracket}

\def\bp{\langle A \rangle}

\begin{document}
 \vspace*{-0.6cm}

\begin{center}\Large{\bf{Skin Structures in Scalar Curvature Geometry}}\\
\smallskip
{\small{by}}\\
\smallskip
\large{\bf{Joachim Lohkamp}}
\end{center}

\vspace{0.2cm}
\noindent Mathematisches Institut, Universit\"at M\"unster, Einsteinstrasse 62, Germany\\
 {\small{\emph{e-mail:  j.lohkamp@uni-muenster.de}}}
\vspace{0.4cm}

 {\small{\center \tableofcontents}

{\contentsline {subsection}{\numberline {}References}{. .}}}

\medskip
%{\footnotesize \center \tableofcontents}

\setcounter{section}{1}
\renewcommand{\thesubsection}{\thesection}
\subsection{Introduction} \label{introduction}
\medskip

In this paper we develop techniques to derive obstruction and structure theories for scalar curvature constraints,  without imposing dimensional or topological restrictions  on the underlying manifold.  To this end, we use skin structures* to merge minimal hypersurface methods with surgery style arguments.\\

We start with an informal review of typical problems and results in scalar curvature geometry to also introduce non-experts to the subject.
Then, in the next two sections, we formulate the main theorems proved in this paper, recollect some basic tools and describe the main steps in the proofs.\\

*Skin structures have been developed in the forerunner papers  [L1] and [L2]. The present one is the third paper in this series. For the readers convenience, we recall some skin structural language from [L1] and [L2], in Ch.s \ref{mre} and \ref{ma}. For a global perspective, the reader may wish to consult the overview of this series, in [L1], Ch.1.1.

\subsubsection{Classical and New Techniques} \label{ov1}
\bigskip

The scalar curvature, $scal(g)$ of a Riemannian manifold $(M,g)$,  is one of the basic curvature notions, besides sectional and Ricci curvature:\\

$\bullet$ \, Scalar curvature measures the infinitesimal growth of the volume of distance balls: the volume $vol_{(M,g)}(B_{\ve}(p))$ of a
        distance ball $B_{\ve}(p)$ of radius $\ve > 0$ in  $(M,g)$ of dimension $n$, $p \in M$ relates to the Euclidean volume on $\R^n$ according to \[
 \frac{vol_{(M.g)}(B_{\ve}(p))}{vol_{(\R^n,g_{Eucl.})}(B_{\ve}(p))} =  1 -   c_n \cdot scal(g)(p) \cdot \ve^2 + O(\ve^4),\]
 for $\ve >0$ and $c_n:= 1/6(n+2)$, cf.[G1], Ch.1. This makes $scal(g)$  the simplest extension of
the Gau\ss\ curvature for surfaces to higher dimensions.\\

$\bullet$ \, The assignment $g  \mapsto scal(g)$ can be concisely characterized as the unique second order differential operator acting on
metrics to functions, both on $M$, which is linear in second derivatives of $g$ and equivariant under the natural action of diffeomorphisms of $M$ on metrics and functions,
normalized to $scal(S^n)\equiv1$, cf.[G1], Ch.1.\\

$\bullet$ \, Scalar curvature occurs in a number of distinguished places in geometry and physics, like that of Einstein's field equations: \[Ric(g)_{\mu
        \nu} - {\textstyle \frac{1}{2}}\cdot scal(g) \cdot g_{\mu \nu} = \kappa \cdot T_{\mu \nu},\]  where $Ric$ is the Ricci curvature, $T$  the stress-energy tensor and $\kappa >0$ a
        coupling constant. The impact of $scal(g)$ becomes even more obvious from Hilbert's way to derive these field equations, namely, from minimizing the action  \[W[g]:=\int
        \big(scal(g) + {\cal{L}}_M \big)\cdot \sqrt{det(g)} \, dV_g,\]  where now ${\cal{L}}_M$ represents the matter/energy fields, cf.[D].\\

In view of such a prominence, it is a puzzling fact that significant consequences of scalar curvature bounds are hard to detect. In turn, it is then not
surprising that those implications known are striking and far-reaching.\\

In this paper we introduce new techniques to study scalar curvature constraints. To explain our results and to give the reader an orientation towards the main challenges, we start
with a brief review of the essential techniques one currently knows and indicate where they reach their limits.\\

For this, we first notice that we only need to consider lower curvature bounds since every Riemannian metric can
 be $C^0$-approximated by others with arbitrarily strongly negative scalar curvature cf. [L3], [L4].\\

The key case to understand the impact of lower scalar curvature bounds is that of positive scalar curvature, $scal > 0$. The treatment of other
cases can oftentimes be reduced to, or at least be inspired from, this case, cf.[G1] for a rich exposition.\\

\textbf{A. $\mathbf{Scal >0}$ and Spin Geometry} \,   It was only the Atiyah-Singer Index theorem, developed in the 60ies, that enabled Lichnerowicz to detect the first
known obstruction for the presence of $scal> 0$-metrics in dimensions $\ge 3$.\\

He considered the Dirac operator on a spin manifold and related it to scalar curvature through some Weitzenb\"ock identity. Then he used the Index theorem to show that $scal> 0$ implies a vanishing $\hat{A}$-genus, cf. [Li].\\

 This idea was further developed to derive many significant obstruction and structural results, most notably by Gromov and Lawson (cf. [GL1], [GL2],
[G1] and [LM]). For instance, they proved that so-called enlargeable manifolds, like $T^n \cs M^n$, for any closed spin manifold $M^n$, \emph{cannot} admit any $scal> 0$-metric.\\

Also within the class of spin manifolds, Witten used the same Weitzenb\"ock  identity to derive a short proof of the positive mass conjecture in general relativity, which is closely related to scalar curvature problems, cf. [W],[PT].\\

However, due to the pervasive use of spin structures these techniques do not work in the non-spin case. Towards a general structure theory for $scal >0$-metrics, this is a substantial issue, since the class of non-spin manifolds embraces the majority of all manifolds.\\

\textbf{B. $\mathbf{Scal >0}$ and  Minimal Hypersurfaces}  \,  In the late 70s, a differently flavored approach to scalar curvature constraints, without a priori topological limitations, was initiated by Schoen and
Yau  in [SY1] - [SY3], [S] and broadened by Gromov and Lawson [GL3]. It uses minimal hypersurfaces and started from the following remarkable $scal >0$\emph{-heredity} result:\\

\emph{Any area minimizing hypersurface $(H^n,g_H)  \subset (M^{n+1},g_M)$ within some smooth compact $scal > 0$-manifold $M^{n+1}$ picks up positive scalar curvature from its ambience so that, after canonical conformal deformations, its scalar curvature is also positive.}\\

For the readers convenience and also for later reference, we briefly explain this powerful yet simple inheritance result.\\

The first and the second variation of $Area(H)$ of a hypersurface $H \subset M$ by a $f \cdot \nu$, where $\nu$ is the outward normal vector field of $H$,
$f \in C^\infty(H,\R)$ with $supp \: f \subset reg(H)$ (= the set of
regular points of $H$), are given by {\small \[Area'(f) =  \int_H tr A_H (z) \cdot f(z) \: dVol\]
\[Area''(f) =  \int_H |\nabla_H f|^2  + \left( (tr A_H)^2 -  |A|^2 - Ric(\nu,\nu) \right) \cdot f^2 \: dVol\]}
$tr A_H$ is the mean curvature of $H$, $|A|^2$ is the sum of the squares of principal curvatures of $H$ and $Ric$ the Ricci tensor of $M$, cf. [G1], p.54.\\

Since $H$ is supposed to be area minimizing, we have $Area'(f) = 0 $ and $Area''(f) \ge 0 $. That is $tr A_H = 0$ and this gives {\small \begin{equation}\label{1}  \quad Area''
(f) = \int_{H}|\nabla_H f|^2  -  \left( |A|^2 + Ric(\nu,\nu)  \right) \cdot f^2 \: dVol \ge 0\end{equation}} Now we recall the Gau\ss-Codazzi equations for hypersurfaces {\small
\begin{equation}\label{gc} \quad |A|^2 + Ric(\nu,\nu)  =1/2 \cdot \left(| A |^2 + scal_M - scal_H +(tr A_H)^2 \right)\end{equation}} where $scal_H$ and $scal_M$ denote the scalar
curvature of $H$ and $M$. Since $tr A_H = 0$ we may rewrite (\ref{1}) as follows {\small
\begin{equation} \label{3} \int_H | \nabla f |^2 + \frac{n-2}{4 (n-1)} \cdot scal_H \cdot f^2 d A \, \, \ge \end{equation} \[ \int_H \frac{n}{2 (n-1)} \cdot  |  \nabla f |^2 + \frac{n-
2}{4 (n-1)}\cdot \left( | A |^2 + scal_M \right) \cdot f^2 d A.\]}

The left hand side of $(\ref{3})$ is the variational integral for the first eigenvalue $\lambda^L_1$ of the conformal Laplacian $L_H$ on $H$:\quad $L_H(f) = -\Delta f
+\frac{n-2}{4 (n-1)} \cdot scal_H \cdot f,$ {\small $$\lambda^L_1 = \inf \{\int_H | \nabla f |^2 + \frac{n-2}{4 (n-1)}  \cdot  scal_H  \cdot  f^2 d A \, |\, f \mm{ smooth, supp}\, f
\subset reg(H), |f|_{L^2}=1 \}$$} Thus, when $scal_M > 0$, inequality $(\ref{3})$ shows that $\lambda^L_1 > 0$. This was first realized by Schoen and Yau [SY]. \\

Then, they noted that according to Kazdan and Warner [KW], the transformation law for scalar curvature under conformal transformation by the first eigenfunction
$f_1>0$ can be used to show that
\begin{equation} \label{4} scal(f_1^{4/n-2} \cdot g_H) \cdot f_1^{\frac{n+2}{n-2}} =
L_H(f_1) = \lambda^L_1 \cdot f_1 > 0.\end{equation}    Therefore, we conclude that with $M$, any of its area minimizing hypersurfaces $H$ carries a $scal > 0$-metric, in
the given case, $f_1^{4/n-2} \cdot g_H$.\\

\textbf{C. Towers of Dimensional Descent} \, What made this $scal >0$-heredity so exciting was the perspective to study scalar curvature using strategies of inductive
dimensional descent along towers of nested area minimizing hypersurfaces equipped with inherited curvature constraints until one reaches a space already understood.\\

As an example, we consider the manifolds $T^n \cs M^n$ we encountered earlier, but now we no longer assume that $T^n \cs M^n$ is spin. We start from a hypothetical $scal >0$-metric on $T^n \cs M^n$ and iteratively choose a tower of hypersurfaces each of them admitting a $scal >0$-metric descending to a surface.\\

 Concretely, we get an area minimizing hypersurface $N^{n-1}$ in $T^n \cs M^n$, deform $N^{n-1}$ to $scal >0$ and choose an area minimizer in this geometry and so on.The topology of $T^n \cs M^n$ allows us to select the homology classes so that we finally reach a surface $F^2$ of genus $\ge 1$.  But the Gau\ss-Bonnet theorem says that $F^2$ cannot admit a $scal >0$-metric. Thus there is no $scal >0$-metric on $T^n \cs M^n$.\\

This type of strategy was, and it still is, the only known candidate for general structure theories of scalar curvature beyond the realm of spin geometry.\\

However, as a kind of retributive justice, the applicability of this method now hinges upon the regularity of area minimizers. The infamous problem is that these hypersurfaces can  contain hardly comprehensible singular sets of codimension $\ge 7$.\\

 While, in dimensions $\le 7$, this means that all such hypersurfaces are smooth and, hence, this strategy works flawlessly, the situation changes in higher dimensions. Despite their apparently small size, these singularities have a deep impact on their regular complement and on the applicability of this inductive $scal >0$-heredity argument.
\begin{itemize}
  \item   The hypersurfaces, and typical elliptic operators defined on them, including the conformal Laplacian, degenerate towards the singular set. In this respect, the high codimension of the singularities may be more harmful than useful.
  \item  In towers of such nested singular spaces, essential parts of the area minimizers, or even the entire minimizer, may disappear in the singular set of their forerunners. Then the tower trivializes and we do not gain any information.
  \item  When these nested hypersurfaces remain non-trivial, the issues may pass to dimensions $\le 6$, since area minimizers in singular spaces may carry more subtle  singularities of smaller codimension. This is paralleled by an increasingly delicate analysis on these spaces.
    \item It is unclear how to extract obstruction results from a $scal >0$-heredity that involves singular spaces, since every Riemannian metric on every manifold $M^n$, $n\ge 2$, can be $C^0$-approximated by \emph{positively} scalar curved metrics which are smooth, except for some kind of orbifold singularities of codimension $2$, with angles arbitrarily close to $2 \cdot \pi$. This can be shown using tools of [L4].
\end{itemize}

This sort of problems made the approach intractable in higher dimensions and, as a consequence, there were no generally applicable obstruction results for scalar curvature constraints beyond dimension $7$. \\

\textbf{D. $\mathbf{Scal >0}$-Heredity with Surgery} \,  In this paper we establish a broader strategy of dimensional descent. Besides a  $scal >0$-heredity, the extended scheme  includes mechanisms to substitute the singularities for regular ends and closures keeping the $scal >0$-condition intact. This way we can regularize the towers of inductive descent. We call this generalized scheme a $scal >0$\emph{-heredity with surgery}. \\

With these methods we resolve the issues with the classical approach and get access to study scalar curvature constraints without the classical topological and dimensional restrictions, we have seen above.\\

Applications, but also extensions of these techniques to cover more general geometric constraints, are treated in separate accounts. To mention some basic results, this approach leads to proofs of the general versions of the positive mass theorems in general relativity, the Penrose inequality and the fact that enlargeable manifolds cannot admit $scal >0$-metrics.\\

 In the next section, Ch.\ref{sor}, we formulate our main results, in Ch.\ref{mre} we give an overview of the proof of these results which broadly relies on skin structural concepts and results from [L1] and [L2].

\subsubsection{Statement of  Results } \label{sor}
\bigskip

Henceforth, we assume we are given  a compact area minimizing hypersurface  $(H^n,g_H) \subset (M^{n+1},g_M)$ within some  smooth $scal \ge  0$-manifold  $M^{n+1}$, with some \emph{non-empty} singular set $\Sigma_H \subset H$, $\Sigma_H \n$.\\

  Our basic desingularization result for singular area minimizers in scalar curvature geometry is the following one. \\

\textbf{Theorem 1} \textbf{(Regular Open Ends)} {\itshape \, For any  $\ve > 0$, there is a smoothly bounded neighborhood $W_{\ve}$ of $\Sigma$ within an $\ve$-neighborhood $U_{\ve}$ of
$\Sigma$ such that $H \setminus W_{\ve}$ admits a smooth metric $g^+_{\ve}$  so that
\begin{enumerate}
\item \quad $g^+_{\ve}$ is conformal to $g_H$, \item  \quad $scal (g^+_{\ve}) > 0$ on $H \setminus W_{\ve}$, \item  \quad $\p W_{\ve}$ has positive mean curvature.
\end{enumerate}}

The positive mean curvature $\p W_{\ve}$ has a deflecting effect on area minimizers \emph{within} $(H \setminus W_{\ve},g^+_{\ve})$. In fact, these minimizers in $(H \setminus W_{\ve},g^+_{\ve})$ never touch $\p W_{\ve}$. To them, $H \setminus \Sigma$ appears to be smooth and compact, whereas $\Sigma$  resides beyond an \emph{inner horizon}, namely $\p W_{\ve}$, within $H \setminus \Sigma$, such area minimizers never pass. More precisely, we have\\

\textbf{Theorem 2} \textbf{(Smooth Ambience I)} {\itshape \, For any  closed smooth hypersurface $L \subset H \setminus \Sigma$ there is some $\ve_L > 0$ such that for
any $\ve \in (0,\ve_L)$ we get an area minimizing hypersurface $L_{\ve} \subset (H \setminus W_{\ve},g^+_{\ve})$ so that $$L_{\ve} \mm{ is  homologous to } L \mm{ and }
L_{\ve}  \cap \p W_{\ve} = \emptyset.$$ Thus $L_{\ve}$ is an area minimizer within a smooth $scal >0$-curved ambient space.}\\

 Since $\bigcup_{\ve > 0}(H \setminus W_{\ve},g^+_{\ve}) = H \setminus \Sigma,$ we observe that the $(H \setminus W_{\ve},g^+_{\ve})$ eventually capture any detail
of $H$ compactly supported in $H \setminus \Sigma$.\\

 \textbf{Remark 1} \, Theorem 2 can be applied in situations where we can iteratively select homology classes, so that relevant information can be located away from the singular set. A sample application of Theorem 2 is the Riemannian version of the positive mass theorem (and of the Penrose inequality) where we reproduce the asymptotic flatness in each inductive loop.\\

  In more general cases, all representatives of some $n-1$-dimensional homology class of $H$ may hit $\Sigma$. To handle this case, we deform $H$ a third time, but in a non-conformal way, to get some additional control over $H$ towards the open ends.\\

\textbf{Theorem 3} \textbf{(Regular Closed Doublings)} {\itshape\, For any  $\ve > 0$,  there is a smoothly bounded neighborhood $V_{\ve}$ of $\Sigma$ contained in an
$\ve$-neighborhood $U_{\ve}$ of $\Sigma$ so that the doubling
\[H^D_{\ve}:= H \setminus  V_{\ve}  \cup_\sim H \setminus  V_{\ve} ,\mm{ where } \sim \mm{ means gluing along } \p V_{\ve}\]
 admits a smooth $scal>0$-metric $g^*_{\ve}$  with
\begin{enumerate}
 \item$g^*_{\ve}|_{H \setminus V_{\ve}}$ is conformal to the induced metric on $H \setminus V_{\ve}\subset M$, except for a kind of
    generalized warped product deformation (arbitrarily) near $\p V_{\ve}$.
    \item $H^D_{\ve}$ is mirror symmetric relative $\p V_{\ve}$ and, hence, $\p V_{\ve}$ is totally geodesic.\\
\end{enumerate}}

As an application we get an extension of Theorem 2 where the resulting minimizers in $H$ approximate any given homology class of the singular total space $H$ while they are properly supported in  $H \setminus \Sigma$.\\

\textbf{Theorem 4} \textbf{(Smooth Ambience II)} {\itshape \,  Let $\beta \in H_{n-1}(H;\Z)$ be given. Then, for any $\ve
> 0$, we can find a hypersurface $N_{\beta,\ve} \subset H \setminus  V_{\ve}$, so that
\begin{enumerate}
\item  $N_{\beta,\ve}$ is an integral current and mass minimizing under perturbations with support in the interior of $H \setminus  V_{\ve}$
\item $N_{\beta,\ve}$ hits $\p V_{\ve}$ perpendicularly, $\p N_{\beta,\ve} \subset \p V_{\ve}$  is a minimal hypersurface and $\p N_{\beta,\ve} \subset N^D_{\beta,\ve}$ is
    totally geodesic.
\item  $N_{\beta,\ve}$  represents the image of $\beta$ under the projection map $\pi_*$ composed with the isomorphism ${I_{H,V_{\ve}}}$ between $H_*$ and $H^{IC}_*$ on
    the couple  $(H,V_{\ve})$:
\[\begin{CD}
\beta \in H_{n-1}(H;\Z) @>{\pi_*}>>  H_{n-1}(H,V_{\ve};\Z) @>{I_{H,V_{\ve}}}>> [N_{\beta,\ve}] \in H^{IC}_{n-1}(H,V_{\ve};\Z).
\end{CD}\]
\end{enumerate}}

 \textbf{Remark 2} \, $H_*$ and $H^{IC}_*$ denote the ordinary singular homology respectively the homology of integral currents. The excision axiom shows that $(H,V_{\ve})$ can be viewed as a smooth couple and thus ${I_{H,V_{\ve}}}$ becomes an isomorphism, cf.[Dp] or [H]. We can use $I_{H,V_{\ve}}$ to alternatively interpret  $N_{\beta,\ve}$ as a  representative of  $H_{n-1}(H,V_{\ve};\Z)$
or, when we append the isomorphism $ H_{n-1}(H,V_{\ve};\Z) \cong  \widetilde{H}_{n-1}(H/V_{\ve};\Z)$, as a representative of  $\widetilde{H}_{n-1}(H/V_{\ve};\Z)$.\\

When $\beta$ can be represented by some hypersurface $P \in \beta$ with $P \cap \Sigma \v$, then we have $\p N_{\beta,\ve} \v$. This reproduces Theorem 2.\\

 \textbf{Conclusion} \, We can apply the deformation of Theorems 2 respectively 4 to replace the use of singular towers of inductive descent for that of  \emph{smooth towers} with unreachable boundary or smooth replacements for $\Sigma$, in the case of doublings.

\subsubsection{Outline of the Strategy} \label{mre}

We first motivate and sketch the contents of this paper, review some of the main ingredients and finally we give a chronological overview of the later chapters.\\

\textbf{A. Open Ends versus $scal >0$} \, To understand the origins of the approach we have chosen in this paper and the problems we encounter, we first look at what happens when we use the conformal deformation via eigenfunctions $u>0$ of the conformal Laplacian $L_H$ of \ref{ov1}.B above.\\

On a singular hypersurface $H$ we get  a $scal >0$-metric also singular along $\Sigma \subset H$. One finds that $u>0$ has minimal growth towards $\Sigma$, when compared to other solutions locally defined around $\Sigma$, and the resulting metric space $(H,u^{4/(n-2)} \cdot g_H)$ is again compact. The notorious problem is that there is no reason why an area minimizing hypersurface within $(H,u^{4/(n-2)} \cdot g_H)$ should not reach or even partially shrink to $\Sigma$.\\

A suggestive idea is to deform $u$, in particular, in a neighborhood of $\Sigma$, into another function $v>0$ so that $(H,v^{4/(n-2)} \cdot g_H)$ is still a $scal >0$-metric, but with
a \emph{spreading open end} towards $\Sigma$. This end structure keeps area minimizers within $H \setminus \Sigma$ from approaching $\Sigma$. In other words, $\Sigma$  lies beyond some \emph{inner horizon}. The classical codimension $\ge 3$ surgeries along submanifolds of $scal >0$-manifolds keeping $scal >0$ make this a plausible option, cf.[GL1].\\

However, such an end structure can only be accomplished at the expense of sacrificing some amount of positive scalar curvature. Here we reach the point where we need controls and estimates, in a situation of hardly any control over $\Sigma$ and the behavior of $H$ near $\Sigma$, not to mention that of $L_H$ or $u$ one such unruly grounds.\\

 At least,  there are tangent cones approximating $H$ around any given singular point as partially linearized and, hence simplified, local models. This suggests that shifting the core of the support of the intended deformation towards $\Sigma$ leads to a better control. \\

But, there is an adversarial effect: the equation $-\Delta u
+\frac{n-2}{4 (n-1)} \cdot scal_H \cdot u = \lambda_1 \cdot u$ is not scaling invariant. This means, when we approach $\Sigma$, a rescaling of ever smaller balls around singular points to unit size
 shrinks  the value $\lambda_1$ and eventually it converges to zero. But it is this value that stands for the available amount of $scal >0$ we may invest to build an open end. In other words, we gradually lose the flexibilty to build open ends without leaving the $scal >0$-condition.\\

 Summarizing, at this stage, we are missing a control over the asymptotic analysis of $u$ and cannot improve the situation from approaching the singular set.\\

  In what follows we describe a way out of both problems. We will also observe that these two problems are tightly coupled. This analysis, like the entire approach, is based on skin structures we describe now.\\

\textbf{B. Skin Structures} \, On any area minimizing hypersurface $H$ we consider a skin transform $\bp_H$. That is, a non-negative measurable function defined on $H \setminus \Sigma$ that
commutes with the convergence of sequences of area minimizers and satisfies the following axioms
\begin{itemize}
    \item $\bp_H \ge |A_H|$ and for any $f \in C^\infty(H \setminus \Sigma,\R)$ compactly supported in $H \setminus \Sigma$ and some $\tau = \tau(\bp,H) \in (0,1)$ we have the Hardy type inequality
\begin{equation}\label{hi}
\int_H|\nabla f|^2  + |A_H|^2 \cdot f^2 dA \ge \tau \cdot \int_H \bp_H^2\cdot f^2 dA.
\end{equation}
    \item $\bp_H \equiv 0$, if $H \subset M$ is totally geodesic, otherwise, $\bp_H$ is  strictly positive.
    \item When $H$ is not totally geodesic, we define $\delta_{\bp}:=1/\bp$, the \emph{${\bp}$-distance}. It is $L_{\bp}$-Lipschitz regular, for some constant
        $L_{\bp}=L(\bp,n)>0$:
        \[|\delta_{\bp}(p)- \delta_{\bp}(q)|   \le L_{\bp} \cdot d(p,q), \mm{ for } p,q \in  H \setminus \Sigma. \]
\end{itemize}
The role of $\delta_{\bp}$ is that of a generalized distance function to $\Sigma$. Up to some normalization, $\delta_{\bp}(p)$ can be thought as the guessed distance to $\Sigma$
applicable even when $\Sigma \v$. On the other hand, $\bp$ can be viewed as a particular type of averaged version $|A|$. Namely, each level set $|A|^{-1}(c)$, $c >0$, of $|A|$ is transformed into a regularized hypersurface $\M_c$  spanned over  $|A|^{-1}(c)$ using some fixed algorithm/procedure. Then $\bp$ is defined as the canonical extension of the definition $\bp^{-1}(c):=\M_c$.\\

In [L1] we described procedures to realize such a function. Indeed, we only exploit the few axioms above and for the remainder of this paper we may choose any such skin transform $\bp$. \\

\textbf{C. Skin Uniformity and Potential Theory} \, As a first application of this concepts, in [L1], we have shown that
 $H \setminus \Sigma$ is a \emph{skin uniform space}: for any two $p,q \in H$, there is a rectifiable path $\gamma: [a,b] \ra H$, for some $a <b$, with $\gamma(( a,b))\subset H \setminus \Sigma$,
so that for any given skin transform $\bp$, there is some $s_H \ge 1$ with \begin{itemize}
    \item \quad  $l(\gamma)  \le s_H \cdot d_{g_H}(p,q)$,
    \item \quad $l_{min}(\gamma(z)) \le s_H \cdot \delta_{\bp}(z)$, for any $z \in \gamma_{p,q}.$
    \end{itemize}
This also implies that  $H \setminus \Sigma$ is a \emph{uniform space} as considered in quasi-conformal geometry,cf.[He]. Such uniformity properties can be thought as boundary regularity results for $\Sigma$, viewed as the boundary of $H \setminus \Sigma$.\\

It is the stronger skin uniformity we need to understand the analysis of elliptic operators on $H \setminus \Sigma$ towards  $\Sigma$. This is explained at length in [L2].  Specifically, we get insights in the potential theory  of \emph{skin adapted operators}, that is, linear second order elliptic operators $L$, so that
\emph{\begin{itemize}
     \item $L$ does not degenerate faster than $\bp^2$ when we approach $\Sigma$.
     \item $L$ satisfies a weak coercivity condition: we require that there is a supersolution $s
         >0$ of $L \, u=0$ and some $\ve >0$, so that: $L \, s \ge \ve \cdot \bp^2 \cdot s.$\end{itemize}}
For instance, we obtain controls for the asymptotic behavior of solutions $u>0$ of $L \, w=0$, when we approach $\Sigma$.  For the Martin boundary  $\p_M(H \setminus \Sigma,L)$ of $L$ on $H \setminus \Sigma$ we have [L2],Th.5: $\p_M(H \setminus \Sigma,L)$ is \emph{homeomorphic} to $\widehat{\Sigma}$ and all Martin boundary points are minimal. ($\widehat{A}$ denotes the one-point compactification of the space $A$.)\\

  Thus,  $u >0$ solves $L \, w=0$ if and only if it can be written as a Martin integral
\begin{equation}\label{mi2}
u(x) = u_\mu(x) = \int_{\widehat{\Sigma}} k_L(x;y) \, d \mu(y).
\end{equation}
$k_L(x;y)$ is the Martin kernel, $\mu=\mu_u$ a Radon measure on $\widehat{\Sigma}$ uniquely associated to $u$.\\

 \textbf{D. Conformal Laplacians} \, The Hardy inequality (\ref{hi}) implies  the following counterpart to (\ref{4}) above, asserting that when $\Sigma_H \n$ and $H \subset M$, for some manifold with $scal \ge  0$, there is a positive $\lambda > 0$ with \begin{equation}
\label{x} \int_H f \cdot L_H f \,  dA \, \ge \, \lambda \cdot \int_H \bp^2\cdot f^2 dA,\end{equation}
for any smooth $f$ with $supp(f) \subset H \setminus \Sigma$. That is, when $M$ has $scal \ge 0$, we know that $L_H$ is  skin adapted.\\

 For $\Sigma \n$, this is a proper and non-trivial refinement of the argument for closed smooth $H$, since $\bp \ra \infty$, when we approach $\Sigma$. Also, it is important that, for $\Sigma \n$, we merely need $scal_M \ge 0$ since, in the approach chosen in this paper, we frequently work with the conformal Laplacian on tangent cones (in the flat $\R^{n+1}$) of $H \subset M$.\\

There is a largest such $\lambda > 0$, the generalized principal eigenvalue $\lambda^{\bp}_{H,L_H} >0$ of the $\bp$-weighted version of $L_H$, the \emph{skin conformal Laplacian} $\delta_{\bp}^2 \cdot L_H$.\\

 By contrast to the case of a compact and smooth hypersurface, we get from [L2],Th.s 5, 9 and 10: For any $\lambda <\lambda^{\bp}_{H,L_H}$, there is a wealth of eigenfunctions $u >0$ of $\delta_{\bp}^2 \cdot L_H$, that is, solutions of
 \begin{equation} \label{y}  - \triangle u  + \frac{n-2}{4 (n-1)}  \cdot scal_H \cdot  u  = \lambda  \cdot  \bp^2  \cdot u.\end{equation}
Towards scalar curvature geometry, we notice that any solution $u>0$ of $L_{H,\lambda} \, w=0$ is $C^{2,\alpha}$-regular, for $\alpha \in (0,1)$, since $scal_H$ is smooth and $\bp^2$ is locally Lipschitz continuous. Thus the scalar curvature of the metric $u^{4/n-2} \cdot g_H$ is well-defined and,  for $\lambda \in (0,\lambda^{\bp}_{H,L_H})$, it has $scal >0$. To couple (\ref{y}) to Martin theory, we write this as \[L_{H,\lambda} \, u = 0,\mm{ for the operator } L_{H,\lambda}:= - \triangle   + \frac{n-2}{4 (n-1)}  \cdot scal_H  -\lambda  \cdot  \bp^2 \]
We observe that $L_{H,\lambda}$ is skin adapted for $\lambda <\lambda^{\bp}_{H,L_H}$ and any iterated tangent cone $C$ inherits the skin adaptedness for $L_{C,\lambda}$ from that of  $L_{H,\lambda}$ on $H$.\\

 \textbf{E. Growth Estimates and Bending Effects} \,  For sufficiently small $\lambda >0$, we can also derive explicit estimates for the growth rates of solutions towards $0 \in C$, cf.  Prop.\ref{ddd} resp. [L2], Th.11. \\

These estimates show that, for  solutions $u_\mu$ associated to measures $\mu$ on $\widehat{\sigma_C}$, with substantial support in $0$,  $u_\mu^{4/n-2} \cdot g_C$ carries an inner horizon shielding the tip $0$ against area minimizers in $C \setminus \{0\}$ (where $\sigma_C$ is the singular set of $C$). For the general case of $\Sigma \subset H$, these estimates suggest to consider, in a sense, \emph{equidistributed} positive Radon measures $\mu$ on $\Sigma$ and their associated solutions
\begin{equation}\label{mic}
u_\mu(x) := \int_{\widehat{\Sigma}} k_{L_{H,\lambda}}(x;y) \, d \mu(y).
\end{equation}
We expect that $(H \setminus \Sigma,u_\mu^{4/n-2} \cdot g_H)$ contains an inner horizon shielding
$\Sigma$, retaining area minimizers within $H \setminus \Sigma$ from approaching $\Sigma$. However, structural inhomogeneities of $H$ along $\Sigma$ make it a delicate task to properly describe or analyze  equidistributed Radon measures on $\Sigma$. Thus, we proceed in a slightly different way and establish a reduction to solvable cases where $\Sigma$ and $H$ are simplicial complexes.\\

\textbf{F. Nanostructures on  $H \setminus \Sigma$} \,  To this end, we use again skin structures  to define large families of small and uniformly shaped deformations. They are compactly supported in $H \setminus \Sigma$ and placed closely to $\Sigma$. It will be suggestive to think of them as \emph{nanostructures} on  $H \setminus \Sigma$. Similar as in nanotechnology, it is the uniform shape and the very small but not infinitesimal size of our deformations we need to get fine details and auto-arranging patterns, not available at larger scales. \\

 What we do with these structures is to break solving the condition - $scal >0$ with spreading open ends - into accessible pieces, combined with an assembly process to obtain such a metric on $H \setminus \Sigma$.\\

Concretely, the basic observation we start from is the  scaling invariance of equation (\ref{y}): the value $\lambda>0$ does not change when we scale $H$ by any constant $c>0$, since all terms transform equally.\\

 Thus, different from the situation in subsection A above, we can zoom into $\Sigma$ and keep the positive bounds on the scalar curvature of $u^{4/n-2} \cdot g_H$ to benefit from the approximation of $H$ by the (partially) linearized model spaces of tangent cones to control (and actually to define) the bending processes. \\

We proceed in two moves, where we use this scaling invariance in various steps analytically but also to build geometric rig to localize the problems:
\begin{itemize}
\item We first consider a \emph{downstream} inheritance of solutions on iterated tangent cones.
We fix  a small $\lambda>0$, some solution $\Phi>0$ of $L_{H,\lambda} \ w=0$ on $H \setminus \Sigma$ and observe that $\Phi$ induces solutions $\Phi(C)$ of $L_{C,\lambda} \ w=0$ on all tangent cones $C^n=\R^{n-k} \times C^k$, for the, again, skin adapted operator $L_{C,\lambda}$.\\

Now we locally deform the $\Phi(C^n)$ along the singular symmetry axes $\sigma^{\cs}_{C^n}=\R^{n-k}\times\{0\}$: we add a solution $v_{\mu}$, where $\mu$ on $\R^{n-k}$, is the prototype of an equidistributed measure along the symmetry axes of $C^n$, namely, the $n-k$-dimensional Lebesgue measure, to get a partial horizon shielding $\sigma^{\cs}_{C^n} \subset \sigma_{C^n}$.

\item The second, technically more delicate step is an \emph{upstream} process.  We select a family of iterated tangent cones represented by a finite tree
$\T$ of iterated blow-ups around singular points, so that each branch ends with cones with only isolated singularities.  Now, closely to $\Sigma$, we define natural pipeline systems $\sp\subset H \setminus \Sigma$ of coordinate patches, the pipes $\st$, surrounding $\Sigma$. Similarly, we get  compatible pipelines $\sp$  surrounding the singular symmetry axes $\sigma^{\cs}_{C^n}$ of each cone in $\T$.\\

On the analytic side we observe and exploit that the $\sp$ match  the scaling invariance of the equation (\ref{y}). We use $\sp$ as a ground plan to appropriately localize the solutions $v_{\mu}$ on each of the cones $C^n \in \T$, fitting into some pipe $\st \subset C^n$ and to transfer these pieces to specified places on $H \setminus \Sigma$ near $\Sigma \subset H$.\\

 The result is a conformal deformation of $H \setminus \Sigma$ that gives us a metric with a global horizon around $\Sigma$, located as closely to $\Sigma$ as we want,  and with $scal >0$ on $H$, on the regular side of the horizon.
\end{itemize}

The overall strategy may be interpreted as a reduction to the case of simplicial complexes. The soul of $\sp$ is a finite simplicial complex $\Sigma^\Delta$, approximating $\Sigma$, whose simplices are pieces of symmetry axes  $\sigma^{\cs}_{C^n}$ of those cones $C^n \in \T$. \\

\textbf{G. Organization of the Paper} \,  In Ch.\ref{fm0} and Ch.\ref{fm1},
we introduce the blow-up trees $\T$ and the associated systems $\sp$. This uses skin structures, in a way independent of  [L1] and [L2], to  make efficient use of tangent cone approximations. Then, in Ch.\ref{smo0}, we define and study locally shielding deformations along the pipes in the cones of $\T$. We assemble them along $\sp$ to define the \emph{horizon metric} $g_{\T}$.\\

Thus the focus of this paper is on the detailed portrait of one extensive construction that eventually leads us to the metric $g_{\T}$.  It is then a matter of understanding  the influence of the parameters of horizon metrics to select them appropriately to make the $g_{\T}$ actually shield $\Sigma$.\\

\textbf{Ch.\ref{fm0}} \,  We start with an upgrade of classical tangent cone approximation results to broader so-called freezing results in Ch.\ref{na} and \ref{na2}. With these techniques, we define self-similar configurations of ball covers of $\Sigma$, with some common approximation accuracy of these balls by balls in some tangent cones, in Ch.\ref{coverings}. \\

 An important aspect, discussed in Ch.\ref{sish}, is that the singular directions of the balls in these cover show some auto-alignment properties comparable to compass needles in common field. We employ this in Ch.\ref{sish1}, to define the basic model of a pipeline.\\

\textbf{Ch.\ref{fm1}} \, Next we observe that tangent cones are a yet too weak concept to gain control over the planned pipeline structures on $H \setminus \Sigma$. To resolve this issue, we introduce in Ch.\ref{gr}, a stabilization through other no longer properly tangent area minimizing cones, the \emph{almost tangent cones}. With them we can turn compactness results for minimizing currents into estimates for the shape of the $\st$ and $\sp$. This is even necessary to properly define any of these structures.\\

Then, in Ch.\ref{alt} and \ref{grq}, we inductively use self-similar covers for $\sigma \cap \p B_1(0) \cap C$, for the singular set $\sigma \subset C$, where $C$ is an almost tangent cone, and so on. The balls in the cover of $\Sigma$ and all the used iterated (almost) tangent cones can be thought as the branching points of the blow-up tree $\T$.\\

\textbf{Ch.\ref{smo0}} \, Here we merge the latter pipeline techniques with results for skin adapted operators from [L2], briefly summarized in Ch.\ref{ma},  and reach the $scal >0$ - heredity with surgery. The argument occupies Ch.\ref{smo} and Ch.\ref{smo0}.\\

In Ch.\ref{ma}-\ref{lll}, we use the theory for $L_{C,\lambda}$, for small $\lambda>0$ on all  cones $C$ in the tree $\T$ to stepwise define deformations of a given solution $v_0>0$ supported in a single $\sp$.
  Following the given tree $\T$ this leads us, in Ch.\ref{lll},  to the definition of the horizon metric $g_{\T}$ on  $H \setminus \Sigma$ and formulate its global shielding effect along $\Sigma$, proved in Ch.\ref{smo}.   In Ch.\ref{smooth} we infer the Theorems from further refinements of the shape of the horizon. \\

  \textbf{Ch.\ref{smo}} \, In Ch.\ref{geom} and \ref{wm} we consider the shielding effects of the individual elements of $g_{\T}$ in the total pipeline $\st$. An interesting detail is that this step uses of the skin uniformity of $H$ and $C$ another time. In this place it is exploited to understand the transition of shielding effects between nodes in the tree $\T$. \\

   Then, in Ch.\ref{smoo}, we exploit the controls over the $\st$ and $\sp$ to find that the shielding effects  combine additively, up to controllable defects.  Choosing the tree $\T$ appropriately, these defects can be made arbitrarily small and we get the asserted shielding effect along $\Sigma$.

\subsubsection{Setup and Notations} \label{sn}

Throughout this paper, $H^n \subset M^{n+1}$ will be a connected locally area minimizing hypersurface without boundary in some $C^\infty$-smooth Riemannian
$n+1$-manifold $(M,g_M)$. $g_H$ is the induced metric on the smooth portion of $H$.\\

In technical terms, $H$ is a locally mass minimizing, integer multiplicity rectifiable current of dimension $n$ without boundary.
The partial regularity theory for these minimizers says that $H$ is a smooth hypersurface except for some singular set $\Sigma_H$ of Hausdorff-dimension $\le n-7$.\\

We consider the following classes of complete area minimizing hypersurfaces.
\begin{description}
    \item[${\cal{H}}^c_n$]:    $H^n \subset M^{n+1}$ is a compact and connected hypersurface without boundary.
    \item[${\cal{H}}^{\R}_n$]:  $(M,g_M) = (\R^{n+1},g_{\R^{n+1}})$, $H$ is an oriented boundary of some open set $A \subset \R^{n+1}$ with $0 \in H$. In particular, $H$ is non-compact and complete.
     \item[${\cal{SH}}^{\R}_n$]:      ${\cal{SH}}^{\R}_n \subset  {\cal{H}}^{\R}_n$ is the subset of those hypersurfaces singular, at least, in $0$.\\
\end{description}
The main class of hypersurfaces we study in this paper is given by
\[{\cal{H}}_n:= {\cal{H}}^c_n \cup {\cal{H}}^{\R}_n \, \mm{ and }\, {\cal{H}} :=\bigcup_{n \ge 1} {\cal{H}}_n.\]
${\cal{H}}_n$ is closed under blow-ups. That is, the limit of converging subsequences under scaling by a diverging sequence of real numbers belongs to ${\cal{H}}^{\R}_n$.\\

In dimensional reduction arguments we also consider more particular classes of minimal hypersurfaces:
\begin{description}
    \item[$\mathcal{C}_{n}$]:   $\mathcal{C}_{n} \subset {\cal{H}}^{\R}_n$ is the space of area minimizing $n$-cones in $\R^{n+1}$ with tip in $0$.
    \item[$\mathcal{SC}_{n}$]:    $\mathcal{SC}_n \subset\mathcal{C}_n$ is the subset of cones singular, at least, in $0$
    \item[$\mathcal{K}_{n-1}$]: For any area minimizing cone $C \subset \R^{n+1}$ with tip $0$, we get the non-minimizing minimal hypersurface $S_C$ in the unit sphere
\[S_C:= \p B_1(0) \cap C \subset S^n \subset  \R^{n+1} \mm{ and set } {\cal{K}}_{n-1}:= \{ S_C\,| \, C \in {\mathcal{C}_{n}}\},\] where $\mathcal{C}_{n}$ is the space of area
minimizing $n$-cones in $\R^{n+1}$ with tip in $0$. We write ${\cal{K}}= \bigcup_{n \ge 1} {\cal{K}}_{n-1}$, for the space of all such hypersurfaces $S_C$.
\end{description}

$A = A_H$ denotes the second fundamental form of $H \subset M$, $|A|$ is the norm or length of $A$ and $\RH$ is the set of positive regular values of $|A|$. We
    denote $|A|$-level sets, for $c >0$, by  $\B_c := |A|^{-1}(c)  \subset H$.\\

For the practical use of a visual distinction, we typically write $\Sigma_H$ for the (compact) singular set of any $H \in {\cal{H}}^c_n$, but $\sigma_C$ for the (non-compact) singular set of a cone $C \in \mathcal{SC}_{n}$. Also we consider the \emph{one-point compactification} $\widehat{C}$ of $C$ respectively  $\widehat{\sigma_C}$ of $\sigma_C$ and denote the point at infinity by $\infty_C$.\\

Since we are interested in intrinsically defined concepts and problems, we use that integer multiplicity rectifiable currents of dimension $n$ without boundary admit local decomposition  into oriented boundaries of open sets in the ambient space, cf.[L1], appendix, part IV:\\

Any locally mass minimizing current $T \in \mathcal{R}_n(U)$, with $\p T =0$, in $U$ can be decomposed in $U$ into oriented boundaries $\p M_i$, for open sets $M_i \subset
U$, $i \in \Z$, $M_i \subset M_{i+1}$, each locally area minimizing in $U$
 so that for any open $W \subset \! \subset U$
\begin{equation}\label{lde}
 T\llcorner U= \sum_{i \in \Z} \p \llb M_i \rrb \llcorner U \,\mm{ and } \, {\bf{M}}_W (T\llcorner U)= \sum_{i \in \Z} {\bf{M}}_W (\p \llb M_i \rrb \llcorner U).
\end{equation}
 This allows us to consider most technical aspects in the framework of oriented boundaries. In chapters 2 and 3 we can therefore assume that the given $H$ actually is such a oriented boundary. Only in Ch.\ref{coverings}, where this is not self-explaining, we give some additional comments.

\bigskip
\setcounter{section}{2}
\renewcommand{\thesubsection}{\thesection}

\setcounter{section}{2}
\renewcommand{\thesubsection}{\thesection}
\subsection{Inherent Self-Similarities} \label{fm0}

\bigskip

Our main objectives in this chapter is to describe self-similar substructures along $\Sigma \subset H$, like self-similar ball covers of $\Sigma$. \\

They simplify the organization of the deformations we install close to $\Sigma$, reducing the amount of information, concerning $\Sigma$, we need to take into account. They help us to circumvent peculiarities of non-uniform and non-unique tangent cone approximations.

\subsubsection{Tangent Cone Approximations} \label{na}

 When we blow-up around any singular point $p \in \Sigma$, that is we scale $H$ by ever larger constants, we observe that, for $\tau \ra \infty$,
 $\tau \cdot H$ (sub)converges to some area minimizing cone, a tangent cone, which is a partially linearized and simplified local model for $H$.\\

However these approximations are more delicate than tangent space approximations of smooth manifolds. There could be infinitely many distinct tangent cones in a given $p \in \Sigma$. Also the space
of tangent cones in $p \in \Sigma$ and the accuracy of the cone approximation of $H$  may change discontinuously when we vary the base point $p \in \Sigma$.\\

Here we describe how we can recoup some uniform control from the use of compactness results from the the theory of area minimizers and elliptic regularity theory.  We use these results to make these cone approximations applicable also for the analysis of the conformal  Laplacian on $H$.\\

For starters, we recall a classical formulation, cf. [F],4.3.16, [Gi], 9.3 or [Si],37.4.

\begin{proposition} \emph{\textbf{(Tangent Cones)}} \label{flat-norm-approx}
\,  Let $H^n \subset M^{n+1}$ be an area minimizing hypersurface and $\tau_m \to +\infty$ a sequence of positive real numbers. Then, for every $p \in H$ we find a
subsequence $\tau_{m_k}$ and an area minimizing cone $C_p \subset \R^{n+1}$ such that for any given open $U \subset \R^{n+1}$ with compact closure the {\it flat norm}
$\db_U$   converges to zero:
$$
\db_U (\tau_{m_k} \cdot H, C_p) \to 0.
$$
When $\overline U$ contains only smooth points of $C_p$, the $\db_U$-convergence implies {\it compact $C^l$-convergence}, for any $l \ge 0$.
\end{proposition}

 \begin{remark} \textbf{(Iterated Tangent Cones)}\label{itc} \, When we  blow-up around a point $p_0 \in \Sigma_H$ we get a tangent cone $C$. Scaling around $0 \in \sigma_{C}$ merely reproduces $C$,  but blow-ups around singular points $p_1\neq 0 \in \sigma_{C_1}$ generate for some \emph{iterated tangent cones}:  area minimizing cones $C'$, where $C'$  is a Riemannian product $\R \times C^{n-1}$ for some area minimizing cone $C^{n-1} \subset \R^n$. This iteration ends with a cone $\R^m \times C^{n-m} $, for some $C^{n-m} \subset \R^{n-m}$ singular only in $0$. It is a basic fact that $n-m \ge 8$, since lower dimensional minimizing  cones are regular, cf.[Gi], Ch.11.
 \qed
 \end{remark}

The first step to improve this result uses the regularity theory within the space of tangent cones. \\

To this end we additionally introduce the notations ${\cal T}_H \subset \mathcal{SC}_n$  for the set of  tangent cones of points in $\Sigma \subset H$ and, specifically, ${\cal T}_p \subset
    \mathcal{SC}_n$ the set of  tangent cones in the point $p \in \Sigma \subset H$. Since we usually do not have any preferred directions on $M$, the sets ${\cal T}_H \subset \mathcal{SC}_n$ and ${\cal T}_p \subset \mathcal{SC}_n$ are homogenized: with
each $C$ we also include all images under linear isometries of $\R^{n+1}$ into the respective space.\\

 The same flat norm convergence results as in \ref{flat-norm-approx} and the definition of tangent cones give the following compactness results.

\begin{lemma} \emph{\textbf{(Compact Cone Spaces)}} \label{cone-compactness}
 \begin{enumerate}
\item The spaces $\mathcal{C}_{n}$ and ${\cal{H}}^{\R}_n$ are compact - under compact convergence in the flat norm topology. $\mathcal{K}_{n}$ is compact in flat norm topology.
\item There is a constant $d_n >0$ so that for $C \in \mathcal{C}_{n}$:\,
\[\db_{B_1(0)} (C, Y) < d_n, \mm{ for some hyperplane } Y \subset \R^{n+1} \Leftrightarrow C \mm{ is non-singular.}\]
\item  $\mathcal{SC}_n \subset\mathcal{C}_n$ and ${\cal{SH}}^{\R}_n \subset {\cal{H}}^{\R}_n$ are closed (and hence compact). Consequently $\overline {\cal T}_H \subset
    \mathcal{SC}_n$, also we have $\overline {\cal T}_p =  {\cal T}_p$.\\
\end{enumerate}
\end{lemma}

\begin{remark} \textbf{(Almost Minimizers)}\,  The minimal hypersurfaces $S_C=\p B_1(0) \cap C$, for any area minimizing cone $C$, are \emph{almost minimizers}. The $S_C$ are critical points of the area functional, but they are not area minimizing and they cannot even be stable, since $Ric_{S^n} >0$, cf.[L1],Appendix, Prop.5.2 for a broader discussion. But being an almost minimizer is sufficient to infer that the  $S_C \in \mathcal{K}_n$ share the regularity theory of area minimizers. In particular, all arguments for area minimizers used in this chapter including those involving the skin transforms equally apply to them.\\
\end{remark}

 \begin{remark} \textbf{(Notions of Compactness)}\label{noc} \, The general compactness results for currents refer to \emph{sequential compactness} relative the topology given by the family of flat pseudometrics, cf.[Si2],Ch.31. But, each $N \in \mathcal{K}_{n}$ has a unique representation as the compact metric completion of the manifold $N \setminus \Sigma_N$. Thus, on $\mathcal{K}_{n}$, we have a proper flat \emph{metric} and, hence, sequential compactness implies compactness i.e. each open cover of $\mathcal{K}_{n}$ contains a finite one.\\
 \end{remark}

\textbf{Freezing Effects near Singular Points} \,
Now we turn to the question how tangent cones can serve as smooth approximative local models near a singular point  $p \in \Sigma \subset H$. For this, we measure the
non-tangential approachability of $p$.

\begin{definition} \label{pen} \quad For $H \subset M$ and any $\rho >0$, we define the \textbf{regular pencil} $\P(z,\rho)$, pointing to $z \in \Sigma$
\begin{equation}\label{ii}
 \P(z,\rho)=\P_H(z,\rho):= \{x \in H \setminus \Sigma \,|\, \delta_{\bp} (x) > \rho \cdot d_{g_H}(x,z)\}.
\end{equation}
and the, generally singular, complement in $H$:
\[\P^*(z,\rho):= H \setminus \P(z,\rho)=\overline{\{x \in H \setminus \Sigma \,|\, \delta_{\bp} (x) \le \rho \cdot d_{g_H}(x,z)\}}^H.\] $arctan(\rho^{-1})$ is a measure for the size of the (solid) opening angle of $\P(z,\rho)$.\\

Also we define the notion of $\bp$-distance tubes
\begin{equation}\label{inn}
 \I(\rho)= \I(H,\rho):=\{x \in H \setminus \Sigma \,|\, \delta_{\bp}(x) < \rho \} \mm{ and its exterior complement }
\end{equation}
 \[\E(\rho):= H \setminus \I(\rho),\mm{ for any } \rho >0.\]
\end{definition}

 The following result formally describes the phenomenon that the twist of these pencils slows down, it literally freezes, when we approach $z$, and we get a fine non-tangential approximation through tangent cones.

\begin{proposition} \emph{\textbf{(Freezing Effect)}}  \label{freez} Let $p \in \Sigma_H \subset H$, some $\ve > 0$ and a pair $R \gg 1 \gg r >0$ be given.
 Then we have the following results
\begin{itemize}
    \item \emph{\textbf{Flat-norm Version}} \, There is some $\tau_{\ve , R , r,p} \gg 1$   such that for every $\tau \ge \tau_{\ve , R , r,p}$ there is a tangent cone
        $C_p^\tau$ of $H$ at $p$ with:
$$\tau \cdot (H \cap B_{R/\tau}(p) \setminus B_{r/\tau} (p)) \text{ is $\ve$-close in flat norm to } C^\tau_p\cap B_R(0) \setminus B_r (0 ).$$
    \item \emph{\textbf{$\mathbf{C^k}$-norm Version}} \,  For any additionally given $1 \gg \delta >0$ and $k \in \Z^{\ge 0}$, we can find a $\tau_{\ve , R , r, \delta,p,k} \ge
        \tau_{\ve , R , r,p}$ such that for every $\tau\ge \tau_{\ve , R , r, \delta,p,k}$
    \[\tau \cdot  B_{R/\tau}(z) \setminus B_{r/\tau} (z) \cap  \P(z,\delta) \subset \tau \cdot H\] can be written as a smooth section
   $\Gamma_\tau  \mm{ with } |\Gamma_\tau|_{C^k} < \ve$  of the normal bundle of
\[ B_R(0)\setminus B_r (0) \cap  \P(0,\delta) \subset C_p^\tau.\]
\end{itemize}
        \end{proposition}

 \begin{remark} \textbf{(Boundary Adjustments)}\label{ba} \,   Strictly speaking,  $\p \Gamma_\tau$ will only gradually approach $\p B_R(0) \cup \p
B_r (0)$, for $\tau \ra\infty$. However, it is only the portion, away from both $\p \Gamma_\tau$ and  $\p B_R(0) \cup \p B_r (0)$, we are interested in. Thus we may
always adjust the definitions near the boundary to our needs and henceforth omit this negligible detail.\\
 \end{remark}

{\bf Proof of \ref{freez}} \, The \textbf{flat norm version} largely follows from \ref{flat-norm-approx}. We assume  we could not find a $\tau_{\delta , R , r,p}$ as asserted, that is, there is
a sequence $\tau_i \ra \infty$, for $i \ra \infty$, so that for none of the $\tau_i$ there is a tangent cone $\delta$-close in flat norm to $\tau_i \cdot (H \cap B_{R/\tau}(p)
\setminus B_{r/\tau} (p))$.\\

 But \ref{flat-norm-approx} says, whatever the sequence $\tau_i \ra \infty$ looks like there is a subsequence
converging to some tangent cone $C$. Therefore,  \ref{flat-norm-approx} leads to a contradiction, since it gives us some $i_0 \gg 1$ so that $\tau_{i_0} \cdot (H \cap
B_{R/\tau}(p) \setminus B_{r/\tau} (p))$ is $\delta$-close to $C \cap B_R(0) \setminus B_r (0 ).$\\

For the \textbf{$C^k$-version} of the result we note that \ref{flat-norm-approx} implies the result for a single tangent cone. To get estimates
 not depending on the given cone we use the compactness of the cone spaces \ref{cone-compactness} and the naturality of $\bp$.\\

 Namely, assume again we could not find a $\tau_{\ve , R , r,\delta,p,k}$ as asserted. Then there is
a sequence $\tau_i \ra \infty$, for $i \ra \infty$, so that for none of the $\tau_i$ there is a tangent cone so that ${\tau_i} \cdot H$ can locally be written as a section
\[\Gamma_{\tau_i} \mm{ with } |\Gamma_{\tau_i}|_{C^k} < \delta\] of the normal bundle over $B_R(0)\setminus B_r (0) \cap \P(0,a) \subset C^{\tau_i}_p.$ But due to
the naturality of $\bp$ and the compactness of $\mathcal{SC}_n$, we reach the following contradiction:\\

 The sequence of truncated cones $B_R(0)\setminus B_r (0) \cap \P(0,a) \subset
C^{\tau_i}_p$ subconverges in flat norm, and from this in $C^k$-norm, to $B_R(0)\setminus B_r (0) \cap \P(0,a) \subset C$, for some tangent cone $C \in
\mathcal{SC}_n$. But for this cone, the assertion, e.g. for $\delta/2$, will be true. \qed

\subsubsection{Induced Entire Solutions and $\D$-Maps} \label{na2}
\bigskip

Here we show that the freezing effect extends to the behavior of solutions of naturally defined elliptic equations $L \, w=0$ near singular points. As an important detail we observe  that solutions on $H$ can asymptotically be approximated by solutions, not only locally but globally defined on tangent cones. When $L$ is skin adapted, this make these induced solutions amenable to a transparent potential theory.\\

\textbf{$\D$-maps}  \, We recall some terminology, from [L1],Ch.2.1, to concisely describe the local transfer of operators and functions defined on $H \setminus \Sigma$ to (the regular portion of) its tangent cones, back and forth.\\

We consider a concrete setup we encounter in compactness arguments: a sequence of pairs $H^n_i \subset M^{n+1}_i$, where the $M^{n+1}_i$ are complete Riemannian manifolds, with base point $p_i \in M_i$ and the $H^n_i$ are area minimizing hypersurfaces, with $p_i \in H_i$.\\

Now assume that the  $M_i$ are compactly $C^k$-converging, for some $k \ge 5$ to a limit manifold $M$, so that $p_i \ra p \in M$, for $i \ra \infty$, so that we have diffeomorphisms
\[\Psi_i : B_R(p_i) \ra B_R(p), \mm{ so that } |(\Psi_i)_*(g_{M_i}) - g_{M}|_{C^k} \ra 0 \mm{ on } B_R(p).\]
 for any given radius $R >0$ and all sufficiently large $i \ge i_R$.\\

Next we assume that $\p H_i \cap B_R(p_i) \v$. Then, basic compactness results show that the $H_i$
subconverge  in flat norm to a limit area minimizer $H \subset M$, with $p \in H$:
\[\Psi_i (B_R(p_i) \cap H_i) \ra B_R(p) \cap H, \mm{ for } i \ra \infty.\]

The two cases we consider later on are both related to tangent cones:
\begin{description}
  \item[\textbf{(BP)}] \textbf{Blow-up around a given point}\, \emph{$M=\R^{n+1}$, and $M_i:= \tau_i  \cdot M_0$ and $H_i:=\tau_i \cdot H_0$, for some sequence $\tau_i \ra \infty$, for $i \ra \infty$.  Then $H$ becomes one of the tangent cones $C$ of $H$ in $p$. Recall that tangent cones are usually not uniquely determined.}
  \item[\textbf{(CC)}] \textbf{Cone Convergence}\,  \emph{$M=M_i=\R^{n+1}$. Then we consider area minimizing cones $H_i, H \in {\cal{C}}_n$ and  compact convergence or equivalently global convergence of (sub)sequences of $S_{C_i}\in \mathcal{K}_{n-1}$}.\\
\end{description}

Next, if $B_R(p) \cap H$ is \emph{smooth}, this rather weak flat norm convergence implies a  smooth convergence:  we consider the normal bundle of $B_R(p) \cap H$ by $\nu$. Then, for $i$ large enough, the $\Psi_i (B_R(p_i) \cap H_i)$ become local and smooth sections $\Gamma_i$ of $\nu$.\\

 The regularity theory makes these local sections canonical identifying diffeomorphisms
and the flat norm convergence implies a $C^k$-convergence of the sections $\Gamma_i$ to the zero section which equals $B_R(p) \cap H$. We obviously can extend this notion to arbitrary open sets  $D \subset H$ with compact closures in $H \setminus \Sigma_{H}$ and their counterpart $D_i \subset H_i$ . We introduce the following mnemonic to name these canonically defined sections of the normal bundle.

\begin{definition} \label{idmap} \emph{\textbf{($\D$\textbf{-map}) }}\,  For sufficiently, large $i$, we call the uniquely determined section  \[\D_{H \ra H_i} := \Gamma_i : D   \ra  \Psi_i(D_i)\]  of $\nu$ over
$D \cap H$  the asymptotic  \textbf{identification map}, briefly the $\D$\textbf{-map},  written $\D$.\\
\end{definition}

In this terminology,  we  restate the $C^k$-version of \ref{freez}: for  $\tau\ge \tau_{\ve , R , r, \delta,p,k}$, we have
\[\D=\D^\tau : B_R(0)\setminus B_r (0) \cap  \P(0,\delta) \ra \tau \cdot  B_{R/\tau}(z) \setminus B_{r/\tau} (z) \cap  \P_{\tau \cdot M}(z,\delta)\]
\[\mm{satisfies } |\D^\tau - id_{C^\tau_p}|_{C^k(B_R(0)\setminus B_r (0) \cap  \P_{\R^{n+1}}(0,\delta))} \le \ve,\] for the norm on sections of the normal bundle and where $id_{C^\tau_p}$ is just the zero section over ${C^\tau_p}$.  In a given context, the $\D$-map is unambiguously defined. \\

\textbf{Freezing of Entire Solutions} \quad On these hypersurfaces, we consider elliptic problems with the naturality property that their solutions transfer to solutions of similar problems under blow-ups.
More generally, we assume there is an assignment $F$ of operators $H \mapsto F(H)$, for any $H \in {\cal{H}}$, so that $F(H)$ commutes with
convergence of the underlying spaces. Then we call the assignment, or just the operator, natural. The conformal Laplacian is our basic example, cf.[L2], Ch.5.2 for details.\\

Also, we typically consider natural and \emph{skin adapted} operators, cf.Ch.\ref{mre}C. and [L2], Ch.3.4 for details, to have a well-controlled potential theory, including simple integral representations for all positive solutions.\\

Now we will see that the freezing effect extends to the analysis of such operators. In the $\D$-map setup, we equip each of the $H_i$ with some positive solution $u_i >0$ of  $L(H_i)\, w=0$ on $H_i \setminus \Sigma_{H_i}$ and consider the two cases (BP) and (CC) above.
\begin{description}
  \item[\textbf{(BP)}] Due to naturality of $L$, $u_i$ also solves $L(\tau_i \cdot H_0)\, w=0$. From a diagonal  sequence argument we may choose an exhausting  sequence of smoothly bounded domains so that a suitably normalized* subsequence compactly $C^{2,\alpha}$-converges to a solution $u >0$ of $L(C)\, w=0$ defined on $C \setminus \sigma_C$.
  \item[\textbf{(CC)}]  Similarly, a suitably normalized*  subsequence of the $u_i$  compactly $C^{2,\alpha}$-converges to a solution $u >0$ of $L(H)\, w=0$ defined on $H \setminus \Sigma_{H}$.
\end{description}
(*Since $L$ is a linear operator we may multiply each $u_i$ by some suitable positive constant $c_i >0$ to ensure a local convergence to a finite and positive limit. Formally, in \ref{freezfu} and \ref{freezfu2}, we just consider quotients and compare them to $1$.)\\

To justify the assertions,  notice that the $\D$-maps eventually provide an almost isometric relation between arbitrarily large domains with compact closure in $H\setminus
\Sigma_{H}$ and the corresponding domains on $H_i \setminus \Sigma_{H_i}$.\\

Then, on these domains, the coefficients of $L(H_i )$ are arbitrarily close to the $\D$-map pull-back of the coefficients of $L(H)$ in H\"older norm. Thus
$u \circ \D$ solves elliptic equations with coefficients arbitrary close, in H\"older-norm, to those of $L(H)$. \\

Since weak solutions of $L  \, w = 0$ are $C^{2,\alpha}$-regular,  we
get \emph{uniform} constants for the elliptic regularity estimates and Harnack inequalities for positive solutions of these families of equations, cf.[GT]. This gives us compactly converging subsequences of these locally defined almost-solutions $u_i \circ \D$ of $L(H) \, w = 0$, for $i \ra \infty$.\\

Summarizing, any such limit is a positive $C^{2,\alpha}$-regular solving  of $L(H) \, w = 0$ is an \textbf{entire solution}, that is, it is defined on all of $H \setminus \Sigma_{H}$ and we call it a solution \textbf{induced} from the $u_i$ on  $H_i \setminus \Sigma_{H_i}$. And we observe a freezing of the coupling between entire solutions on $H$ and the $u_i$, in the sense of \ref{freez}: for a given $\ve >0$ and \emph{any} $i \gg 1$ and for any given open subset $D$ with compact closure $\overline{D} \subset  H \setminus \Sigma_H$, there is \emph{some} entire solution $u>0$ of  $L(H) \, w = 0$ with
\[   |u_i  \circ \D /u -1|_{C^{2,\alpha}(D)}  \le \ve.\]
 The reason is, as in \ref{freez}, that otherwise we could find a subsequence $u_{i_k} \circ \D$ that not even compactly subconverges to an entire solution.\\

  More precisely, we get the following upgrade of \ref{freez} in the (BP)-case. The counterpart for the (CC)-case is formulated in \ref{freezfu2} below.

\begin{corollary} \emph{\textbf{(BP-Freezing Effects)}}  \label{freezfu} Let $p \in \Sigma_H \subset H$,  $\ve, \delta >0$ and $R \gg 1 \gg r >0$ be given. Then, there is a large $\tau^*(L,\ve,\delta, R , r, p)>0$ so that:\\

\noindent For any real $\tau\ge \tau^*$ and any entire  solution $v >0$ of $L(H) \, w = 0$, there is a tangent cone $C^\tau_p$, with the following properties

\begin{itemize}

\item  \emph{\textbf{Spatial Freezing}}  \, $C^\tau_p\cap B_R(0) \setminus B_r (0 )$ is $\ve$-close in flat norm to $\tau \cdot (H \cap B_{R/\tau}(p) \setminus B_{r/\tau} (p))$ and for the $\D$-map from $C^\tau_p$ to $\tau \cdot H$:
\[\D_{C^\tau_p \ra \tau \cdot H} : B_R(0)\setminus B_r (0) \cap  \P(0,\delta) \ra \tau \cdot  B_{R/ \tau}(p) \setminus B_{r/\tau} (p) \cap  \P_{\tau \cdot H}(p,\delta) \mm{ we have }\]
\begin{equation}\label{ffefu}
 |\D - id_{C^\tau_p}|_{C^5(B_R(0)\setminus B_r (0) \cap  \P(0,\delta) )} \le \ve.
\end{equation}
\item \emph{\textbf{Functional Freezing}}  \, On $B_R(0)\setminus B_r (0) \cap  \P(0,\delta)$ we have for some \textbf{entire solution} $u>0$  of $L(C^\tau_p) \, w = 0$ defined on $C^\tau_p \setminus \sigma_{C^\tau_p}$:
\begin{equation}\label{ffefua}
 |u \circ \D / v-1|_{C^{2,\alpha}(B_R(0)\setminus B_r (0) \cap  \P(0,\delta) )} \le \ve,
\end{equation}
for given $\alpha \in (0,1)$.
\end{itemize}
\end{corollary}

Henceforth we choose $\tau^*_{L,\ve , \delta, R , r, p}$ as the \emph{smallest} value, so that \ref{freezfu} holds and think of its inverse as the \emph{start radius} for balls centered in $p \in \Sigma_H$ from that on \emph{all} balls of smaller radius are non-tangentially approximated by tangent cones and so that after rescaling to some unit size, the
non-tangential accuracy expressed in  (\ref{ffefu}) and (\ref{ffefua}) is smaller than $\ve$.\\

We formalize this observation. For our notational convenience,  we choose some custom coupling of the a priori independent size parameters to a single accuracy parameter.

\begin{definition} \emph{\textbf{(Approximation Accuracy)}} \label{fi} \, For  $H \in {\cal{H}}_n$ and $L$ given, we call the $\ve$ parameter of  \ref{freezfu}, the
\textbf{accuracy}  of the $\D$-map respectively of the approximation and  \[\mathcal{R}_\ve(L,p):= 1/\tau^*_{L,\ve, \ve^2,\frac{1}{\sqrt{\ve}}, \sqrt{\ve}, p} \mm{ is the \textbf{$\ve$-accuracy radius} } \mm{ in } p \in \Sigma_{H}.\]
An approximation  of high accuracy means one with a very small $\ve$.
\end{definition}

More explicitly, for any ball $B_\rho(p) \subset H$ around $p$ with radius $\rho \le\mathcal{R}_\ve(p)$  the  (BP)-freezing holds with accuracy $\ve$ once we scaled
 it to a ball of radius $R$: After scaling by $1/(\sqrt{\ve} \cdot \rho)$ we have for $B_R(0)\setminus B_r (0) \cap  \P(0,\ve^2)$, with $R =\frac{1}{\sqrt{\ve}}$ and $r=\sqrt{\ve}$ that the approximation quality, also near the far ends $\p B_R$ (and  $\p B_r$), is $\sqrt{\ve}=\ve/\sqrt{\ve}$.\\

The freezing effect \ref{freezfu} can be enhanced to a stable version, valid also for merely \emph{almost tangent} cones, cf.\ref{bow1} below. For this, we also use $\D$-maps between hypersurfaces $H,H' \in {\cal{K}}_n$ and the corresponding functional freezing result for the (CC)-case:

\begin{corollary} \emph{\textbf{(CC-Freezing Effects)}}  \label{freezfu2} Let  $\ve, \delta >0$ and $R \gg 1 \gg r >0$ be given. Then there is an $\zeta(L,\ve,\delta , R , r)>0$ so that for $C \in {\cal{C}}_n$, any $C' \in {\cal{C}}_n$ with $d_H(S_C,S_{C'}) \le \zeta$ and any entire solution $u_{C'} >0$ of $L(C') \, w = 0$ on $C'\setminus \sigma_{C'}$, we have an entire   solution $u_{C,C'}>0$ of $L(C) \, w = 0$ on $C\setminus \sigma_{C}$ so that
\begin{equation}\label{ffefu}
 |\D - id_{C}|_{C^5(B_R(0)\setminus B_r (0) \cap  \P(0,\delta) )} \le \ve, \mm{ and }
\end{equation}
\begin{equation}\label{ffefu2}
 |u_{C'}  \circ \D /u_{C,C'} -1|_{C^{2,\alpha}(B_R(0)\setminus B_r (0) \cap  \P(0,\delta) )}  \le \ve,
\end{equation}
for given $\alpha \in (0,1)$.
\end{corollary}

{\bf Proof} \quad  This follows the same lines of contradiction arguments as before.  The independence of $\ve$ from $C$ uses the compactness of the space ${\cal{C}}_n$. \qed

\subsubsection{Self-Similar Covers of $\Sigma$ } \label{coverings}

\bigskip

We describe self-similar ball covers of $\Sigma \subset H$ and observe a self-regulating alignment of the singular directions of these balls.
We start with some preliminary considerations.\\

\textbf{Reduction to Oriented Boundaries}  \,  To define these ball covers, we consider $\Sigma$ and $H$ as subsets of $M$ to exploit its locally nearly Euclidean geometry of  $M$ in some combinatorial arguments used to control the covering numbers.\\

To this end we review our setup. Any of the currents $H \in {\cal{H}}_n$ admits a local decomposition into oriented boundaries of an open set in $M$, formally stated as (\ref{lde}). This can be accomplished through a cover of  $M$ by sufficiently small open balls $B[1],... B[k]$ so that in each of these balls we have such a decomposition into sheets of oriented boundaries. When $M$ is compact,  there are only finitely many such sheets since they all contribute a lower bounded amount of area (better mass) to the bounded total area. Also  for $M=\R^{n+1}$, $H$ already is an oriented boundary, by definition.\\

When we now choose some ball $B_r(p) \subset B[i] \subset M$ around some point $p \in \Sigma$ of radius $r>0$ we observe that $B_r(p)$ usually intersects not only the sheet $p$ belongs to, but potentially also some others and this means balls from different sheets may intersect.  This problem can easily be solved: the radii $r$ are  chosen much smaller than the distance between any two sheets in any of the $B[i]$, cf. also [L1], appendix IV. (We observe that higher multiplicity sheets do not matter in this context.)\qed

\textbf{Intrinsic versus Extrinsic Radii} \,  With this reduction, all radii of balls in $H^n$,  which appear in \ref{add} and  \ref{brop} below, are understood as being measured in terms of the extrinsic metric $g_M$ of the ambient space $M^{n+1}$. However, for our further applications, from Ch.\ref{sish} onwards, the intrinsically measured length becomes more appropriate. This does not cause any issues since the freezing result \ref{freez} and the skin uniformity of $H \setminus \Sigma$, cf. [L1], Ch.4, show that the two concepts asymptotically coincide:\\

 For $H \in {\cal{H}}^{\R}_n$, we recall from  [L1], Lemma 4.8: For any $H \in
{\cal{H}}^{\R}_n$ and any two points $p,q \in H \subset \R^{n+1}$ we have
 \begin{equation}\label{e1e0}
 d_{g_H}(p,q)  \le  c^{\R}_n \cdot  d_{g_{\R^{n+1}}}(p,q) \le c^{\R}_n \cdot d_{g_H}(p,q).
 \end{equation}
 for some constant $c^{\R}_n \ge 1$ depending only on the dimension $n$.\\

 For a cone $C \subset \R^{n+1}$ with tip $0$,  $p \in C$ we have $d_{g_{\R^{n+1}}}(p,0) =  d_{g_C}(p,0)$. For general $H \in {\cal{H}}^{\R}_n$, this latter equality, (\ref{e1e0}) and \ref{freez} show that we get a uniform convergence
 \begin{equation}\label{e2e}
 |d_{g_{\tau \cdot H}}(x,z)-d_{g_{\tau \cdot M}}(x,z)| \ra 0, \mm{ for } \tau \ra \infty,  \mm{ and } x \in  \tau \cdot (H \cap B_{R/\tau}(z) \setminus B_{r/\tau} (z)).
 \end{equation}
 This carries over from $H \in {\cal{H}}^{\R}_n$ to any $H \in {\cal{H}}_n$, using the  reduction to oriented boundaries explained earlier in this section.  Indeed we only apply this argument in a range of distances where we have quantitative constraints on the accuracy of the $\D$-maps, cf.\ref{add} and the setup of \ref{brop}. \qed

Now we explore the idea that $H$ is coarsely self-similar when we zoom into $\Sigma$ or while we travel along this singular set. We write $\Sigma$ as a disjoint union of subsets sorted according to the $\ve$-accuracy radii of the respective balls.\\

In \ref{add} and  \ref{brop} all radii are understood as being measured in terms of the extrinsic metric $g_M$

\begin{definition} \emph{\textbf{(Accuracy Decompositions  of $\Sigma$)}}\label{add} \, For $H \in {\cal{H}}_n$ or $H \in {\cal{K}}_n$ and $L$ a given natural skin adapted operator, we define the \textbf{accuracy decomposition} of $\Sigma = \dot{\bigcup}_{i\ge 1}\Sigma_i,$, more precisely $\Sigma_H = \dot{\bigcup}_{i\ge 1} (\Sigma_H)_i$, for a given accuracy $\ve >0$, as follows:
\[ \Sigma_1 = \{ p \in \Sigma \,|\, \mathcal{R}_\ve(L, p) >  10^3 \} \mm{ and } \Sigma_i = \{ p \in \Sigma \,|\, 10^3 \cdot  s^{i-2} \ge \mathcal{R}_\ve(L, p) >  10^3 \cdot s^{i-1}\},\]
for $i \ge 2$, where $s \in (0,10^{-3})$ is a given parameter.
\end{definition}

 \begin{remark}  \label{oo} \, 1. The decomposition largely depends on $H \setminus \Sigma$ near $\Sigma$ and not only the point set $\Sigma$ alone. It may happen that for all but finitely many $i$: $\Sigma_i \v$.\\

  Also we observe for any $H \in {\cal{K}}_n$ and any general* $H \in {\cal{H}}^c_n$, that for any $i_0$ there is a small $\ve(i_0) >0$ so that
$\Sigma_i \v$ for any $i \le i_0$.\\

 *In turn, in the case of a singular cone $C \in {\cal{C}}_n$ we always have $0 \in \Sigma_1$. Thus in the rather exotic and clearly non-analytic case where
$H \in {\cal{H}}^c_n$ is locally exactly isometric to such a cone,  we also had $0 \in \Sigma_i$, for some fixed $i(H) \ge 1$ and any $\ve>0$. But this can be generically be excluded and it never happens for $H \in {\cal{K}}_n$. \\

2. In general, the $\Sigma_i$ will be uncountable but each of them contains a countable dense
subset $\Sigma^{*}_i \subset \Sigma_i$ and we set $\Sigma^{*}:=\bigcup_{i \ge 1}\Sigma^{*}_i \subset \Sigma$. We write the elements of $\Sigma^*_i$ accordingly as
$a_{i,1}, a_{i,2},..$.  We use the $\Sigma^*_i$ as an index set to argue via ordinary induction, in place of a less elementary transfinite induction.\\

Specifically, we assign the ball $\overline{B_{s^{i}}(p)} \subset M^n$ to any of the points $p \in \Sigma^*_i$ and set
\begin{itemize}
  \item  ${\cal{B}}_i:= \{\overline{B_{s^{i}}(p)}\,|\, p \in \Sigma^*_i\}$ and ${\cal{B}}:=\dot{\bigcup}_{i\ge 1}{\cal{B}}_i$
  \item  Then ${\cal{B}}_i$ is a cover of $\Sigma_i$ and, hence,  ${\cal{B}}$ is a cover of $\Sigma$.  Namely, given any point $x  \in \Sigma_i$  there is $p_x  \in \Sigma^*_i$ with $d_M(x,p_x) < s^{i}/2$ and thus $x \in \overline{B_{s^{i}}(p_x)}$.
\end{itemize} \qed\end{remark}

Now we describe self-similar subcovers ${\cal{A}}_i = {\cal{A}}_i \subset {\cal{B}}_i$, for each $\Sigma_i$.

\begin{proposition} \emph{\textbf{(Self-Similar Covers)}} \label{brop} \, Let $H \in {\cal{H}}_n$ or $H \in\mathcal{K}_{n}$  and $L$ be given. Then, for any \textbf{accuracy} $\ve>0$ and some $s_0(M) \in (0,10^{-3})$ respectively $s_0(n) \in (0,10^{-3})$, we have: For any
\textbf{self-similarity factor} $s \in (0,s_0]$ the following holds:
\begin{enumerate}
\item For each $i \ge 1$, there is a locally finite family ${\cal{A}}[i](\Sigma_i,\ve,s)$ of  closed balls
\[{\cal{A}}[i] = \{\overline{B_{s^i}(p)} \subset H\,|\, p \in A[i]\}, \mm{ for some discrete set } A[i] \subset \Sigma^*_i.\]
\item For any $k \ge 1$ we have: $\bigcup_{i \le k} {\cal{A}}[i]$ covers $\bigcup_{i \le k} \Sigma_i$.
\item For  $H \in {\cal{H}}^c_n$  or $H \in\mathcal{K}_{n}$ each of the families ${\cal{A}}[i]$ is finite, but there can still be infinitely many $i$ with ${\cal{A}}[i]\n$
 \item  $q \notin \overline{B_{s^i}(p)}$, $p \notin \overline{B_{s^k}(q)}$, for any two different $\overline{B_{s^i}(p)} \in {\cal{A}}[i]$, $\overline{B_{s^k}(q)} \in {\cal{A}}[k]$.
 \item There is a constant $c(n)$, depending only on the dimension, so that ${\cal{A}}[i]$ splits into  disjoint families ${\cal{A}}[i,1],...,{\cal{A}}[i,c]$ with
 \begin{equation}\label{10}
 B_{10 \cdot s^i}(p)
     \cap B_{10 \cdot s^i}(q) \v, \mm{ for } \overline{B_{s^i}(p)}, \overline{B_{s^i}(q)}, \mm{ in the same family }{\cal{A}}[i,j]
 \end{equation}
\end{enumerate}
The  ${\cal{A}}[i]$  are called the \textbf{layers} of the cover ${\cal{A}}(\Sigma_H,\ve,s):=\bigcup_{i \ge 1}{\cal{A}}[i]((\Sigma_H)_i,\ve,s)$ of $\Sigma$ by balls around the points of $A:=\bigcup_{i \ge 1}A[i]=\bigcup_{i \ge 1, c \ge j \ge 1}A[i,j]$. \\

The ${\cal{A}}[i,1],...,{\cal{A}}[i,c]$ are the \textbf{sublayers} of
the layer ${\cal{A}}[i]$. Finally, in the case where $\Sigma \v$, we set ${\cal{A}} :\v$ and $A:\v$.
\end{proposition}

 \textbf{Proof} \quad We inductively define a selection map $\kappa$:
\[\kappa=(\kappa_1,\kappa_2):
\Sigma^{*} \ra\Z^{\ge 0} \times \Z^{> 0} \mm{ from its restrictions } \kappa[i]: \Sigma^{*}_i \ra\Z^{\ge 0} \times \{i\}.\]  The value $\kappa[i]_1$ says that, when $\kappa[i]_1(p)=0$, the ball
$\overline{B_{s^{i}}(p)}$ is deleted and otherwise it is placed in a family of balls labelled by $\kappa[i](p)=(\kappa[i]_1(p),i)$.\\

The $\kappa[i]$ are defined from two nested inductions: over $i$, passing from $\Sigma^*_i$ to $\Sigma^*_{i+1}$ and, secondly, over $j$, for the elements $a_{i,j} \in
 \mathcal A[i]$, for $i$ fixed.\\

When $\Sigma^{*}_i \v$, there is no map $\kappa[i]$ to define and we continue with  $\Sigma^{*}_{i+1}$ and set $A[i] \v, {\cal{A}}[i] \v$. In particular, to start  the induction over $i$ we define $\kappa_1$ on the \emph{first non-empty} $\Sigma^*_i$:\\

\textbf{Start on $\Sigma^*_i$}: \,For $a_{i,1} \in \Sigma^*_i$, we set $\kappa[i]_1(a_{i,1}):=1$.\\

 \textbf{Step on $\Sigma^*_i$}: \,We assume $\kappa_1$ has been defined for $a_{i,j}$, $j \le m$. Then we set
\[\kappa[i]_1(a_{i,m+1}):=0,\, \mm{ if }\, a_{i,m+1} \in \textstyle \bigcup_{j \le m} \overline{B_{s^i}(a_{i,j})},\,\mm{ and otherwise,}\]
 \[\kappa[i]_1(a_{i,m+1}):= \min(\{j \le m\,|\,d(a_{i,j},a_{i,m+1}) >  10 \cdot s^i + 10 \cdot s^i\} \cup \{m+1\}).\]\

Now, in the second induction layer, we assume  $\kappa_a$ on $\Sigma^*_a$ has been defined for any $a \le i$. Then the next step is definition of $\kappa_{i+1}$ on $\Sigma^*_{i+1}$:\\

 \textbf{Start on $\Sigma^*_{i+1}$}: \,We start with  $a_{i+1,1} \in \Sigma^*_{i+1}$  and set $\kappa[i+1]_1(a_{i+1,1}):=0$, when
\[a_{i+1,1} \in  \textstyle \bigcup_{a \le i} \{\overline{B_{s^a}(a_{a,j})} \,|\, a_{a,j} \in \Sigma^*_a,\kappa[a]_1(a_{a,j}) \ge 1 \}\,\mm{ and } \kappa[i+1]_1(a_{i+1,1}):=1  \mm{ otherwise.}\]

 \textbf{Step on $\Sigma^*_{i+1}$}: \,We assume $\kappa[i+1]$ has been defined for $a_{i+1,j}$, $j \le m$. Then we set
$\kappa[i+1]_1(a_{i+1,m+1}):=0,$ when \[a_{i+1,m+1} \in \textstyle \bigcup_{j \le m} \overline{B_{s^{i+1}}(a_{i+1,j})} \cup \bigcup_{a \le i} \{\overline{B_{s^a}(a_{a,j})} \,|\,
a_{a,j} \in \Sigma^*_a,\kappa[a]_1(a_{a,j}) \ge 1 \}, \] and otherwise,
 \[ \kappa[i+1]_1(a_{i+1,m+1}):= \min(\{j \le m\,|\,d(a_{i+1,j},a_{i+1,m+1}) >  10 \cdot s^{i+1} + 10 \cdot s^{i+1}\} \cup \{m+1\}).\]

This completes the definition of $\kappa$ on $\Sigma^{*} $. Now we set  for $i,j \ge 1$:
 \[{\cal{A}}[i,j] := \{\overline{B_{s^i}(p)}\,|\, p \in \Sigma^*_i, \kappa[i]_1(p)=j\},\, {\cal{A}} := \bigcup_{i\ge 1} {\cal{A}}[i]:= \bigcup_{i,j\ge 1} {\cal{A}}[i,j]\]
 and the corresponding sets of center points \[A[i,j]:= \{p \in \Sigma^*_i\,|\, \kappa[i]_1(p)=j\},\, A := \bigcup_{i\ge 1} A[i]:= \bigcup_{i,j\ge 1} A[i,j].\]

Next we choose $s_0 >0$ small enough so that $s_0^{-1} \cdot M^{n+1}$ locally looks nearly like the flat $\R^{n+1}$ and the following holds for any $z \in s_0^{-1} \cdot M$:
\begin{equation}\ \label{vol} k_1 \le Vol(B_{1/3}(z))\mm{ and } Vol(B_{30}(z))  \le k_2,\end{equation}
for some constants $k_i(n) >0,  i=1,2$. When $M$ is compact, such an $s_0$ certainly exists, and when $M$ is non-compact it already is a Euclidean space.\\

We claim that for $s \in (0,s_0]$ and any $i \ge 1$: \begin{equation}\label{defc}c(n):= (\mm{\emph{the smallest integer}}  \ge k_2/k_1) \, \Rightarrow \, {\cal{A}}[i,j] \v,\mm{ for } j
> c.\end{equation}

 Otherwise, we had some $\overline{B_{s^i}(p)} \in {\cal{A}}[i,c+1]$ with $B_{10 \cdot s^i}(p) \cap B_{10 \cdot s^i}(x_m) \n$
for at least $c$ different $x_m \in A[i,j]$ and we get the following estimates, when we scale by $1/s^i  \ge 1/s_0$.  In $1/s^i \cdot M$ the $(c+1)$ balls $B_{1/3}(x_m)$ and $B_{1/3}(p)$ are pairwise disjoint and satisfy
\[(c+1) \cdot k_1 \le \sum Vol(B_{1/3}(x_m))+ Vol(B_{1/3}(p))=\]
\[Vol \big(\bigcup B_{1/3}(x_m) \cup B_{1/3}(p)\big) \le Vol(B_{30}(p)) \le k_2.\]

But this contradicts the definition of $c$, and, hence, ${\cal{A}}[i,c+1] \v$. \\

In particular, any ball of radius $\le 2 \cdot s^i$ around some point of $H$, intersects at most $c$ balls $\overline{B_{s^i}(p_k)}$, $p_k \in {\cal{A}}[i]$. From this local finiteness, we now observe that, for any $k \ge 1$, $\bigcup_{i \le k} {\cal{A}}[i]$ covers $\bigcup_{i \le k} \Sigma_i$.\\

Namely, if we had some  $q \in  \Sigma_k$ with
 $q \in \Sigma_k \setminus   \bigcup_{p \in \bigcup_{i \le k} A_i} \overline{B_{s^i}(p)}$, we first notice from the local finiteness of each ${\cal{A}}[i]$ that $\bigcup_{p \in \bigcup_{i \le k} A_i} \overline{B_{s^i}(p)}$ is a closed subset of $H$. Thus $q $ belongs to the open complement.\\

Then, since $\Sigma^*_k \subset \Sigma_k$ is dense, there is also some $q'= a_{k,m}\in  \Sigma^*_k \setminus \bigcup_{p \in \bigcup_{i \le k} A_i} \overline{B_{s^i}(p)}$. But such an $a_{k,m}$ cannot exist: for $\kappa[i]_1(a_{k,m})=0$, $a_{k,m}\in \bigcup_{p \in  \bigcup_{i \le k} A_i} \overline{B_{s^i}(p)}$, a contradiction. Otherwise, $\kappa[i]_1(a_{k,m})=j>0$ and thus $a_{k,m} \in A[k,j]$ and this again shows $a_{k,m}\in \overline{B_{s^i}(a_{k,m})} \subset  \bigcup_{p \in \bigcup_{i \le k} A_i} \overline{B_{s^i}(p)}$. Thus there was also no point $q \in  \Sigma_k \setminus \bigcup_{p \in\bigcup_{i \le k} A_i} \overline{B_{s^i}(p)}$.\\

Finally, due to the compactness of $H \in {\cal{H}}^c_n$  or $H \in\mathcal{K}_{n}$, the disjointness of the balls $\overline{B_{s^i/3}(p)}$, for any $\overline{B_{s^i}(p)}$ in ${\cal{A}}[i]$, and their uniformly lower bounded volume contribution $v_i>0$, expressed in (\ref{vol}),  each of the covers ${\cal{A}}[i]$ is finite. Thus one could even find more explicit bounds for their cardinality and the growth rate for increasing $i$. (But this does neither mean that ${\cal{A}}$ is finite nor it admits a finite subcover.)\\

The remaining assertions follow readily from the definitions. In particular, we note that for  $H \in\mathcal{K}_{n}$,  $s_0$ depends only on the dimension. \qed

\subsubsection{Auto-Alignments of Covers} \label{sish}
\bigskip

We generally think of both parameters, the accuracy $\ve>0$ and the self-similarity factor $s>0$, in the construction of ${\cal{A}}$, as being extremely small. Then we observe that the singular directions of the balls in the cover ${\cal{A}}$ of \ref{brop} show some self-regulating alignment properties:
\begin{itemize}
  \item  For any two adjacent balls $\overline{B_{s^i}(p)}$ and $\overline{B_{s^i}(q)}$ in one layer ${\cal{A}}_i$, the singular and highly curved directions in any two tangent
  cones $C_p$  and $C_q$, in $p$ and $q$, align like compass needles in a magnetic field.
  \item  For pairs of adjacent balls in different layers the smaller ball in some  ${\cal{A}}_{i_2}$ is always located very closely to the singular and highly curved
directions emanating from the center of the larger ball  in ${\cal{A}}_{i_1}$, for ${i_2}>{i_1}$.
\end{itemize}

We describe this in terms of pencils $\P$ directly on $H$. Recall that   $|\delta_{\bp}(p)- \delta_{\bp}(q)|   \le L_{\bp} \cdot d(p,q)$, for $p,q \in  H \setminus \Sigma$.
For notational simplicity we assume that $L_{\bp} = 1$.  Also, without further notification, we henceforth use the equivalences between intrinsic and extrinsic distances discussed in the beginning of Ch.\ref{coverings}, (\ref{e1e0}) and (\ref{e2e}), above, to view the ball covers as being intrinsically defined.

\begin{proposition}\emph{\textbf{(Auto-Alignments)}}\label{comp1}\, For the  self-similar covers ${\cal{A}}$ of \ref{brop} we have
\begin{itemize}
  \item \emph{\textbf{Alignment within ${\cal{A}}_i$}} \, For any two $p \neq q \in A[i]$ with $d_{g_H}(p,q) \le 10 \cdot s^i$ and any  $\kappa, \mu \in (0,1)$:
  \begin{equation}\label{al}
 B_{5 \cdot s^i}(q)\cap \left(\P^*(p,\kappa \cdot \mu/20) \cup B_{ \kappa \cdot s^i/10}(p)\right) \subset \P^*(q,\kappa)\cup B_{\mu \cdot s^i}(q).
\end{equation}
  \item \emph{\textbf{Alignment of ${\cal{A}}_{i_1}$ relative ${\cal{A}}_{i_2}$}} \, For any $\kappa \in (0,1)$, there is an $s[\kappa, M]\in (0,s_0(M))$, so that for $s \in (0,s[\kappa, M])$,  $i_2 >i_1$ and $p \in A[i_2]$  and $q \in A[i_1]$:
\begin{equation}\label{al1}
 B_{5 \cdot s^{{i_2}}}(p)  \subset   \P^*(q,\kappa).
\end{equation}
\end{itemize}
\end{proposition}

 \textbf{Proof} \, We start with the alignment relation within ${\cal{A}}_i$ and distinguish between the two cases $x \in B_{5 \cdot s^i}(q)\cap  \P^*(p,\kappa \cdot \mu/20)$ and $x \in B_{5 \cdot s^i}(q)\cap B_{ \kappa \cdot s^i/10}(p)$.\\

We first note that, for $x \in \P^*(p,\kappa \cdot \mu/20) \cap B_{5 \cdot s^i}(q)$ we have $\delta_{\bp}(x) \le 15 \cdot \kappa \cdot \mu\cdot s^i/20<\kappa \cdot \mu\cdot s^i$, since  $d_{g_H}(p,q) \le 10 \cdot s^i$ and $d_{g_H}(x,q) \le 5 \cdot s^i$.\\

For $x \in  B_{5 \cdot s^i}(q) \setminus B_{\mu  \cdot  s^i}(q)$ this implies $\delta_{\bp}(x) \le   \kappa  \cdot d_{g_H}(q,x)$. The remaining option  is that $x \in B_{\mu \cdot  s^i}(q)$ and thus we conclude
 \begin{equation}\label{a0}
\P^*(p, \kappa \cdot \mu/20) \cap B_{5 \cdot s^i}(q) \subset \P^*(q,\kappa)\cup B_{\mu \cdot s^i}(q).
\end{equation}
Secondly, for $x \in B_{5 \cdot s^i}(q)\cap B_{ \kappa \cdot s^i/10}(p)$ we have from $\delta_{\bp} \le dist(\cdot,\Sigma)$
\[\delta_{\bp}(x) \le  \kappa  \cdot s^i/10 \mm{ and } d_{g_H}(p,x) \le  \kappa  \cdot s^i/10 < s^i/10.\]
Measured from $q$, this gives us $s^i/2\le  d_{g_H}(q,x)$, since $d_{g_H}(p,q) > s^i$. Therefore, we have $\delta_{\bp}(x) \le  \kappa  \cdot s^i/10 <   \kappa \cdot d_{g_H}(q,x)$ and this shows $B_{5 \cdot s^i}(q)\cap B_{\kappa \cdot s^i/10}(p) \subset   \P^*(q, \kappa)$. We combine this with (\ref{a0}) to get (\ref{al}).\\

For the alignment relation of ${\cal{A}}_{i_1}$ relative ${\cal{A}}_{i_2}$, we note from $\delta_{\bp} \le dist(\cdot,\Sigma)$ that $\delta_{\bp}(x) \le 5 \cdot s^{i_2}$, for $x \in B_{5 \cdot s^{i_2}}(p)$. Moreover, from ${i_2} \ge {i_1}+1>1$, $s < 10^{-3}$ and $d_{g_H}(p,q) > s^{i_1}$ we have:
\[ s^{i_1}/2 \le s^{i_1}- 5 \cdot s^{i_2}  \le  d_{g_H}(x,q), \mm{ for any } x \in B_{5 \cdot s^{i_2}}(p).\]
Thus,  for $s \in (0,s[\kappa])$ with $s[\kappa]:=\min\{\kappa/10,s_0\}\le 10^{-3}$, we have for $x  \in B_{5 \cdot s^{i_2}}(p)$:
 \[\delta_{\bp}(x) \le 5 \cdot s^{i_2} \le  5 \cdot s^{i_1+1}\le 10 \cdot s \cdot d_{g_H}(x,q) \le \kappa \cdot d_{g_H}(x,q).\]
In other words, $B_{5 \cdot s^{i_2}}(p) \subset \P^*(q,\kappa)$, for  $s \in (0, s[\kappa])$.  \qed

\subsubsection{Separators and $\bp$-Pipelines} \label{sish1}
\bigskip

The alignment relations \ref{comp1} suggest to consider difference sets of nested pairs of sets
 \begin{equation}\label{q}
 \Q[p,R,\alpha,r]:= B_R(p)\cap  \left(\P^*(p, \alpha)\cup B_r(p)\right), \mm{ for }  R >r >0, \alpha  >0.
 \end{equation}
For each layer ${\cal{A}}_{i}$ of the self-similar covers ${\cal{A}}$ we have $c(n)$ sublayers. Now we define more specific  pairs of radii and solid angles, we later assign sublayerwise,  to build  controllably disjoint families of geometric objects $\s$ locally \emph{surrounding} $\Sigma_H$, cf. Ch.\ref{smoo}.\\

The following concept of separators introduces control parameters for these difference sets, the $\bp$-pipes we discuss and exploit below. The size of separators is reflected in uniform estimates for the deformations we later place into $\bp$-pipes, and their inductive refinements the $\bp$-multi-pipes.

 \begin{definition} \emph{\textbf{(Separator Sequences)}} \label{sse} \, We choose a finite collection $\Theta=\Theta_n$ of real numbers $\pi_n[a,b] \in (0, \frac{1}{100})$, $a=1, ,...,c(n)$, $b= 0,1$ with
 \[\pi_n[1,0]>   \pi_n[1,1] > ....>\pi_n[c(n),0] >  \pi_n[c(n),1]\, \mm{ and }\]
 \begin{equation}\label{bd} \frac{\pi_n[a,0]}{10} > \pi_n[a,1] \mm{ and }
\frac{\pi_n[a,1]^2}{100} > \pi_n[a+1,0]
 \end{equation}
We call such a sequence of numbers $\Theta_n$, an $n$-dimensional \textbf{separator}  (sequence).
\end{definition}

In inductive arguments, we also consider tuples $\widehat{\Theta}=(\Theta_n,.. ,\Theta_7)$ of separators and we call $\Theta_m$ a \emph{small separator} when $\pi_m[1,0]\ll1$ and we call it a \emph{broad separator} when $\pi_n[a,1]/\pi_n[a,0]\ll 1$ and $\pi_n[a+1,0]/\pi_n[a,1]\ll 1$, for all $a$, where the meaning of being $\ll 1$ depends on the context.

\begin{definition} \emph{\textbf{($\bp$-Pipes)}} \label{sp} \, For a self-similar cover ${\cal{A}}$, a separator $\Theta$ and any $p \in A[i,j] \subset A \subset \Sigma_H$, $i \in \Z^{\ge1}$ and $1 \le j \le c(n)$, we set
\begin{equation}\label{s}
\s[\Theta, s, p]= \s[\pi_n[j,0],\pi_n[j,1], s, i, j, p]:=  ...
\end{equation}
\[ \Q[p,5 \cdot s^i,\pi_n[j,0] ,\pi_n[j,0]\cdot s^i] \setminus \Q[p,5 \cdot s^i,\pi_n[j,1] ,\pi_n[j,1]\cdot s^i]=...\]\vspace{-0.3cm}
\[B_{5 \cdot s^i}(p)\cap  \left(\P^*(p, \pi_n[j,0])\cup B_{\pi_n[j,0] \cdot s^i}(p)\right) \setminus  B_{5 \cdot s^i}(p)\cap  \left(\P^*(p, \pi_n[j,1])\cup B_{\pi_n[j,1] \cdot s^i}(p)\right).\] \\
We call these difference sets $\s \subset H \setminus \Sigma$ the \textbf{$\bp$-pipes}. \\
\end{definition}

The pipes $\s[\Theta, s, p]$  inherit and amplify the auto-alignment properties \ref{comp1} of $\cal{A}$:

\begin{proposition} \emph{\textbf{(Auto-Aligned Pipes)}}\label{comm}\, For any separator $\Theta$ and $H \in \mathcal{H}^c_{n}$, for some ambient manifold $M^{n+1}$, there is a $s[\pi_n[c(n),1], M] >0$, so that for a self-similar cover ${\cal{A}}$ with self-similarity factor $s \in (0,s[\pi_n[c(n),1], M])$, we have for any two $p,q \in A$,  $p \neq q$:  \[ \s[\Theta, s, p] \, \cap \,  \s[\Theta, s, q]\v.\]
For $H=S_C  \in\mathcal{K}_{n}$, the value $s[\Theta] := s[\pi_n[c(n),1], S^{n+1}]$ depends only on  $\Theta$.
\end{proposition}
($\Theta$ determines the two actual control parameters $\pi_n[c(n),1]$ and $n$.)\\

\textbf{Proof} \, For any two $p,q \in A[i]$, \ref{brop}(v) condition (\ref{10}) shows that we only need to care about the case where $q \in A[i,a+m]$ and $p \in A[i,a+k]$, for some
$k > m \ge 1$.\\

In this case, $\frac{\pi_n[a+m,1]^2}{100} > \pi_n[a+k,0]$ of \ref{sse}(\ref{bd}) shows  from \ref{comp1}(\ref{al})   that
\[ B_{5 \cdot s^i}(q)\cap  \s[\Theta, s, p] \subset B_{5 \cdot s^i}(q)\cap \left(\P^*(p,\pi_n[a+k,0]) \cup B_{\pi_n[a+k,0]\cdot s^i}(p)\right)... \]
\[...\subset \P^*(q,\pi_n[a+m,1])\cup B_{\pi_n[a+m,1] \cdot s^i}(q), \mm{ and thus }\s[\Theta, s, p]  \cap \s[\Theta, s, q] \v.\]\

For $p \in A[{i_2}]$, $q \in A[{i_1}]$, ${i_2} >{i_1}$, with $d_{g_H}(p,q) \le 10 \cdot s^{i_1}$ and $s \in (0,s[\pi_n[c(n),1],M])$, we have from (\ref{al1}): \[\s[\Theta, s, p]\subset B_{5 \cdot s^{i_2}}(p) \subset \P^*(q, s[\pi_n[c(n),1]]),\]
and hence $\s[\Theta, s, p] \, \cap \, \s[\Theta, s, q] \v$. \qed\\

\textbf{Finite Open Subcovers} \, We want to get a better control over the alignments of $\bp$-pipes to transfer $\bp$-pipelines between different spaces via $\D$-maps. For this we use compactness results, like that  for the singular set but also of the space of minimal cones ${\cal{C}}_n$. Eventually, this brings us to the introduction of multi-pipes, we reach in the next chapter.\\

As a first step, we note that self-similar covers by closed balls of \ref{brop}, are generally infinite and do not admit finite subcovers. However,  for doubled radii  of  the balls in the cover $\cal{A}$, the interiors of the new balls form a cover of $\Sigma_H$. Then, for $H \in {\cal{H}}^c_n$ or $H \in\mathcal{K}_{n}$, the compactness of $\Sigma_H$ allows us to select a finite subcover.\\

More concretely, there is a constant $i^\cs(H,\ve,s) \ge 1$ so that
\begin{enumerate}
\item ${\cal{A}}^\cs:=\bigcup_{i^\cs \ge i \ge 1}{\cal{A}}^\cs[i]$, for  ${\cal{A}}^\cs[i]:= \{B_{2 \cdot s^i}(p)\,|\, p \in A[i]\}$, already covers $\Sigma_H$.
\item There is a constant $c(n)$, depending only on the dimension, so that each layer ${\cal{A}}^\cs[i]$ splits into  disjoint sublayers ${\cal{A}}^\cs[i,1],...,{\cal{A}}^\cs[i,c]$ with \[B_{10 \cdot s^i}(p)
     \cap B_{10 \cdot s^i}(q) \v, \mm{ for } B_{2 \cdot s^i}(p), B_{2 \cdot s^i}(q), \mm{ in the same sublayer }{\cal{A}}^\cs[i,j].\]
\end{enumerate}

We still enlarge the balls in ${\cal{A}}^\cs$ and define a cover ${\cal{A}}^*(\Sigma_H,\ve,s)$ of $\Sigma_H$ as:
 \begin{equation}\label{ae}
 {\cal{A}}^*:=\bigcup_{i^\cs \ge i \ge 1}{\cal{A}}^*[i], \mm{ for } {\cal{A}}^*[i]:=\{B_{5 \cdot s^i}(p)\,|\, p \in A[i]\}   \mm{ and  write }  A^* :=\bigcup_{i^\cs \ge i \ge 1}A[i].
\end{equation}
 As a formal side note, there are obviously many possible choices to select ${\cal{A}}$ and ${\cal{A}}^*$ upon a given choice of $H,\ve,s$. To rule out ambiguities in subsequent constructions which are based on the chosen covers, we always choose one such cover and keep it fixed. Similarly, in all further constructions we choose one fixed sample. Thereby these constructions formally depend only on $\Sigma \subset H,\ve,s$.\\

A notable advantage of working with this finite covers is that, different from the full infinite union, $\dot\bigcup_{p \in A^*}\s[\Theta, s, p]$ remains disjoint from small metric distance tubes  and $\bp$-distance tubes of $\Sigma_H$:

 \begin{lemma} \emph{\textbf{(Basic $\bp$-Pipelines)}} \label{ab}\,  For $H \in\mathcal{K}_{n}$ and some separator $\Theta$,  there  are some minimal $\rho[0]({\cal{A}}^*(\Sigma_H),\Theta)  >0$ and a maximal  $\rho[1]({\cal{A}}^*(\Sigma_H),\Theta) >0$, so that
  \begin{equation}\label{tu0}
  \dot\bigcup_{p \in A^*(\Sigma_H)} \s[\Theta, s, p] \subset  \overline{U_{\rho[0]/L_{\bp}}(\Sigma_H)} \subset \overline{\I(H,\rho[0])}
 \end{equation}
 \begin{equation}\label{tu}
 U_{\rho[1]/L_{\bp}}(\Sigma_H) \cap \dot\bigcup_{p \in A^*(\Sigma_H)}\s[\Theta, s, p]  =  \I(H,\rho[1]) \cap \dot\bigcup_{p \in A^*(\Sigma_H)}\s[\Theta, s, p] \v.
 \end{equation}
 \begin{equation}\label{ttu}
 U_{\rho[1]/L_{\bp}}(\Sigma_H)  \subset \I(H,\rho[1]) \subset \bigcup_{p \in A^*(\Sigma_H)} B_{2 \cdot s^i}(p)
 \end{equation}
 We call  $\rho[0]$ the \textbf{outer width} and $\rho[1]$ the \textbf{inner width} of the \textbf{$\bp$-pipeline} \[\sp  :=\dot\bigcup_{p \in A^*(\Sigma_H)}\s[\Theta, s, p].\]
\end{lemma}

 \textbf{Proof} \quad For (\ref{tu0}) we choose $\rho[0]/L_{\bp}$ as the maximum of the radii of all balls in $ {\cal{A}}^*$.\\

Similarly, for  (\ref{tu}), we note that for each $\s[\Theta, s, p]$, $p \in A^*$, we have some $\rho_p>0$ so that $\I(\rho_p) \cap \s[\Theta, s, p] \v$. Since $\cs A^* < \infty$, we can just choose $\rho=\min_{p \in A^*} \rho_p >0$ and find $ \I(H,\rho) \cap \dot\bigcup_{p \in A^*(\Sigma_H)}\s[\Theta, s, p] \v$. Since, $\delta_{\bp} \le L_{\bp}  \cdot dist(x,\Sigma)$, we have $U_{\rho/L_{\bp}}(\Sigma_H) \subset  \I(H,\rho)$ and this implies $U_{\rho/L_{\bp}}(\Sigma_H) \cap \dot\bigcup_{p \in A^*(\Sigma_H)}\s[\Theta, s, p] \v$.\\

The last property (\ref{ttu}) for the inner width can be accomplished employing the properness of $\bp$ on compact hypersurfaces.\qed

The relations  (\ref{tu})  and  (\ref{ttu})  will later be paraphrased saying that $\Sigma$, and even $ \I(H,\rho[1])$, are surrounded by the $\bp$-pipeline $\sp$.  For the cone $C^{n+1}$ associated to $H=S_C$,  these relations mean that the cone $C(\sp) \subset C^{n+1}$ over  $\sp$ satisfies   \begin{equation}\label{vv}
C(\sp) \subset \P^*(0, \rho[0]) \setminus  \P^*(0, \rho[1]) \subset C^{n+1},
\end{equation}
that is,  the pipe concept is compatible with dimensional reductions.

\begin{corollary}   \emph{\textbf{(Truncated Cones within Pipes)}} \label{tcw} \,
For any $S_C \in \mathcal{K}_{n}$ and an $\bp$-pipeline $\sp \subset S_C$, the \textbf{truncated cone} over $\sp$ \[C_{tr}(\sp,i,s):= C(\sp)  \cap (B_{5 \cdot s^i}(0) \setminus B_{\rho[1] \cdot s^i}(0)) \subset  C,\] for some $i \in \Z^{>0}$,  $s \in (0,s_0]$
is contained in a \textbf{carrier pipe} $\s \subset C$:
\begin{equation}\label{ttr}
C_{tr}(\sp) \subset   \s[\pi_{n+1}[j,0],\pi_{n+1}[j,1], s, i, j, 0],
\end{equation}
for a separator $\Theta_{n+1}$ with $\pi_{n+1}[j,0]=\rho[0]$ and $\pi_{n+1}[j,1]=\rho[1]$
\end{corollary}
\bigskip

\textbf{Custom Coupling to Cone Approximations} \, The relations for $\bp$-pipes we have seen so far are motivated and formulated to match the self-similar covers ${\cal{A}}$ of \ref{brop}, but a main property of  ${\cal{A}}$ we have not yet used is the presence of tangent cone approximations of prescribed accuracy $\ve$. Indeed the size of the accuracy parameter $\ve$ affects the choice of ${\cal{A}}$, but it has no influence on the assertions or the chosen parameters $s$ and $\Theta$ in \ref{comp1} - \ref{ab}.\\

This is a recurrent detail. Intuitively speaking, the $\bp$-pipes are independent geometric beings, following their own rules and capable to live under any accuracy regime.
Only when we deliberately want to transfer $\bp$-pipes between different spaces we introduce a custom coupling between the separator $\Theta$ and $\ve$.
Then we observe an obvious but vital compatibility of $\bp$-pipes under $\D$-maps:

\begin{lemma} \emph{\textbf{($\D$-Transfer of $\bp$-Pipes)}} \label{it}\, For any separator $\Theta$, there is an $\ve_0(\Theta)>0$, so that for a self-similar cover $\cal{A}$ of accuracy $\ve \in (0,\ve_0)$ and a suitable tangent cone $C$ in $p \in A[i,j] \subset A \subset \Sigma_H$:
\begin{itemize}
\item \, $\s[\Theta, s, p]\subset s^{-i} \cdot B_{s^i/\sqrt{10 \cdot \ve}}(p)\setminus B_{s^i \cdot \sqrt{10 \cdot \ve}} (p) \cap  \P(0,100 \cdot\ve^2)$
\item \, $\Delta(\D(\s[\Theta, s, 0]), \s[\Theta, s, p]) \subset U_{3 \cdot \ve}(\p (\s[\Theta, s, p]))$,
\end{itemize}
 where $\Delta(\cdot,\cdot)$ is the difference set and $U_{3 \ve}(\cdot)$ the distance neighborhood.
\end{lemma}

 \textbf{Proof} \quad This readily follows from the definition of the accuracy \ref{fi}, the freezing effects and the Lipschitz continuity of $\delta_{\bp}$.\qed

This also explains, why we choose the pipes $\s[\pi_n[j,0],\pi_n[j,0], s, i, j, p]$  to surround pieces of $\Sigma$, in place of the full tube $\Q[p,5 \cdot s^i,\pi_n[j,0] ,\pi_n[j,0]\cdot s^i]$, since $\bp$-pipes are supported away from $\Sigma$. Again, this compatibility completely fails for pipes similarly defined using the metric distance to the singular set, in place of the $\bp$-distance.

  \begin{remark}  \textbf{($\D(\s)$ versus $\s$)} \label{itr}\, \ref{it} allows us  to simplify the terminology: we henceforth also call the $\D$-image of a $\bp$-pipe, for sufficiently high accuracy that depends on the context, again a $\bp$-pipe. Indeed, in our applications, we are interested in subsets of the $\s$ supported away from $\p(\s)$. Therefore,
 the small deviation of the $\D$-image $\D(\s)$ of a pipe $\s \subset C$   to the  intrinsically defined one $\s \subset H$  is irrelevant. This also applies to compositions of several $\D$-maps. \qed
\end{remark}

\setcounter{section}{3}
\renewcommand{\thesubsection}{\thesection}
\subsection{Stable Trees and Multi-Pipes on $H \setminus \Sigma$} \label{fm1}

\bigskip

We use the pipeline widths of \ref{ab} to inductively build pipes with surrounding $\Sigma$, called $\bp$-multi-pipes. They have an internal substructure of nested pipes and cones over pipelines, we iteratively pull-back from tangent cones around singular points.\\

 Different from the definition of $\bp$-pipes, in \ref{sp}, their definition and control can only be achieved very closely to $\Sigma$. This makes it suggestive to think of $\bp$-multi-pipes as nanostructures on $H \setminus \Sigma$.\\

A hurdle in the construction of $\bp$-multi-pipes is that the width of $\bp$-pipelines depends  discontinuously on variations of $H \in\mathcal{K}_{n}$. Namely, the singular set $\Sigma_H$, but also the consecutive selections of ${\cal{A}}$ and of ${\cal{A}}^*$, generally change abruptly under even slight variations of $H$ and we cannot expect common positive lower bounds for $\rho[i]({\cal{A}}^*(\Sigma_H),\Theta)$ on $\mathcal{K}_{n}$ or common upper bounds for the ratio $\rho[0]/\rho[1]$.\\

To resolve this issue we consider \emph{shared} covers and pipelines $\dot\bigcup_{p \in A^*(\Sigma_H)}\s[\Theta, s, p]$ for a given $H$, with other $H' \in\mathcal{K}_{n}$, sufficiently close to $H$.  These shared pipelines, also share common localized bounds for the widths. \\

The technical twist is to widen the range of tangent cones through \emph{almost tangent cones}. They are still area minimizers, but not properly tangential. Adjoining them allows us to employ the compactness of $\mathcal{K}_{n}$ to globalize these width estimates.\\

\subsubsection{Stabilizations and Almost Tangent Cones} \label{gr}
\bigskip

For any single ball $B_{5 \cdot s^i}(p) \in {\cal{A}}^*(H)$ for some $H \in\mathcal{K}_{n}$ we have a controlled (BP)-freezing approximation \ref{freezfu} of this ball via
tangent cones in $p \in \Sigma$. From the (CC)-freezing \ref{freezfu2} we know that for any other $H'= S_{C'} \in\mathcal{K}_{n}$ with $\db(S_C,S_{C'})\ll 1$ small enough,
we have fine approximations of $H= S_{C}$ near $p$.\\

  Since the cover ${\cal{A}}^*(H)$ is finite this applies uniformly to all balls in ${\cal{A}}^*(H)$ and we observe from \ref{freezfu} and \ref{freezfu2} that the tangent cone approximations of $H$, for given radii of balls, can also be transferred to $H'$. Explicitly, this reads as follows.

\begin{lemma}   \emph{\textbf{(Stable Covers)}} \label{bow0} \,  For any accuracy $\ve >0$, $s \in (0,s_0(n))$ and any $H \in\mathcal{K}_{n}$, there is a small $\eta_H(\ve,s)$,  so that the cover ${\cal{A}}^*$ can be shared with any $H' \in \mathcal{K}_{n}$, $\eta_H$-close to $H$ in the following sense:\\

For any $H' \in \mathcal{K}_{n}$ with $\db(H,H') <\eta_H$ and any $p \in {\cal{A}}^*[i] \subset {\cal{A}}^*$ there is a tangent cone $C$ of $H$ in $p$, so that the following composition $\D^{cp}$ of the $\D$-maps from subsets of $C$ to $s^{-i} \cdot H$ and $s^{-i} \cdot H$ to $s^{-i} \cdot H'$ has accuracy $\ve$:
\[\D_{C \ra s^{-i} \cdot H'}^{cp}: B_\frac{1}{\sqrt{\ve}}(0)\setminus B_{\sqrt{\ve}}(0) \cap  \P(0,\ve^2) \xrightarrow{\D_{C \ra s^{-i} \cdot H}}  s^{-i} \cdot   B_\frac{s^i}{\sqrt{\ve}}(p) \setminus B_{s^i \cdot \sqrt{\ve}} (p) \cap  \P(p,\ve^2)\]\vspace{-0.6cm}
\begin{equation}\label{ii}...\xrightarrow{\D_{s^{-i} \cdot H \ra s^{-i} \cdot H'}} \D_{C \ra s^{-i} \cdot H'}^{cp}(B_\frac{1}{\sqrt{\ve}}(0)\setminus B_{\sqrt{\ve}}(0) \cap  \P(0,\ve^2)) \subset s^{-i} \cdot H'.\end{equation}\vspace{-0.2cm}
\end{lemma}

Now we use stable covers and the compactness of ${\cal{K}}_n$ to get finite subcovers of $\{B_{\eta_H}(H) \,|\, H \in {\cal{K}}_n\}$ of ${\cal{K}}_n$ to resolve the issue with the missing uniform bound for inner width $\rho$ in (\ref{tu}) through the use of a discrete grid of covers.

\begin{proposition}   \emph{\textbf{(Stable Grids)}} \label{bow1} \,   For any accuracy $\ve >0$, a separator $\Theta_n$, $s \in (0,s[\Theta])$ and $H \in \mathcal{K}_{n}$, we choose the finite cover
${\cal{A}}^*(\Sigma_H,\ve,s)$ of $\Sigma_H$ described in (\ref{ae}) and a \textbf{stabilization accuracy} $\ve^*\in (0,\ve)$.\\

Then there is a finite subset $\mathcal{L}_n(\Theta_n,\ve,\ve^*,s) \subset \mathcal{K}_{n}$, we call a \textbf{stable grid}, such that for any $H' \in\mathcal{K}_{n}$ there is an $H \in \mathcal{L}$ with $\db(H,H') <\Delta(H)$, for some  density function $\Delta(H,\Theta,\ve,\ve^*,s) >0$ on $\mathcal{K}_{n}$, small enough so that:
\begin{itemize}
\item There are common widths $\rho^*[i](H,\Theta,\ve,\ve^*,s) >0$, $i=0,1$, for these $H' \in\mathcal{K}_{n}$, chosen minimal for $i=0$ and maximal for $i=1$, with the following properties: \[\D :  \E(H, \rho^*[1]/10^3) \ra   \E(H',\rho^*[1]/10^3) \mm{ is well-defined with accuracy  } \ve^* \mm{ and }\]
    \begin{equation}\label{tt0}
 \D_{H \ra H'}\left(\dot\bigcup_{p \in A^*(\Sigma_H)}\s[\Theta, s, p]\right) \subset \overline{\I(H',\rho^*[0])}
\end{equation}
\begin{equation}\label{tt}
\I(H',\rho^*[1]) \cap  \D_{H \ra H'}\left(\dot\bigcup_{p \in A^*(\Sigma_H)}\s[\Theta, s, p]\right)\v.
\end{equation}
\item For any $B_{5 \cdot s^i}(p) \in {\cal{A}}^*[i](H) \subset {\cal{A}}^*(\Sigma_H,\ve,s)$ and a suitable tangent cone of $H$ in $p$, the accuracy of  the composed $\D$-maps $\D^{cp}$ (\ref{ii})  into $s^{-i} \cdot H'$ is  $ 2\cdot \ve$.
     \item For any $H \in \mathcal{L}$, we have
\begin{equation}\label{pr}
 U_{\rho^*[1]/L_{\bp}}(\Sigma_H)  \subset \I(H,\rho^*[1]) \subset \bigcup_{p \in A^*(\Sigma_H)} B_{2 \cdot s^i}(p).
\end{equation}
\end{itemize}
We call $\rho[{\mathcal{L}}_n,0]:=\max_{H \in {\mathcal{L}}_n}\rho^*[0](H)<\infty$ the \textbf{stable outer width} and $\rho[{\mathcal{L}}_n,1]:=\min_{H \in {\mathcal{L}}_n}\rho^*[1](H)>0$ the \textbf{stable inner width}  of the stable grid ${\mathcal{L}}_n$.
\end{proposition}
  \begin{remark}    \label{re}\, 1. We observe from \ref{oo}.1. that for any given separator $\Theta$ and compact $H \in {\cal{H}}_n$ or $H \in\mathcal{K}_{n}$:
\begin{equation}\label{li}
 \rho[{\mathcal{L}}_n,0] \ra 0,  \, \rho[{\mathcal{L}}_n,1] \ra 0 ,  \mm{ for } \ve \ra 0,
\end{equation}
In turn, independent of the chosen $\ve>0$, we have from the definition already of a single of $\bp$-pipes \ref{sp}
\begin{equation}\label{li2}
\rho[{\mathcal{L}}_n,0]/\rho[{\mathcal{L}}_n,1] \ra \infty,  \mm{ for } \pi_n[a,0]/\pi_n[a,1]\ra \infty,
\end{equation}
for $a=1, ,...,c(n)$.\\

2. The statements make sense also for $\Sigma_H \v$. Then $\ve,s$ are formal parameters for an empty cover, but we get the non-trivial constraint $\ve^*\in (0,\ve)$ for the stable grid.\\

3.   As in \ref{it} and \ref{itr} we observe that $\D$-images of a basic pipeline for some $H$ in the stable grid  $\D\left(\dot\bigcup_{p \in A^*(\Sigma_H)}\s[\Theta, s, p]\right)$ can be used as a $\bp$-pipeline in $H'$, and to keep the notations as simple as possible we also denote such $\D$-images by $\sp$.\\

4.  Since the $\D$-map does not run from a tangent cone of $H'$ to $H'$, but from $H$ to $H'$, we observe that $\D$-images of $\sp$ need not to be close to  unions of $\bp$-pipes in $H'$. For instance, when $\Sigma_{H'}\v$  we had no genuine $\bp$-pipe in $H'$ at all.
\qed
\end{remark}

 \textbf{Proof} \,  We know from \ref{ab} that for any  $H \in\mathcal{K}_{n}$ there is a $\rho[1]_H >0$, so that
 \begin{equation}\label{tu1}\I(H,\rho[1]_H) \cap \dot\bigcup_{p \in A^*(\Sigma_H)}\s[\Theta, s, p] \v.
 \end{equation}
 From this, \ref{bow0} shows that for any $\ve^\bullet>0$ there is an $\eta_H(\ve^\bullet)>0$ small enough so that for any $H^\bullet \in\mathcal{K}_{n}$ with
 $\db(H,H^\bullet) <\eta_H$ the map $\D :  \I(H, \rho[1]_H/10^3)^c \ra   \I(H^\bullet,\rho[1]_H/10^3)^c$ is well-defined, cf.\ref{ba}, with accuracy  $\ve^\bullet$ and
\begin{equation}\label{tt}
\I(H^\bullet, \rho[1]_H/10) \cap  \D\left(\dot\bigcup_{p \in A^*(\Sigma_H)}\s[\Theta, s, p]\right)\v,
\end{equation}
and,  for $\ve^\bullet < \ve^n$, the chain rule also shows that for $B_{5 \cdot s^i}(p) \in {\cal{A}}^*[i](H) \subset {\cal{A}}^*(\Sigma_H,\ve,s)$ and some tangent cone of $H$ in $p$, we may assume that the accuracy of $\D^{cp}$ into $s^{-i} \cdot H'$ is  $2 \cdot \ve$.
Now we cover ${\cal{K}}_n$ by the open balls $B_{\eta_H(\ve^\bullet)}(H)$, for $H \in {\cal{K}}_n$. From the compactness of ${\cal{K}}_n$, we get a finite subcover $B_{\eta_{H_1}(\ve^\bullet)}(H_1),...,B_{\eta_{H_i}(\ve^\bullet)}(H_i)$. We choose $\mathcal{L}:=\{H_1,...,H_i\}$ and $\rho^*[1]:=\min_{m=1,..,i} \rho[1]_{H_m}/10$. The argument for the outer width is completely similar. \qed

Since the whole setup is designed to study natural operators we also want to understand how skin adaptedness is inherited to non-tangent cones from the stable grid.\\

In our sample case, the conformal Laplacian $L(H):= - \triangle   + \frac{n-2}{4 (n-1)}  \cdot scal_H$,  we do not need to care about this details: there is a common lower bound $a_n > 0$, so that $L(C)$ is  $a_n$-skin adapted for any $C \in \mathcal{SC}_{n}$. This follows from skin adaptedness of $L(C)$ for any single $C$ in [L2],Th.10,  since the ambient space is flat. Then we get a uniform positive bound for all these cones, from the compactness of $\mathcal{K}_{n}$. \\

For more general natural operators we have

\begin{corollary}   \emph{\textbf{(Functionally Stable Grids)}} \label{bow2} \,
For a natural operator $L$, $a$-skin adapted on $X  \in \mathcal{H}^c_n$,  an $n-1$-dimensional separator $\Theta_{n-1}$, some $a>0$ and $\ve>0$ we have a stable grid $\mathcal{L}_{n-1}(L,\Theta_{n-1},\ve,\ve^*,s) \subset \mathcal{K}_{n-1}$  so that in addition to the conditions and conclusions of \ref{bow1}:\\

For any tangent cone $C'$ of $X$, $H'=S_{C'}$  and, for any $H =S_C\in \mathcal{L}$ so that $\db(H,H') < \Delta$, for some density function $\Delta(L, H,\Theta_{n-1}, \ve,\ve^*,s)>0$, we have:
\begin{itemize}
  \item $L(C)$ is $a/2$-skin adapted.
  \item For any entire solution $u_{C'} >0$ of $L(C') \, v = 0$ on $C'\setminus \sigma_{C'}$, we have an entire solution
  $u_{C,C'}>0$ of $L(C) \, v = 0$ on $C\setminus \sigma_{C}$ so that
\begin{equation}\label{fsg}
 |u_{C'}  \circ \D /u_{C,C'} -1|_{C^{2,\alpha}(B_\frac{1}{\sqrt{\ve}}(0)\setminus B_{\sqrt{\ve}}(0) \cap  \P(0,\ve^2))}  \le \ve,
\end{equation}
for any given $\alpha \in (0,1)$.
\end{itemize}
\end{corollary}
\textbf{Proof} \,  We only need to care about the additionally asserted $a/2$-skin adaptedness of $L(C)$. For the second assertion we use \ref{freezfu2}.\\

For a natural operator $L$, $a$-skin adapted on $X  \in \mathcal{H}^c_n$, $a>0$, the inheritance of skin adaptedness [L2],Lemma 5.11 implies that $L(C)$ on any of its tangent cones $C$ is again $a$-skin adapted.  \\

There is an $\eta^+(L) >0$ so that, if there is an $H'=S_{C'}$,  $C' \in \mathcal{C}_n$ with $\db(H,H') <\eta^+$ and $L(C')$ is $a$-skin adapted, that is,
  there is a supersolution $s >0$ of  $L \, s \ge a \cdot \bp^2 \cdot s$, for some $a >0$ on  $C'$, then $L(C)$ is $a/2$-skin adapted.
 This also follows as in [L2],Lemma 5.8 from comparison of the variational integrals for the first eigenvalue and from the compactness of $\mathcal{K}_{n-1}$ we can choose  $\eta^+(L) >0$ independently of $C$. \qed

Typically, we fix one natural skin adapted operator $L$ in our applications. \ref{bow2} shows that we may omit mentioning $L$ but assume that the stable grid and $\Delta$
of \ref{bow1} have already been chosen tight enough to also serve in \ref{bow2}. This also applies to the following notion, we get and henceforth use the validity of the
 corresponding  functional version without stating it explicitly

\begin{definition}\emph{\textbf{(Almost Tangent Cones)}} \label{atc} \,  For  $X  \in \mathcal{H}^c_n$ and  accuracies $\ve > \ve^* >0$, a separator $\Theta_n$, $s \in (0,s[\Theta])$ stable grid $\mathcal{L}_{n-1}(\Theta_{n-1},\ve,\ve^*,s) \subset \mathcal{K}_{n-1}$, we call a cone $\CC \in \mathcal{C}_n$ with $S_{\CCC}\in \mathcal{L}$ and $\db(S_{\CCC},S_{C'}) < \Delta$, for some tangent cone $C'$ of $X$, and $\Delta$ as in \ref{bow1} resp. \ref{bow2}, an \textbf{$\mathcal{L}_{n-1}$-almost tangent cone} of $X$.\\
\end{definition}

\subsubsection{Stable Blow-Up Trees $\T$} \label{alt}

\bigskip

We merge the latter stabilization processes with the classical iterated blow-ups  around singular points in dimensional reduction arguments. We first describe trees of
iterated blow-ups around centers of balls belonging to finite self-similar covers. Then we define the subclass of \emph{separator trees} of trees with separators $\widehat{\Theta}$ compatible with these blow-ups. \\

 \textbf{Iterated $\mathcal{L}$-almost Tangent Cones}\label{itc2} \,  For an $\mathcal{L}$-almost tangent cone $\CC$ of $X  \in \mathcal{H}^c_n$ we choose a singular point $p \in \sigma_{\CCC}$, $ p \neq 0$. Then any tangent cone $C'$ of $\CC$ in $p$  is a Riemannian product $\R \times C^{n-1}$ for some area minimizing cone $C^{n-1} \subset \R^n$. \\

 Now we choose an almost tangent cone $\CC^*  \in \mathcal{SC}_{n-1}$ with $S_{\CCC^*} \in \mathcal{L}_{n-2}$ and so that $\db(S_{\CCC^*},S_{C^{n-1}}) < \eta_{n-2}^*$  and consider  $\R \times \CC^*$.\\

 This can be iterated until we reach a cone $\R^{n-k} \times \CC ^k$, for some $\CC^k \subset \R^{k+1}$ singular only in $0$. These cones are  $\mathcal{L}$-almost tangent cones of their  predecessor $\mathcal{L}$-almost tangent cone and for a natural operator $L$, $a$-skin adapted on $X  \in \mathcal{H}^c_n$, $a>0$, this means that
$L(\R^{n-k} \times \CC ^k)$ is $a/2^{n-k}$-skin adapted.\\

\textbf{Branch Nodes and Trees} \, We select  trees $\T$ each consisting of  finite families N[k] of \emph{branching nodes} defined inductively, in each dimension
$i \le n$ running from $k=n$ to some $k \ge 7$.  The root node is $X$. Subsequent nodes are unit spheres in almost tangent cones.\\

 Once a branching node in some N[k] has been reached, we select stable covers for this node and approximating \emph{almost tangent cones}. This stabilization step
 is labelled S[k]. The unit sphere in the selected cones are the next generation N[k-1] of nodes.\\

 We continue iteratively and define branches of subsequent series of branching nodes. Finally, we reach a branch node so that all its associated almost tangent cones are singular only in $0$, for some $k \ge 7$. The dimension of some of these terminal cones may be  $>7$ and usually it varies depending on the chosen branch we followed in $\T$.\\

 In the following definition we also introduce a $c(m)$-fold branching on each N[m+1]-node. This is a preparation to introduce compatible separators, cf.\ref{ts}-\ref{set} below:

\begin{definition}\emph{\textbf{(Blow-Up Trees)}} \label{se} \,  For  $X  \in \mathcal{H}^c_n$ and any two families $\widehat{\ve},\widehat{s}$ of  accuracy parameters and self-similarity factors
\begin{itemize}\label{e}
  \item $\widehat{\ve}= (\ve_n; (\ve_{n-1}[\kappa_{n-1}],\ve^*_{n-1}[\kappa_{n-1}])_{\kappa_{n-1}=1,...c(n-1)};...;(\ve_7[\kappa_7],\ve^*_7[\kappa_7])_{\kappa_7=1,...c(7)} ),$
  \item $\widehat{s}=\left(s_n; s_{n-1}[\kappa_{n-1}]_{\kappa_{n-1}=1,...c(n-1)};....;s_7[\kappa_7]_{\kappa_7=1,...c(7)})\right)$,
\end{itemize}
  where $c(m)$ denotes the $m$-dimensional covering number of \ref{brop}, we inductively define a \textbf{stable blow-up tree} $\T(X,\widehat{\ve},\widehat{s})$ of nodes in $\mathcal{K}$, starting
  from $X  \in \mathcal{H}^c_n$, running downwards from $n$.
The inductive construction of a tree $\T$ reads as follows:
\begin{itemize}
  \item \emph{\textbf{N[n] Root Node}}\, The only N[n] node is $X$. We choose a finite self-similar cover ${\cal{A}}^*(\Sigma_X,\ve_n,s_n)$ of $\Sigma_X \subset X$ with tangent cones $C^n_x \in {\cal{C}}_n$,  $x \in  A^*(\Sigma_X,\ve_n,s_n)$.
  \item \emph{\textbf{S[n]  Stabilization}}\, We approximate any given tangent cones $C^n_x$ by $c(n-1)$ almost tangent cones  $\CC^n[\kappa_{n-1}]_x$, $\kappa_{n-1}=1,...c(n-1)$ with $S_{\CCC[\kappa_{n-1}]_x}\in \mathcal{L}_{n-1}[\kappa_{n-1}]$, a stable grid for parameters $\ve_{n-1}[\kappa_{n-1}]$,$\ve^*_{n-1}[\kappa_{n-1}]$ and $s_{n-1}[\kappa_{n-1}]$.
  \item \emph{\textbf{N[n-1]  Branch Nodes}}\, The $S_{\CCC[\kappa_{n-1}]_x}$ are the N[n-1] nodes and we have two cases to further proceed from:
  \begin{enumerate}
    \item When $\Sigma_{S_{\CCC[\kappa_{n-1}]_x}} \v$, the process for this nodes stops here.
    \item Otherwise, we choose a finite cover ${\cal{A}}_{\kappa_{n-1}}^*(\Sigma_{S_{\CCC[\kappa_{n-1}]_x}},\ve_{n-1},s_{n-1})$ of  the singular set $\Sigma_{S_{\CCC[\kappa_{n-1}]_x}} \subset S_{\CCC[\kappa_{n-1}]_x}$ with tangent cones $C^{n-1}[\kappa_{n-1}]_p \in {\cal{C}}_{n-1}$,  for any $p \in  A^*(\Sigma_{S_{\CCC[\kappa_{n-1}]_x}},\ve_{n-1},s_{n-1})$.
  \end{enumerate}
  \item \emph{\textbf{S[n-1] Stabilization}} \, We approximate each $C^{n-1}[\kappa_{n-1}]_p$ by $c(n-2)$ almost tangent cones  $\CC^{n-1}[\kappa_{n-1}][\kappa_{n-2}]_x$, $\kappa_{n-2}=1,...c(n-2)$ with $S_{\CCC^{n-1}[\kappa_{n-1}][\kappa_{n-2}]_x}\in \mathcal{L}_{n-2}[\kappa_{n-2}]$, a stable grid for parameters $\ve_{n-1}[\kappa_{n-2}]$,$\ve^*_{n-1}[\kappa_{n-2}]$ and $s_{n-1}[\kappa_{n-2}]$.
  \item \emph{\textbf{Inductive Steps}}\,  Similarly, we define branch nodes \emph{\textbf{N[k]}} and the stabilization \emph{\textbf{S[k]}}, for $7 \le k < n-1$, until we reach terminal cones singular only in $0$.\\
\end{itemize}
\end{definition}

  \begin{remark}  \textbf{(Tree Structure of $\T$)} \label{ts}\, It may happen that two cones in different branches of $\T(X,\widehat{\ve},\widehat{s})$ coincide as elements in the space of minimal cones. For instance, one may choose $c(m)$ copies of the same almost tangent cone $\CC^{n-1}$ in each stabilization step $S[m]$. But they (more precisely, their unit spheres) do not coincide when viewed as nodes in $\T$. Indeed, we shall  label and treat each of them in a rather different way, when we consider separator trees below. Thus each $\T(X,\widehat{\ve},\widehat{s})$ is a proper tree without loops. \qed
\end{remark}

\textbf{Nested Separators} \, We use such  trees $\T$ as guidelines for the assembly of multi-pipes. We start from $\bp$-pipes in terminal nodes, build pipelines and iterated cones over pipelines until we reach the root node $X$.\\

In particular, we consider cones $C(\sp)$ over $\bp$-pipelines $\sp$ in $S_C=H \in \mathcal{K}_{n}$ of \ref{bow1} within the associated cone $C^{n+1}$. In this inductive process we employ the uniform control, we get from the stable widths  $\rho[{\mathcal{L}}_n,0]$ and $\rho[{\mathcal{L}}_n,1]$, for the $S_C=H \in \mathcal{K}_{n}$. They imply that the truncated cones  $C_{tr}(\sp)$ of \ref{tcw}, fit into \emph{uniformly} sized $\bp$-pipes in the cone $C^{n+1}$.\\

The choice of appropriate separators on the cones $C^{n+1}$ depends on these stable widths for the given separator $\Theta_n$.  We can inductively find a compatible
 separator $\Theta_{n+1}$ and also prescribe some of the properties of $\Theta_{n+1}$ by means of a given one $\Theta^\bullet_{n+1}$.

\begin{lemma}   \emph{\textbf{(Compatible Separators)}} \label{bow3} \,  Let $\Theta_n$ and $\Theta^\bullet_{n+1}$ be a given pair of separators
and $\gamma>0$ some small constant. Then there are accuracies $\ve_\kappa$ and self-similarity factors $s_\kappa$ with
 \[\gamma > \ve_1>  ... >\ve_{c(n+1)} >0, \,\, \,  \gamma> s_1> ... >s_{c(n+1)} >0\]
 and stable grids ${\cal{L}}[\kappa]:={\cal{L}}_n(\Theta_n,\ve_\kappa,\ve_\kappa^*,s_\kappa)  \subset \mathcal{K}_{n},$  $\kappa=1,...,c(n+1),$
 so that for some separator $\Theta_{n+1}$ the following relations are satisfied
 \begin{description}
 \item[$\mathbf{(R^1_n)}$]  \quad  $\pi_{n+1}[\kappa,0] \ge10 \cdot \rho[{\mathcal{L}}[\kappa], 0] \,\mm{ and }\, \rho[{\mathcal{L}}[\kappa],1] \ge 10 \cdot \pi_{n+1}[\kappa,1]$.
\item[$\mathbf{(R^2_n)}$] \quad $\rho[{\mathcal{L}}[\kappa],1] \ge 10 \cdot \rho[{\mathcal{L}}[\kappa+1],0] $    and  $\rho[{\mathcal{L}}[\kappa],0] \ge 10 \cdot \rho[{\mathcal{L}}[\kappa],1]$.
\item[$\mathbf{(R^3_n)}$]  \quad  $\pi_{n+1}[\kappa,1]/\pi_{n+1}[\kappa+1,0] \ge \pi^\bullet_{n+1}[\kappa,1]/\pi^\bullet_{n+1}[\kappa+1,0]$

    \quad  \qquad $\pi_{n+1}[\kappa,0]/\pi_{n+1}[\kappa,1] \ge  \pi^\bullet_{n+1}[\kappa,0]/\pi^\bullet_{n+1}[\kappa,1]. $
 \end{description}
\end{lemma}

 \textbf{Proof} \,  This largely follows from a review of the notions and results we have seen so far.  Here we explain the order of how they are applied. We choose the parameters
by an  induction over $\kappa$.\\

The $(R^3_n)$ conditions can obviously always be satisfied and readjusted to the needs of $(R^1_n)$ and $(R^2_n)$.\\

 Now, we start with any  $\pi_{n+1}[\kappa,0] \le \pi^\bullet_{n+1}[\kappa,0]$ and choose $\ve _1$ small enough to ensure, from \ref{re}(\ref{li}) that $\pi_{n+1}[1,0] \ge10 \cdot \rho[{\mathcal{L}}[1], 0]$. Then we choose
$\pi_{n+1}[1,1]>0$ small enough to satisfy both, $\rho[{\mathcal{L}}[1],1] \ge 10 \cdot \pi_{n+1}[1,1]$ and the second $(R^3_n)$ condition. Also, the second condition in  $(R^2_n)$ follows directly from the definition \ref{sse}(\ref{bd}). This concludes the initial induction step for $\kappa=1$.\\

Once we accomplished step $\kappa=m$, we proceed as follows: we choose $\pi_{n+1}[2,0]>0$ to the needs of the first $(R^3_n)$ condition and then we choose
$\ve_{m+1}/\ve_m>0$ small enough to also get the first conditions in $(R^1_n)$ and $(R^2_n)$. From this point on we proceed as in the case $\kappa=1$ we described before. \qed\\

When the cone $C$ that contains the truncated cone $C_{tr}(\sp)$  is an almost tangent cone $C=\CC$  for some $H^* \in \mathcal{K}_{n+1}$ a carrier pipe $\s$ of $\ref{tcw}(\ref{ttr})$, with $C_{tr}(\sp) \subset \s$, is nearly isometrically mapped under $\D^{cp}$ to a genuine one in $H^*$, cf.\ref{it} and \ref{bow0}. \\

 This way we can also transfer the finer substructure of $C_{tr}(\sp)$ onto $H^* $. Obviously this amounts to having a high accuracy for the cover of $H^* $ to ensure that this $\D^{cp}$-map is not only well-defined but also close to an isometry. We inductively iterate this process to get the following extension of \ref{bow3}, where  $(R^3_m)$ ensures the  compatibility when we reach a node along different branches of potentially varying length.

\begin{corollary}\emph{\textbf{(Separator Trees)}} \label{set} \, For  $X  \in \mathcal{H}^c_n$ and any triple $\widehat{\Theta}^\bullet,\widehat{\ve}^\bullet,\widehat{s}^\bullet$ of separators, accuracy parameters and self-similarity factors, we find another such triple
  $\widehat{\Theta},\widehat{\ve},\widehat{s}$ and a tree $\T_{\widehat{\Theta}}(X,\widehat{\ve},\widehat{s})$, we call a \textbf{separator tree}, so that
\begin{itemize}\label{e2}
  \item $\widehat{\ve} \le \widehat{\ve}^\bullet$ and $\widehat{s} \le \widehat{s}^\bullet$, i.e. the inequality holds for each component.
\item The relations $(R^1_m)$, $(R^2_m)$ and $(R^3_m)$ hold for any $7 \le m  \le n$.
\item The $\D^{cp}$-maps are supported on a neighborhood of all pipes with separator $\widehat{\Theta}$.
\end{itemize}
\end{corollary}
Also we may assume that the functional version, in the sense of \ref{bow2}, holds. That is, there are entire solutions induced from any given $v>0$ solving $L \, w=0$ on $X \setminus \Sigma_X$ on any of the product cones $\R^{n-m-1} \times C^{m+1}$, we get from the nodes $S_C$ of $\T$ in N[m], approximating the $\D$-pull backs of $v$ with $\ve_m$-accuracy.\qed

\subsubsection{Heights and Assembly of Multi-Pipes} \label{grq}
\bigskip

Now we build generalized $\bp$-pipes with internal structures matching iterated blow-ups and almost tangent cone approximations. We call them  $\bp$-multi-pipes or simply  multi-pipes $\st$.\\

 Different from typical blow-up arguments, where we push information forward to  blow-up geometries, the $\st$, in a given $X \in {\cal{H}}^c$, are unions of iterated \emph{pull-back} of simpler multi-pipes defined in iterated tangent cones via $\D$-maps within a tree $\T$. The assembly of the $\st$ is a nested iteration of the definition of basic pipelines in \ref{ab}.\\

On this way back to $X$, it will be important that the accuracy of the tangent cone approximation \emph{gradually} increases. If we simply choose a tree where all
components of $\widehat{\ve},\widehat{s}$ are very small, we observe that the stable inner and outer widths shrink too fast to allow a proper pull-back of pipes. This issue resolved from choosing separator trees $\T_{\widehat{\Theta}}(X,\widehat{\ve},\widehat{s})$.\\

\textbf{Inductive Definition of Multi-Pipes $\st$}  \, \label{dd}  Thus we choose any \emph{separator tree} $\T$ over $X$  and build multi-pipes $\st$ and multi-pipelines $\sp$ inductively until we reach $X$. This induction is \emph{not} over the dimension $m$ of the branch node N[m], but over a parameter, the \emph{height} of a node in $\T$, that measures the distance to terminal nodes.

\begin{definition}\emph{\textbf{(Heights of Terminal Order)}} \label{lt} \, On a given tree $\T(X,\widehat{\ve},\widehat{s})$ we consider the following order $\preceq$ relating pairs of nodes:
\begin{itemize}
  \item For two nodes $A \in N[k]$, $B \in N[m]$, for $k \le m$, where  $A$ belongs to the branch emanating from $B$, we write $A \preceq B$. Two nodes $P,Q \in \T$  are called   \textbf{comparable} when either $P \preceq Q$ or $Q \preceq P$.
\item The \textbf{branching number} $b(A,B)$   between two comparable nodes $A \in N[k]$, $B \in N[m]$, that is, the number of blow-ups, and stabilizations not counted separately, occurring while we travel between them, is $b(A,B)=|m-k|$.
\item  A branch node of $\T$ is called a \textbf{terminal node}, when all its
 $\mathcal{L}$-almost tangent cones, we selected in the definition of  $\T$, are singular only in $0$.
  \item An $\mathcal{L}$-almost tangent cone with unit sphere belonging to the nodes of $\T$ is called a \textbf{terminal cone}, when it is singular only in $0$.
\item The maximal branching number between a node $A$ and  comparable terminal nodes
 \[\mathbf{h}(A) := \max\{b(A,B)\,|\, B \mm{ is a terminal node comparable to } A \}\]
 is called the \textbf{height} $\mathbf{h}(A)$ of the node $A$.
\end{itemize}
\end{definition}
\smallskip

\begin{remark} \label{sah}  \quad We collect some observations readily drawn from these definitions.\\
1.  A non-terminal node may also possess terminal cones, but at least one of the cones is non-terminal.\\
2. Terminal nodes are those of height $0$ and with only isolated singular points, seen from $\Sigma \prec \sigma$ in \ref{s} below. $X$ is the unique node of largest height $\mathbf{h}(X) \le n-7$.\\
3. The singular set of nodes of height $\mathbf{h}$ has dimension $\le \mathbf{h}$, but its dimension may also be strictly smaller.\\
4.  For any branch node $A$ of $\T$ that belongs to N[m], we have $m \ge  \mathbf{h}(A)+7$. But this may be far from an equality and the deviation depends on the branch of $\T$.\qed
 \end{remark}

\textbf{Base Step 1 - Height 0} \, The induction starts from a terminal node $V$ of $\T$ in  N[m] with $\Sigma_V \n$. Then the almost tangent cones $\CC^m_p$ of $V$ in $\T$ are singular only in $0$. We  choose the $\bp$-pipes $\s\subset \CC^m_p$ for some separator $\Theta_m$ and $s \in (0,s[\pi_m[c(m),1], M])$.
In the present context, we rename $\s$  to $\s_0$ and consider
\[\s_0[\pi_m[j,0],\pi_m[j,1], s, i, j, 0] \mm{ as a \textbf{multi-pipe of height 0}.}\]

We map this $\s_0$ via a composition of the $\D$-map to the tangent cone $C^m_p$ in  $\T(X,\widehat{\Theta},\widehat{\ve},\widehat{s})$, and the $\D$-map from this cone to $V$. For $\ve \in (0,\ve_0)$ as in \ref{it} and some stabilization accuracy $\ve^*\in (0,\ve)$ we consider the disjoint union of these pipes for each ball of the cover ${\cal{A}}^*$. This defines a $\bp$-pipeline $\sp$,cf.\ref{ab} and \ref{itr}, we now call
  \[\sp_0  :=\dot\bigcup_{p \in A[i,j]^*(\Sigma_V)}\D^{cp}(\s_0[\pi_m[j,0],\pi_m[j,1], s, i, j, 0]) \mm{ a \textbf{multi-pipeline  of height 0}.}\]\

\textbf{Base Step 2 - Height 1} \,  We continue with the first non-trivial step, the definition of multi-pipes of height 1: we have two cases
 \begin{itemize}
   \item  $V=X$ and thus we have $m=n$. Then, we keep $\sp_0$ and the process stops.
   \item   $V=S_{\CC^{m+1}_x}$ for an almost tangent cone $\CC^{m+1}_x \in {\cal{C}}_{m+1}$ of another branch node $V'$ in N[m+1]. Then, we consider the cone $C(\sp_0) \subset \CC^{m+1}_x$ over $\sp_0$. Since $\T$ is a separator tree we know that for some $j \in \{1,...,c(m+1)\}$
\begin{equation}\label{bi}
C(\sp_0) \subset \P^*(0, \pi_{m+1}[j,0]) \setminus  \P^*(0, \pi_{m+1}[j,1]) \subset \CC^{m+1}_x.
\end{equation}
 \end{itemize}

This means that, after suitable truncation, $C(\sp_0)$ fits into a carrier $\bp$-pipe  of $\CC^{m+1}_x$. This motivates the following definition of multi-pipes one may think of as wheels composed of hubs  and spokes:

\begin{definition}\emph{\textbf{(Multi-Pipes of Height 1)}} \label{st}  The  set  $\st_1(\sp_0) \subset C^{m+1}[x]$  defined by
\[\st_1(\sp_0):=\st_{1}(\sp_{0})_{spokes} \, \cup \, \st_{1}(\sp_{0})_{hub}, \mm{ where }\]
\begin{itemize}
  \item $\st_{1}(\sp_{0})_{spokes}:=C(\sp_0) \cap B_{5 \cdot s_{m+1}^i}(0)\setminus  B_{\pi_{m+1}[j,1] \cdot s_{m+1}^i}(0)$
  \item $\st_{1}(\sp_{0})_{hub}:=B_{\pi_{m+1}[j,0]\cdot s_{m+1}^i}(0)\setminus \big( \P^*(0, \pi_{m+1}[j,1]) \cup B_{\pi_{m+1}[j,1]\cdot s_{m+1}^i}(0)\big)$
\end{itemize}
is called a \textbf{multi-pipe}  of height $1$ over the multi-pipeline $\sp_0$ of height $0$. We observe $\st_1(\sp_0) \subset \s[\Theta_{m+1},s_{m+1}, 0]$ and call  $\s[\Theta_{m+1},s_{m+1}, 0]$ a \textbf{carrier pipe} of $\st_1(\sp_0)$.\\
\end{definition}

\textbf{Inductive Step -  Height $\mathbf{h}$+1} \,   We may assume that
we have already defined multi-pipes $\st_k$ with
\[\st_k [\Theta_m,s_m,...,\Theta_{m-k},s_{m-k},0] \mm{ within carrier pipes } \s[\Theta_m,s_m,0] \subset C^* \in {\cal{C}}_m\] in all those almost tangent cones $C^*$ where
$S_{C^*}$ has height $k \le \mathbf{h}$. Now let  $V$ in N[m] be a node of height $\mathbf{h} +1$. We  repeat the assembly step above for the carrier pipes $\s[\Theta_{m},s_{m},x]$ in  the almost tangent cones $C^m[x]$ and define the carrier pipeline \[\sp  :=\dot\bigcup_{p \in A[i,j]^*(\Sigma_V)}\D^{cp}(\s[\pi_m[j,0],\pi_m[j,1], s, i, j, 0])  \subset V,\]
where $\D^{cp}$ is the composition of the $\D$-map from the almost tangent cone $\CC^{m}_x$ to a tangent cone $C^m_x$ of $V$ and the $\D$-map from that cone to $V$.\\

Thus we can repeat the definition of a multi-pipe of height 1 for carrier pipes and distinguish between two cases:

\begin{itemize}
  \item If $V=X$, we keep $\sp$ and the process stops.
  \item In the alternative case $S_{C^{m+1}[x]}$, we consider the cone $C(\sp) \subset \CC^{m+1}_x$ over $\sp_0$  and note that for a suitable $j \in \{1,...,c(m+1)\}$:
\begin{equation}\label{bi}
C(\sp) \subset \P^*(0, \pi_{m+1}[j,0]) \setminus  \P^*(0, \pi_{m+1}[j,1]) \subset \CC^{m+1}_x.
\end{equation}
\end{itemize}

Now we insert an inductively already given multi-pipe into this one. The heights of these inserted multi-pipes may vary between $0$ and $\mathbf{h}$. Thus given the definition of multi-pipes $\st_k$ of height $k \le \mathbf{h}$, we extend the definition to the height $\mathbf{h} +1$ case as follows.

\begin{definition} \emph{\textbf{(Multi-Pipes of Height $\mathbf{h}$+1)}} \label{st2} \, Let $\st_k\subset \CC^{m+1}_x$ be an almost tangent cone $C^{m+1}[x]$, with $S_{\CCC_x}$ of height $k$,
of a node of height $\mathbf{h}$, $S_{\CCC_z}$, $x \in A^*(\Sigma_{S_{\CCC_z}})$. Then we define
the \textbf{multi-pipeline $\sp_\mathbf{h}$ of height }$\mathbf{h}$ in  $S_{\CCC_z}$ as
   \begin{equation}\label{dst}
   \sp_{\mathbf{h}}:=\dot\bigcup_{x \in A^*(\Sigma_{S_{\CCC_z}})}\D^{cp}\left(\st_k[\Theta_m,s_m,...,\Theta_{m-k},s_{m-k},0] \right)
 \end{equation}
\[\subset \dot\bigcup_{x \in A^*(\Sigma_{S_{\CCC_z}})}\D^{cp}\left(\s[\Theta_{m+1},s_{m+1},x]\right) \subset S_{\CCC_z}.\]
The subset $\st_{\mathbf{h}+1}(\sp_{\mathbf{h}})\subset \CC^{m+1}_x$, defined by
\[\st_{\mathbf{h}+1}(\sp_{\mathbf{h}}):=\st_{\mathbf{h}+1}(\sp_{\mathbf{h}})_{spokes} \, \cup \, \st_{\mathbf{h}+1}(\sp_{\mathbf{h}})_{hub}, \mm{ where }\]
\begin{itemize}\label{bi3}
  \item $\st_{\mathbf{h}+1}(\sp_{\mathbf{h}})_{spokes}:=C(\sp_\mathbf{h})  \cap B_{5 \cdot s_{m+1}^i}(0)\setminus  B_{\pi_{m+1}[j,1] \cdot s_{m+1}^i}(0)$
  \item $\st_{\mathbf{h}+1}(\sp_{\mathbf{h}})_{hub}:=B_{\pi_{m+1}[j,0]\cdot s_{m+1}^i}(0)\setminus \big( \P^*(0, \pi_{m+1}[j,1]) \cup B_{\pi_{m+1}[j,1]\cdot s_{m+1}^i}(0)\big)$
\end{itemize}
is called a \textbf{multi-pipe}  $\st_{\mathbf{h}+1}$  of height $\mathbf{h}+1$ over the multi-pipeline $\sp_\mathbf{h}$ of height $\mathbf{h}$. Again, to reproduce the induction hypothesis we note that\\

 $\st_{\mathbf{h}+1}(\sp_{\mathbf{h}}) \subset \s[\Theta_{m+1},s_{m+1},x]$ and we call  this $\s$ a \textbf{carrier pipe} of $\st_{\mathbf{h}+1}(\sp_{\mathbf{h}})$.\\
\end{definition}
We observe that $\st_{\mathbf{h}+1}(\sp_{\mathbf{h}})$ resp. its parameters depend only on the dimensions $m+1$ to $m-\mathbf{h}$ or, differently stated, on $m$ and $\mathbf{h}$.\\

The ultimate outcome are multi-pipelines $\sp_{\mathbf{h}(X)}$ of height $\mathbf{h}(X)$ on $X$. They clearly substantially depend on the chosen tree $\T$.   \qed

\begin{remark} \label{tz} \, 1. The distinction between hubs and spokes is made to explain the inductive construction multi-pipe(line)s. But we observe that each of spokes is composition of spokes over hubs and spokes of lower height. Inductively,  this means, $\sp \subset H$ has the structure of a union of $\st$'s, but less structured, of $\D$-map distorted tapered products of a varying number of intervals with an initial pipe for a cone singular only in $0$. We call these raw building blocks the \textbf{elements} of $\sp$.\\

2.  In the definition of trees and multi-pipes we used the lower-dimensional  almost tangent cones $\CC^k$. However with each $\CC^k$ that belongs to some unit sphere $S_{C^*}$ in a cone $C^*$, the product cone $\R\times\CC^k$ can also be considered as an almost tangent cone of $C^*$. Iteratively, $\R^{n-k}\times\CC^k$ is almost tangential to $X$. Moreover, for small separators we can use  $\R^{n-k} \times \CC^k$ as a tapered approximation of a spoke, cf.\ref{tst} for details.\\

3. Hubs and spokes have non-trivial intersections. They cross \textbf{perpendicularly} along their boundaries:
\begin{equation}\label{nti}
\st_{\mathbf{h}+1}(\sp_{\mathbf{h}})_{spokes} \, \cap \, \st_{\mathbf{h}+1}(\sp_{\mathbf{h}})_{hub} = C(\sp_\mathbf{h})  \cap B_{\pi_{m+1}[j,0]\cdot s_{m+1}^i}(0)\setminus   B_{\pi_{m+1}[j,1]\cdot s_{m+1}^i}(0).
\end{equation} \qed \end{remark}

\subsubsection{Multi-Pipes and Singularities} \label{thr}
\bigskip

The definition of $\bp$-multi-pipes made essential use of stable grids and almost tangent cones. The multi-pipes inherit and refine the important
property \ref{bow1}(\ref{tt0}) saying that $\bp$-pipelines stably surround the singular set. Here we discuss some finer details.\\

We consider  $\bp$-pipes $\s[\Theta, s, 0] \subset C$, for  $C \in \mathcal{SC}_{n}$, more closely. They have been used in the definition of multi-pipes $\st$ as supporting carriers. In turn, for suitably chosen parameters, $\s[\Theta, s, 0] \subset C$ is not a subset, but at least surrounded by multi-pipes. This gives us the following easy to use decomposition of multi-pipes controlled by some singular threshold:\\

In dimension $n=7$ we know that $C$ is singular precisely in $0$. For a  single pipe $\s[\Theta, s, 0] \subset C$  we have $a(C):=\sup_{S_C} \bp < \infty$ and for any separator $\Theta$ with $\pi_n[1,0]  < 1/a(C)$  this says $\P^*(p, \pi_n[j,0]) \v$ and  $\s[\Theta, s, 0] = B_{\pi_n[j,0] \cdot s^i}(0)  \setminus B_{\pi_n[j,1] \cdot s^i}(0)$.   Since $\mathcal{SC}_{7}$ is compact and $\bp$ is naturally defined, we also have $a_7:=\sup_{C \in \mathcal{SC}_{7}} a(C)< \infty$.\\

In higher dimensions, $C$ may also be singular outside $0$. In this case, we have $\P^*(p, \pi_n[j,0]) \n$ for all separators $\Theta$. This also implies that, although for a cone $C$ singular only in $0$ we get $a(C) < \infty$, there is no common upper bound for $a(C)$ on the subspace of these cones.  But we can use stable grids of \ref{bow1} to get a replacement for this upper bound.

 \begin{corollary}   \emph{\textbf{(Singular Threshold)}} \label{strs} \, For any accuracy $\ve >0$, a separator $\Theta_n$, $s \in (0,s[\Theta])$ and a stabilization accuracy $\ve^*\in (0,\ve)$, there is an $a_n(\Theta,\ve,\ve^*,s) >0$ so that for any $S_C \in \mathcal{K}_{n}$ and $x \in S_C \setminus \Sigma_{S_C}$, either $\bp(x) \le a_n$ or $x \in \I(S_C,\rho^*[1])$.
\end{corollary}

Clearly, we just choose $a_n:= \rho^*[1]$. The point is the interpretation: we get a decomposition into a part of $C$ with a uniform bound on $\bp$ on $S_C$ and its complement that is compressed into the spoke over a pipeline in $S_C$.\qed

Next we formulate an important relation, symbolically  "$\Sigma\prec \sigma$", that relates the singular set of a given $X \in {\cal{H}}^c_{n+1}$ to that of \emph{almost tangent cones}. In fact, it generally does not hold for proper tangent cones of $X$.\\

 In simple terms, $\Sigma\prec \sigma$ says, $\Sigma_X \cap B_{5 \cdot s^i}(p) \setminus B_{5 \cdot s^i \cdot \delta}(p)$, for $\delta \in (0,1)$, stays in a small neighborhood of the singular set  $\sigma$, of some $\mathcal{L}$-almost tangent cone, when viewed extrinsically in a common ambient space.\\

This implies that not only $\bp$-pipes, but also cones  $C(\sp)$ over $\bp$-pipelines within almost tangent cones, locally surround $\Sigma_X$. Iteratively applied, this shows that also  the multi-pipes surround $\Sigma_X$.\\

We express this relation intrinsically in terms of a disjointness of  complements of the singular sets we may compare to subsets of $X^{n+1}$ via $\D$-maps.

\begin{proposition}\emph{\textbf{(\textbf{Surrounding Relation} "$\Sigma\prec \sigma$")}} \label{s} \, Let $X \in {\cal{H}}^c_{n+1}$,  some stable grid  $\mathcal{L}_n(\Theta,\ve,\ve^*,s) \subset \mathcal{K}_{n}$ for some given accuracy $\ve >0$,  separator $\Theta_n$, self-similarity factor $s \in (0,s[\Theta])$
$\delta \in (0,1)$ be given.\\

Then, there are some $\ve_X, s_X>0$ and a self-similar
cover ${\cal{A}}(\Sigma_X,\ve_X,s_X)$ of $\Sigma_X$, so that for any $x \in A^*(\Sigma_X)$ and any $\mathcal{L}$-almost tangent cone $\CC_x$ in  $x \in  A[j] \subset  A^*(\Sigma_X)$
\begin{equation}\label{ee}
s_X^{-j} \cdot \left(\Sigma_X \cap (B_{5 \cdot s_X^j}(x) \setminus B_{5 \cdot s_X^j\cdot \delta}(x))\right) \cap ... \end{equation}
\[...\D_{C(x) \ra s_X^{-j} \cdot X} \big( C(\D_{S_{\CCC_x} \ra S_{C(x)}}\big(\E(S_{\CCC_x},\rho^*[1](S_{\CCC_x},\Theta_n,\ve,\ve^*,s))\big))\big) \v.\]
and for an additionally given $\Theta^\bullet_{n+1}$ we can choose  $\ve_X, s_X>0$ so that, in addition, the $\Theta_{n+1}$ derived from \ref{bow3} satisfies
\begin{equation}\label{ee0}
 \left(\Sigma_X \cap (B_{5 \cdot s_X^j}(x) \setminus B_{5 \cdot s_X^j\cdot \delta}(x))\right) \cap \s[\Theta_{n+1}, s_X, x] \v\end{equation}
\end{proposition}

\textbf{Proof} \, This is mainly a consequence of the way we defined the stable grid: since ${\cal{L}}$ is finite, there is a common upper bound $\Gamma_{\cal{L}} >0$ in  \ref{bow1} so that for any $H' \in \mathcal{K}_{n}$ with $\db(H,H') <\eta_H$ with  $H \in \cal{L}$:
\[ \bp_{H'} \le \Gamma_{\cal{L}} \mm{ on } \D(\E(S_{\CCC_x},\rho^*[1])) \subset H'.\]
Thus, we also get some $\Gamma_{\cal{L}}^* >0$ so that, when $\ve_X$ is sufficiently small:
\[\bp_X \le \Gamma_{\cal{L}}^* \mm{ on } (B_{5 \cdot s_X^j}(x) \setminus B_{5 \cdot s_X^j\cdot \delta}(x)) \cap \D \big( C(\D\big(\E(S_{\CCC_x},\rho^*[1])\big))\big).\]
The remaining assertions follow easily from \ref{bow3} and the definitions.\qed

There are forerunners of $\Sigma\prec \sigma$ already in Federer's classical inductive proof that the codimension of $\Sigma$ is $\ge 7$, cf.[Gi], Ch.11. In turn,  self-similar covers and $\Sigma\prec \sigma$ can be used to rewrite Federer's argument in a concise way. Indeed, one can inductively estimate the area $\p (\sp)$ and use the isoperimetric inequality to estimate the volume growth of distance tubes around $\Sigma$ from the relation $\delta_{\bp} \le L  \cdot dist(x,\Sigma)$. We leave the elementary details to the interested reader.\\

\setcounter{section}{4}
\renewcommand{\thesubsection}{\thesection}
\subsection{Scal$>0$-Heredity with Surgery} \label{smo0} 
\bigskip

We start from a singular area minimizer $H^n \subset M^{n+1}$ within some smooth compact manifold $M$ with $scal \ge 0$. Our main goal, we accomplish in this and the next chapter, is to deform $(H \setminus \Sigma,g_H)$ to $scal >0$-geometries with horizons shielding $\Sigma$ against intrusions of area minimizing sequences of hypersurfaces within $H \setminus \Sigma$.\\

In this chapter, we define the basic deformations and study their bending effects for model cases of product cones. Then we merge versions of these deformations, localized to match separator trees $\T$, to a horizon metric $g_{\T}$. The proof that the two steps of localization and merging do not harm the shielding effects of $g_{\T}$ is a technically separate topic explored in the next chapter Ch.\ref{smo}.

\subsubsection{Analysis of Conformal Laplacians} \label{ma}
\bigskip

We resume the discussion of the potential theory of conformal Laplacians from Ch.\ref{mre}.C and D for an area minimizing hypersurface $H^n \subset M^{n+1}$, in some manifold $M^{n+1}$ with $scal \ge  0$, where $M^{n+1}$ is either compact or the flat $\R^{n+1}$, and $\Sigma_H \n$.\\

Then $L_H$ is  skin adapted. In particular, there is a positive $\lambda > 0$ with $ \int_H f \cdot L_H f \,  dA \, \ge \, \lambda \cdot \int_H \bp^2\cdot f^2 dA,$
for any smooth $f$ with $supp(f) \subset H \setminus \Sigma$.\\

The largest such $\lambda > 0$ is the generalized principal eigenvalue $\lambda^{\bp}_{H,L_H} >0$ of the  skin conformal Laplacian $\delta_{\bp}^2 \cdot L_H$.  Thus the  operator
 \[ L_{H,\lambda}:= - \triangle   + \frac{n-2}{4 (n-1)}  \cdot scal_H  -\lambda  \cdot  \bp^2 \]
$L_{H,\lambda}$ is skin adapted for $\lambda <\lambda^{\bp}_{H,L_H}$ and under blow-ups, any of the (iterated almost) tangent cones  $C$ inherits the skin adaptedness of $L_{C,\lambda}$ from that of  $L_{H,\lambda}$ on $H$.\\

For $\lambda \in (0,\lambda^{\bp}_{H,L_H})$ and $\alpha \in (0,1)$, the Martin theory for $L_{H,\lambda}$ shows there are  many $C^{2,\alpha}$-solutions $u >0$  of $L_{H,\lambda} \,w=0$ and for each of them we have $scal(u^{4/n-2} \cdot g_H)>0$.\\

As we will see below, choosing a sufficiently small $\lambda >0$ is not only needed to ensure the skin adaptedness of $ L_{H,\lambda}$. It is also necessary  to derive particular estimates for the growth rates of solutions of $L_{C,\lambda} \, w=0$ towards $0 \in C$, on tangent cones  $C$.\\

 Namely, we need growth rates strong enough to get, under conformal deformation, spreading open ends towards $0$. To this end, we shrink $\lambda>0$ to a small size. This matches the observation that such bending processes always sacrifice some of the positivity of $scal$.\\

\textbf{Results for Cones} \, We start with an area minimizing cone  $C^n\subset \R^{n+1}$ which can be written as a product $\R^{n-k} \times C^k$, where $C^k \subset \R^{k+1}$ is again an area minimizing cone equipped with
$g_{C^n}=g_{\R^{n-k}} + g_{C^k}$.  This suggests to consider dimensionally shifted conformal Laplacian on $\{0\} \times C^k$.
\[L_{C^k,n} :=-\Delta - \frac{n-2}{4 (n-1)} \cdot |A|^2,\,  \mm{ for }
 n \ge k \mm{ and } (L_{C^k,n})_\lambda  = L_{C^k,n}-\lambda  \cdot  \bp^2\]
Recall that for $C$ and actually any area minimizing  hypersurface $H \subset \R^{n+1}$, we have
\[scal_H\equiv -|A_H|^2, \mm{ from the Gauss-Codazzi equations.}\]
We summarize the following basic results from [L2],2.25, 3.18, 5.12 and  5.19:

\begin{proposition} \emph{\textbf{(Dimensionally Shifted $L_C$)}} \label{ddd0} \, For a singular area minimizing cone $C^k \subset \R^{k+1}$ and $n \ge k \ge 7$ we have:
 \begin{itemize}
\item $L_{C,n}$ is skin adapted and  the principal eigenvalue can be uniformly lower estimated for all $k$-dimensional cones by a positive constant $\Lambda^*_k >0$ independent of $n \ge k$.
    \item For $\lambda < \Lambda^*_k$, $(L_{C,n})_\lambda$ is also skin adapted.
  \end{itemize}
\end{proposition}

For skin adapted elliptic operators on skin uniform spaces we have a potential theory that, broadly speaking, does not render details not already visible from the geometric shape of the underlying space. This is illustrated by the following result.

\begin{proposition} \emph{\textbf{(Martin Boundaries)}} \label{ddd} \, The Martin boundary $\p_M (C^k,(L_{C,n})_\lambda)$ is homeomorphic to $\widehat{\sigma_C}$, the one-point compactification of the singular set $\sigma_C$. Thus we have
\[\p_M (C^k,(L_{C,n})_\lambda) \cong [0,1]\times \Sigma_{\p B_1(0) \cap C}/\sim,\]
where  $\Sigma_{\p B_1(0) \cap  C}=\p B_1(0) \cap  \sigma_{C}$ is the singular set of $S_C=\p B_1(0) \cap C$, with  $x\sim y$,  if $x , y \in \{0\} \times \Sigma_{\p B_1(0)
        \cap C}$\, or \, $x , y \in\{1\} \times \Sigma_{\p B_1(0) \cap C}$. Also, we get some growth estimates:\\
\begin{enumerate}

 \item The two distinguished points $[0], [1] \in \p_M (C^k,(L_{C,n})_\lambda)$, viewed as functions $\Psi_{C,-}=[0]$ and $\Psi_{C,+}=[1]$  on $C \setminus \sigma$, can be written as
\[\Psi_\pm(\omega,r)=\Psi_{C,\pm}[n,k](\omega,r) = \psi_C(\omega) \cdot r^{\alpha_\pm},\]
for $(\omega,r) \in S_C \setminus \sigma \times \R$, with {\small $\alpha_\pm =\alpha_\pm[n,k]= - \frac{k-2}{2} \pm \sqrt{ \Big( \frac{k-2}{2} \Big)^2 + \mu_{C,(L_{C,n})^\times_\lambda}}$}.
\item   There is a constant $\lambda^*_k \in (0,\Lambda^*_k)$, depending only on $k$, so that  for any
 $n \ge k$:
  \[\mu_{C,(L_{C,n})^\times_\lambda} \ge - 1/3 \cdot \left(\frac{k-2}{2}\right)^2,\mm{ for any }\lambda \le \lambda^*_k.\]
\item For $\lambda \in (0,\lambda^*_k]$, the growth rates of the two solutions $\Psi_\pm[n,k](\omega,r) = \psi_C(\omega) \cdot r^{\alpha_\pm}$ satisfy  $\mu^*_{C,(L_{C,n})^\times_\lambda} <-\eta^*_\lambda <0$ and hence
 {\small \[-\vartheta^*_\lambda >  \alpha_+ \ge - (1- \sqrt{2/3}) \cdot \frac{k-2}{2}>  - \frac{k-2}{2} >- (1+ \sqrt{2/3}) \cdot \frac{k-2}{2} \ge \alpha_- > \vartheta^*_\lambda -(k-2),\]}
 for some constants $\eta^*_\lambda, \vartheta^*_\lambda>0$, depending only on $\lambda$ and $k$.\\
\end{enumerate}
\end{proposition}

\textbf{Distinguished Solutions} \, The growth rates of $\Psi_+[n,n]$ on $\R^{n-k} \times C^k$ towards $\R^{n-k} \times \{0\}$ and of $\Psi_+[n,k]$ on $C^k$ towards $\{0\}$ coincide, but those of $\Psi_-[n,n]$ and $\Psi_-[n,k]$ differ. In what follows, it will be important to understand this asymmetry between $\Psi_+[n,n]$ and $\Psi_-[n,n]$:\\

$\Psi_+[n,n]$ shares the $\R^{n-k}$-translation symmetry with the underlying space $\R^{n-k} \times C^k$. Namely, the point at infinity $\infty \in \widehat{\sigma_{C^n}}$ is the support for the Dirac measure that defines  $\Psi_+[n,n]$ in its Martin integral representation. But $\infty$ does not move under the $\R^{n-k}$-translations and hence we have $\Psi_+[n,n](x,y)=\Psi_+[n,k](y)$, for $(x,y) \in \R^{n-k} \times C^k$.\\

In turn, $\Psi_-[n,n]$ is \emph{not} invariant under $\R^{n-k}$-translations. But we can average  $\Psi_-[n,n]$ on $\R^{n-k} \times C^k$ along the $\R^{n-k}$-factor and get a solution  $\Psi_{\bigodot}[n,k]$. This function is no longer a minimal solution on $\R^{n-k} \times C^k$, but it admits a simple Martin representation, showing its relation to the problem on $C^k$:
 \begin{equation}\label{ias}
 \Psi_{\bigodot}[n,k](a,b)=\int_{\widehat{\sigma_{C^n}}} k((a,b);y) \, d \mu_{\Psi_{\bigodot}[n,k]}(y)=
\end{equation}
\[\int_{\R^{n-k} \times \{0\}} \Psi_-[n,n]((a,b)-(y_1,y_2)) \, d \lambda_{\,\R^{n-k} \times \{0\}}(y_1)= \kappa_n \cdot \Psi_-[n,k](b),\]
for some $\kappa_n >0$, any  $(a,b) \in \R^{n-k} \times C^k$ and where $\lambda_{\,\R^{n-k} \times \{0\}}$ is the Lebesgue measure on $\R^{n-k} \times \{0\} \subset \R^{n-k} \times C^k$, while the complete measure $\mu_{\Psi_{\bigodot}[n,k]}$, on $\widehat{\sigma_{C^n}}$, vanishes on $\widehat{\sigma_{C^n}}\setminus \R^{n-k} \times \{0\}$.\\

\textbf{Spherical Estimates} \, We also need to ensure that the radial growth rate of $\Psi_\pm[n,k]$ remains uniformly lower bounded while we move along $S_C=\p B_1(0) \cap C$ towards $\Sigma_{S_C}$ and so that the bound is also \emph{independent} of $C^k \in {\cal{SC}}_k$.\\

 We get such a bound from the following basic estimates for the spherical component $\psi_C(\omega)$ of  $\Psi_\pm[n,k](\omega,r) = \psi_C(\omega) \cdot r^{\alpha_\pm}$ defined on $S_C \setminus \Sigma_{S_C}$, cf.[L2], Prop.5.20:
\begin{proposition}\emph{\textbf{(Global Harnack Estimates for $\psi_C$)}}   \label{skk} \, For any $n \ge k \ge 7$ and $\Psi_\pm[n,k](\omega,r) = \psi_C(\omega) \cdot r^{\alpha_\pm}$, there are constants $a_{n,k,\lambda} >0$, depending only on $n,k,\lambda$, and $b_{n,k,\lambda,\rho}$ additionally depending on $\rho>0$, but not on the cone $C^k$:
\begin{enumerate}
  \item $|\psi_C|_{L^p(S_C \setminus \Sigma_{S_C})}<\infty$, for $ p< \frac {k-1} {k-3}$.
  \item $|\psi_C|_{L^1(S_C \setminus \Sigma_{S_C})} \le a_{n,k,\lambda} \cdot \inf_{\omega \in S_C \setminus \Sigma_{S_C}}  \psi_C(\omega)$.
  \item $\sup_{\omega \in\E(\rho)}\psi_C(\omega) \le  b_{n,k,\lambda,\rho} \cdot|\psi_C|_{L^1(S_C \setminus \Sigma_{S_C})}$, for $\E(\rho)=\{x \in S_C \,|\, \bp(x) \le \rho^{-1}\}$.
\end{enumerate}
Similarly, $|v|_{L^q(B_1(0) \cap C^k)}<\infty$, for $q< \frac {k} {k-2}$, for any solution $v>0$ of  $(L_{C^k,n})_\lambda \, w=0$.
\end{proposition}

We can use the integrability (i) to gauge frequently used standard solutions on $C$ without selecting a regular basepoint. In particular, for the distinguished solution $\Psi_{C,+}$ and the solutions  $\Psi_{\bigodot}[n,k]$, we may henceforth assume
\begin{equation}\label{gauge}
|\Psi_{C,+}|_{L^1(B_1(0) \cap C^n)}=|\Psi_{\bigodot}[n,k]|_{L^1(B_1(0) \cap C^n)}=1.
\end{equation}
\smallskip

\textbf{Minimal Metrics} \, What makes $\Psi_{C,+}$ particularly interesting is its minimal growth towards the entire singular set $\sigma_C$. This follows from the fact that the associated Radon measure is only supported in $\infty \in \widehat{\sigma_{C^n}}$, cf.[L2],Lemma 4.5 and Prop. 5.8. In turn, this also characterizes $\Psi_{C,+}$ from Martin theory.\\

Now we consider the cases of singular minimizers $H \in {\cal{H}}^c_n$ and $H \in H^n \subset M^{n+1}$:\\

For ${\cal{H}}^{\R}_n$ and $\lambda \in (0,\lambda^{\bp}_{H,L_H})$, we know that $\widehat{\Sigma}=\Sigma \cup \{\infty\}$ and we have a unique solution $\Psi_{H,+}>0$
$|\Psi_{H,+}|_{L^1(B_1(0) \cap C^n)}=1$ that is associated to the Dirac measure in $\infty$. It has minimal growth along $\Sigma_H$.\\

For ${\cal{H}}^c_n$ and $\lambda \in (0,\lambda^{\bp}_{H,L_H})$, we also seek for a solution with minimal growth along $\Sigma_H$. But here we need to argue in a slightly different way.  The point is that due to the skin adaptedness of $L_{H,\lambda}$ there is no regular positive solution with minimal growth along the entire compact set $\Sigma_H$, cf.[L2],5.3.\\

Instead, we use the minimal $L_{H,\lambda}$-Green's function $\Phi(H)(\cdot) =G(q_H,\cdot)$, so that  $\Phi(H)(q^*_H )=1$, for some given $q^*_H \neq q_H \in H \setminus \Sigma$. $\Phi(H)$ has minimal growth along $\Sigma_H$ and we consider the singular  metric $g_H^+:=\Phi(H)^{4/(n-2)} \cdot g_H$. \\

 The singularity of the deformation $\Phi(H)^{4/(n-2)} \cdot g_H$ in $q_H$ does not cause geometric problems. The growth of $\Phi(H)(x)$ near $q_H$ is of order $r^{-(n-2)}$, for  $r=dist(x,q_H)$, and, therefore,  the deformation gives an asymptotic flat end. This means, it remains invisible to area minimizers within $H$, cf. the proof of \ref{merg2} for details.\\

The geometry $\Phi(H)^{4/(n-2)} \cdot g_H$ is essentially invariant under changes of the basepoint: due
to the skin uniformity of $H \setminus \Sigma$ we know from [L2],Th.7, that for any two basepoints  $q_H(i) \in H \setminus \Sigma$,$i=1,2$, the \emph{quotient} $\Phi(H)(1)/\Phi(H)(2)$ of the associated Green's functions $\Phi(H)(i)=G(q_H(i),\cdot)$,  admits a positive \emph{continuous extension} to the entire singular space $H$, outside
the $q_H(i)$.\\

\begin{definition}  \emph{\textbf{(Distinguished Metrics)}} \label{mii} \, We label some particular geometries.
For $\lambda \in (0,\lambda^*_k)$, we consider two $scal>0$-metrics on $\R^{n-k} \times C^k$, $C^k \in {\cal{SC}}_k$:
\begin{itemize}
  \item $g^+_{C^n}:=\Psi^{4/(n-2)}_+[n,n] \cdot (g_{\R^{n-k}} + g_{C^k})$, the \textbf{minimal metric} relative $(L_{C,n})_\lambda$.
  \item $g^{\bigodot}_{C^n}:=\Psi_{\bigodot}^{4/(n-2)}[n,k] \cdot (g_{\R^{n-k}} + g_{C^k})$, the \textbf{axial metric} along $\R^{n-k} \times \{0\} \subset \sigma_C$  relative $(L_{C,n})_\lambda$.
      \end{itemize}
      For singular area minimizing hypersurfaces $H^n \subset M^{n+1}$, with generalized principal eigenvalue $\lambda^{\bp}_{H,L_H} >0$. For $\lambda \in (0,\lambda^{\bp}_{H,L_H})$, we distinguish between two cases and call the $scal>0$-metrics
       \begin{itemize}
\item $g^+_H:=\Psi_+^{4/(n-2)} \cdot g_H$, the \textbf{minimal metric} relative $L_{H,\lambda}$, when ${\cal{H}}^{\R}_n$,
\item $g^+_H:=\Phi(H)^{4/(n-2)} \cdot g_H$, a \textbf{minimal metric} relative $L_{H,\lambda}$, when ${\cal{H}}^c_n$.
\end{itemize}
\end{definition}
Note that both metrics $g^+_{C^n}$ and $g^{\bigodot}_{C^n}$ reproduce under scalings since they also describe conical metrics. This makes the assumption (\ref{gauge}) primarily relevant when we  need to compare the two functions $\Psi_{C,+}$ and $\Psi_{\bigodot}[n,k]$ in later constructions.\\

\textbf{Custom Downstream Processes}  \, As explained in Ch.\ref{mre},F. we use the latter potential theoretic results in a preparatory \emph{downstream
process}. Here we select the analytic parameters of this process to manage the subsequent \emph{upstream process}, in Ch.\ref{smoo}, where localized secondary deformations  for all cones in a suitable separator tree $\T$ are pulled-back to $H$. To this end, we make some custom choices for the eigenvalue problems of the conformal Laplacians and their associated eigenfunctions.\\

For an area minimizing hypersurface $H \subset M$, in some compact manifold with $scal_M \ge  0$ and $\Sigma_H \n$, we consider the following setup:\\

  \textbf{A.} \, We consider $L_{H,\lambda}$ and $(L_{C^k,n})_\lambda$, for all $C=\R^{n-k} \times C^k \in \T$, in a separator tree $\T$, we specify later on,  and $n \ge k \ge 7$. We select the following value for $\lambda$: \[\lambda := \min\{\lambda^{\bp}_{H,L_H},\lambda^*_k, n \ge k \ge 7 \}/2 >0.\]

The chosen $\lambda$ makes the conformal Laplacians, considered above, skin adapted and all potential theoretical results we discussed so far become applicable to $H$ and any of its (almost) tangent cones.\\

 \textbf{B.} \,  We choose a particular initial $scal>0$-conformal deformation of $H \setminus \Sigma$: a  minimal metric $g^+_H=\Phi(H)^{4/(n-2)} \cdot g_H$ relative $L_{H,\lambda}$.\\

Then the induced solutions $\Phi(C)$, for cones $C \in \T$,  also have minimal growth along $\sigma_C$, viewed as an open subset $\sigma_C \subset \widehat{\sigma_C}$, cf.[L2],5.18, and hence, from \ref{ddd}(iv) we get: $\Phi(C)=\Psi_{C,+}$. Thus the metrics $\Phi(C)^{4/(n-2)} \cdot g_C$ are the minimal metrics $g^+_{C}$ and, under blow-ups, they $C^{2,\alpha}$-compactly asymptotically approximate $\Phi(H)^{4/(n-2)} \cdot g_H$. \qed

\subsubsection{Horizons for Minimal Metrics on Cones} \label{e-barriers2}
\bigskip

In this section we build the basic model of an $scal>0$-horizon. All later horizon constructions are localized, deformed or superposed variants of this basic version.\\

We start with a comparison  of the bending effects of the metrics $g^+_{C^n}$ and $g^{\bigodot}_{C^n}$ towards $\R^{n-k} \times \{0\}$ in terms of estimates for the mean curvature of distance tubes of $\R^{n-k} \times \{0\}$.

\begin{lemma} \emph{\textbf{(Bending Effects)}} \label{etl} \, We consider $u= f(\omega) \cdot r^\beta$, for  $\omega \in S_{C^k}:= \p B_1(0)
\cap C^k$,  $r(x):=dist(x, \R^{n-k} \times \{0\})$, for $x \in \R^{n-k} \times C^k$ and some positive $f \in C^{2,\alpha}(S_{C^k} \setminus \Sigma_{S_{C^k}},\R)$.\\

The mean curvature  of distance tubes $T_\varrho:=\R^{n-k} \times (\p B_\varrho(0)
\cap C^k) \subset C^n$ of $\R^{n-k} \times \{0\}$,  of radius $\varrho >0$, with radius measured relative $g_{C^n}$,  is given by
\begin{equation}\label{te}
u^{4/n-2} \cdot tr\, A_{T_\varrho}(u^{4/n-2} \cdot g_{C^n}) =  - \left( (k-1) +    2 \cdot (n-1) \cdot  \beta/(n-2)\right) \cdot \varrho^{-1}
\end{equation}
In particular,  for the two $scal >0$-metrics $g^+_{C^n}$ and  $g^{\bigodot}_{C^n}$ on $C^n \setminus \sigma_{C^n}$, we get
\begin{equation}\label{mmm}
tr\, A_{T_\varrho}(g^+_{C^n})  <0 \,\mm{ and } \, tr\, A_{T_\varrho}(g^{\bigodot}_{C^n})  >0.
\end{equation}
\end{lemma}
Our sign convention is that $\p B_1(0) \subset  (\R^n, g_{Eucl})$ has \emph{negative} mean curvature when viewed from $\R^n \setminus \overline{B_1(0)}$.\\

\textbf{Proof} \quad Recall that the second fundamental form $A_L(g)$ of a submanifold $L$ with respect to some metric $g$ transforms under conformal deformations $g
\ra u^{4/n-2} \cdot g$ according to the formula, cf. [Be],1.163,p.60:
 \begin{equation}\label{mmm}
 A_L(u^{4/n-2} \cdot g)(v,w) = A_L(g)(v,w) - \frac{2}{n-2}\cdot {\cal{N}} (\nabla u / u) \cdot g(v,w),
 \end{equation}
where $ {\cal{N}} (\nabla u / u)$ is the normal component of $\nabla u / u$ with respect to $L$. Since $C^k \subset \R^{k+1}$ is a cone, we have $tr\,A_{T_\varrho}( g ) = -(k-1)/\varrho$, since the $\R^{n-k}$-factor of $T_\varrho$ is totally geodesic.  In turn, for the trace of the second summand the value ${\cal{N}} (\nabla u /u)$, which is easily computed to $2/(n-2)  \cdot \beta/\varrho$, is multiplied  by $n-1$.\\

For the metrics $g^+_{C^n}$ and  $g^{\bigodot}_{C^n}$ on $C^n \setminus \sigma_{C^n}$ this computation reads as follows:
\begin{equation}\label{min}
\Psi^{4/(n-2)}_+[n,k]  \cdot tr\, A_{T_\varrho}(\Psi^{4/(n-2)}_+[n,k] \cdot g_{C^n}) =
\end{equation}
  \[- \left( (k-1) +    2 \cdot (n-1) \cdot  \alpha_+/(n-2)\right) \cdot \varrho^{-1}< -  (k-2)/2 \cdot \varrho^{-1}<0,\]
\begin{equation}\label{max}
 \Psi^{4/(n-2)}_{\bigodot}[n,k]  \cdot tr\, A_{T_\varrho}(\Psi^{4/(n-2)}_{\bigodot}[n,k] \cdot g_{C^n}) =
\end{equation}
  \[- \left( (k-1) +    2 \cdot (n-1) \cdot  \alpha_-/(n-2)\right) \cdot \varrho^{-1}> (k-2)/2 \cdot \varrho^{-1}>0.\] \qed

Thus we observe that $g^{\bigodot}_{C^n}$  shields $\R^{n-k} \times \{0\}$ against the
intrusion of area minimizers, whereas $g^+_{C^n}$ does not.\\

Now we merge $g^+_{C^n}$ and $g^{\bigodot}_{C^n}$ to a radial cut-off version of $g^{\bigodot}_{C^n}$ relative a base metric  $g^+_{C^n}$. That is, we insert a horizon into $g^+_{C^n}$ shielding $\R^{n-k} \times \{0\}$ while keeping $scal >0$.\\

The construction, starting from minimal metrics, may be regarded as a counterpart to the classical $scal >0$-preserving codimension $\ge 3$-surgeries on smooth $scal >0$-manifolds in [GL1]. Recall that in the classical case, we deform the complement of (a tube around) a submanifold $N$ of codimension $\ge 3$, so that it keeps $scal >0$ but we also get a totally geodesic boundary.\\

In both cases we need to sacrifice a portion of the $scal >0$ to spread the geometry towards the singular axes $\R^{n-k} \times \{0\}$  respectively towards $N$.\\

 We proceed in several elementary steps starting from the minimal metric.

\begin{proposition} \emph{\textbf{(Merging $g^+$ and  $g^{\bigodot}$)}}\label{el}\, For any cone $C^n= \R^{n-k} \times C^k$, $C^k \in {\cal{SH}}^{\R}_k$, and any pair $(R,r)$ with $R>r>0$ there is a metric $g^\Box_{C^n}$  with the following properties:
\begin{itemize}
\item $g^\Box_{C^n}$ is a $scal >0$-metric conformal to $g_{C^n}$.
\[g^\Box_{C^n} = \F^2_{R,r}(|y|) \cdot g_{C^n}, \mm{ for some } \F_{R,r} \in C^\infty(\R^{>0},\R^{>0}),  (x,y) \in C^n= \R^{n-k} \times C^k.\]
  \item $g^\Box_{C^n}= g^{\bigodot}_{C^n}$, on $\R^{n-k} \times B_R(0) \setminus B_r(0)$,
  \item  $g^\Box_{C^n}= g^+_{C^n}$, outside $\R^{n-k} \times B_{R^*}(0) \setminus B_{r^*}(0)$,
\end{itemize}
for some pair $(R^*,r^*)$ with $R^* > R$, $r>r^*>0$. $\mathbb{F}_{R,r}$ and $(R^*,r^*)$ depend only on $(R,r)$ and $n,k$, but not on $C^n \in {\cal{C}}_n$.
\end{proposition}

Here, and for later construction steps, we choose a cut-off function $\phi \in C^\infty(\R,[0,1])$ with $\phi \equiv 1$ on $\R^{\le 0}$, $\phi \equiv 0$  on $\R^{\ge 1}$.  \\

\textbf{Proof} \quad   We first consider a couple of conformal deformations \textbf{D1-D4} for one fixed cone $C$, each of them is radial symmetric. We specify the start and end radii later on. It is noteworthy, that the first deformation \textbf{D1} consumes some of the positivity of the scalar curvature, whereas the others do not.\\

\textbf{D1} \, For large $D >0$, we set $h_D(x):=\phi_D(\rho):=\phi(\rho-D)$, where $\rho=dist(x,\R^{n-k} \times \{0\})$ and get
\[( \Psi_+ +h_D\cdot \Psi_{\bigodot})^{(n+2)/(n-2)} \cdot scal\left( (\Psi_+ + h_D\cdot \Psi_{\bigodot})^{4/(n-2)} \cdot g_{C_0^n}\right)=\]
 {\small\[-(\Delta\Psi_+ +h_D \cdot \Delta\Psi_{\bigodot}) + \frac{(n-2)}{4(n-1)} \cdot scal(g_{C_0^n}) \cdot (\Psi_+ + h_D\cdot \Psi_{\bigodot}) - (\Delta h_D  \cdot \Psi_{\bigodot} + 2 \cdot \langle \nabla h_D, \nabla \Psi_{\bigodot} \rangle) \]}
Since the $\Psi_\pm$ are eigenfunctions for the same eigenvalue $\lambda >0$, this equals
 \[= \lambda \cdot \bp^2 \cdot (\Psi_+ + h_D\cdot \Psi_{\bigodot})  - (\Delta h_D  \cdot \Psi_{\bigodot} + 2 \cdot \langle \nabla h_D, \nabla \Psi_{\bigodot} \rangle) \]
\[= \psi_C(\omega) \cdot \Big(\lambda \cdot a^2 \cdot \rho^{-2} \cdot (r^{\alpha_+} + \phi_D\cdot r^{\alpha_-}) -  (\phi''_D + (k-1)/r \cdot \phi'_D) \cdot r^{\alpha_-} - 2 \cdot  \alpha_-  \cdot \phi'_D \cdot r^{\alpha_--1}  \Big).\]
For $D>0$ large enough, this expression is \emph{positive}, as is seen from the following details
\begin{itemize}
  \item $\alpha_+-\alpha_- \ge  2 \cdot \sqrt{2/3}  \cdot \frac{k-2}{2}>3$ and  $\lambda >0$,
  \item $a \ge A_n$, for some constant $A_n >0$, depending only on $n$, since $\delta_{\bp} \le L
        \cdot dist(\cdot,\Sigma)$, for some $L(\bp,n)>0$
  \item $|\phi'_D|, |\phi''_D| \le c(\phi)$, for some $c(\phi) >0$, independent of $D$.
\end{itemize}

\textbf{D2} \, For small $e >0$, we set $h_e(x):=\phi_e(\rho):=1-\phi(\rho/e -1)$ and get the following equation as in step D1:
\[(h_e\cdot \Psi_+ + \Psi_{\bigodot})^{(n+2)/(n-2)} \cdot scal\left( (h_e\cdot \Psi_+ + \Psi_{\bigodot})^{4/(n-2)} \cdot g_{C_0^n}\right)=\]
\[\psi_C(\omega) \cdot \Big(\lambda \cdot a^2 \cdot \rho^{-2} \cdot (\phi_e\cdot r^{\alpha_+} + r^{\alpha_-}) -  (\phi''_e + (k-1)/r \cdot \phi'_e) \cdot r^{\alpha_+} - 2 \cdot  \alpha_+  \cdot \phi'_e \cdot r^{\alpha_+-1}  \Big). \]
This is \emph{positive}, provided $e>0$ is small enough. Here, we use
 \begin{itemize}
   \item $\alpha_+-\alpha_- \ge  2 \cdot \sqrt{2/3}  \cdot \frac{k-2}{2}>3$ and $\lambda >0$
   \item $a \ge A_n$, for  $A_n >0$, depending only on $n$
   \item $|\phi'_e| \cdot e, |\phi''_e| \cdot e^2\le c^*(\phi)$, for some $c^*(\phi) >0$, independent of $e$.
 \end{itemize}

\textbf{D3} \, For any $\gamma \ge 1$, there is a unique $\rho_\gamma >0$, so that  $\gamma \cdot \rho^{\alpha_+}_\gamma =\rho^{\alpha_-}_\gamma$. For $\rho < \rho_\gamma$, we have
 $\gamma \cdot \rho^{\alpha_+} <\rho^{\alpha_-}$ and, for $\rho > \rho_\gamma$, we have $\gamma \cdot \rho^{\alpha_+} >\rho^{\alpha_-}$. Also we observe, that $\rho_\gamma \ra 0$, for $\gamma \ra \infty$.\\

Therefore, for any $\ve \in (0,\rho_\gamma/10)$, we can find a smooth interpolating function $F_{\gamma,\ve}>0$ so that $F_{\gamma,\ve}(\rho)=\gamma \cdot \rho^{\alpha_+}$ for $\rho <\rho_\gamma - \ve$ and
$F_{\gamma,\ve}(\rho)= \rho^{\alpha_-}$ for $\rho>\rho_\gamma +\ve$, with \[(\psi(\omega) \cdot F_{\gamma,\ve})^{(n+2)/(n-2)} \cdot scal\left((\psi(\omega) \cdot F_{\gamma,\ve})^{4/(n-2)} \cdot g_{C_0^n}\right)= \]
\[L_C (\psi_C(\omega) \cdot F_{\gamma,\ve}) = -\psi_C(\omega) \cdot \left(F_{\gamma,\ve}''  + (k-1)/\rho \cdot F_{\gamma,\ve}'\right) + L_{C,\lambda}^\times \psi(\omega) \cdot F_{\gamma,\ve}= \] \[\psi_C(\omega) \cdot \left(- (F_{\gamma,\ve}''  + (k-1)/\rho \cdot F_{\gamma,\ve}') + \mu_{C,L_{C,\lambda}^\times} \cdot F_{\gamma,\ve}\right)>0\]

The latter positivity holds for both component functions $\gamma \cdot \rho^{\alpha_+}$ and $\rho^{\alpha_-}$. An appropriated smoothing $F_{\gamma,\ve}$ of the piecewise smooth function $\min\{\gamma \cdot \rho^{\alpha_+},\rho^{\alpha_-}\}$ can be accomplished employing convolutions with smoothing kernels and subsequent use of cut-off functions placed towards the interval ends of $(\rho_\gamma-\ve,\rho_\gamma+\ve)$.\\

\textbf{D4} \, For  $d \in \Z^{>0}$, we set  $h_d(x):=\phi_d(\rho):= \phi(2^d \cdot \rho -1)$ and we notice for $|\eta| \le \zeta$, for some sufficiently small $\zeta >0$:
\[((1+\eta \cdot h_d)\cdot\Psi_+)^{(n+2)/(n-2)} \cdot scal\left(((1+\eta \cdot h_d)\cdot\Psi_+)^{4/(n-2)} \cdot g_{C_0^n}\right)=\]
\[\psi(\omega) \cdot \Big(\lambda \cdot a^2 \cdot \rho^{-2} \cdot (1+\eta \cdot \phi_d) \cdot \rho^{\alpha_+}  -  \eta \cdot (\phi''_d + (k-1)/\rho \cdot \phi'_d) \cdot \rho^{\alpha_+} - 2 \cdot \eta \cdot \alpha_+  \cdot \phi'_d \cdot \rho^{\alpha_+-1}  \Big) >0, \]
The point is, that for small $\zeta$, the factor $1+\eta \cdot h_d$ merely causes a small perturbation in $C^3$-norm. Also, this $\zeta >0$ can be chosen independently of $d$ since $\Psi_+^{4/(n-2)} \cdot g_{C_0^n}$ is scaling invariant, and, iteratively, because $supp |\nabla h_{d_1}| \cap supp |\nabla h_{d_2}| \v$, for $d_1 \neq d_2$, we infer for any $|\eta| \le \zeta$, $1 \le i < j$:
 \[scal\Big(\Big(\prod_{d=i,...j}(1+\eta \cdot h_d)\cdot\Psi_+\Big)^{4/(n-2)} \cdot g_{C_0^n}\Big) >0.\]

\textbf{Conclusion} \, We seemlessly compose the four deformations \textbf{Dm}, $m=1,..,4$ on disjoint shells $\R^{n-k} \times B_{\rho_m}(0) \setminus B_{\rho_{m-1}}(0)$,  for radii $\rho_0 > \rho_1 > \rho_2 > \rho_3 > \rho_4 >0$. That is we choose the radii stepwise: $\rho_0$ and $\rho_1$ for D1, then we continue between radii $\rho_1$ and $\rho_2$ with \textbf{D2} and similarly we append \textbf{D3} and finally we reach \textbf{D4}. \\

\textbf{D1} and \textbf{D2} are the essential bendings to get a repelling shell around  $\R^{n-k} \times \{0\}$. For any $\Gamma \gg 1$ we can choose the inner radius $\rho_2$ of D2 small enough, relative to $\rho_1$, so that any area minimizing hypersurface passing this shell will have area $>\Gamma$.\\

 The two deformations \textbf{D3} and \textbf{D4} are necessary to further localize the entire deformation to a compact subset, also in $\R^{n-k}$-directions. Then we need to ensure that the metric on $\R^{n-k} \times B_{\rho_4}(0)$ is again the original minimal metric.\\

Finally, we observe that all constants and estimates depend only on $n$, $k$ and the chosen function $h$. Thus the claim actually holds for the same radii for any such cone $C^n \in {\cal{SC}}_n$. \qed

We note that the metrics $g^\Box_{C^n}$ are again invariant under $\R^{n-k}$-translations. Thus from \ref{etl}  we infer

\begin{corollary}  \emph{\textbf{(Horizons along $\{0\} \times C^k$)}} \label{lhi} \, For any cone $C^n= \R^{n-k} \times C^k$, $C^k \in {\cal{SH}}^{\R}_k$ equipped with a metric $g^\Box_{C^n}$ of \ref{el},  there is an area minimizing hypersurface \[F^{n-1} =\R^{n-k} \times F^{k-1}_{C^k}, \mm{ where } F^{k-1}_{C^k} \subset C^k\]  is an  area minimizing hypersurface in $C^k$ equipped with $g^\Box_{C^n}|_{\{0\} \times C^k}$ with
 \begin{equation}\label{bou}
 B_{R}(0) \cap C^k  \subset  F^{k-1}_{C^k} \subset B_{R^*}(0) \cap C^k.
 \end{equation}
\end{corollary}

\textbf{Proof} \quad The estimates \ref{etl} show that the separating hypersurfaces $\R^{n-k} \times \p B_{R}(0)$ and $\R^{n-k} \times \p B_{R^*}(0)$ have positive resp. negative mean curvature. Thus an area minimizing sequence $H_i \subset \R^{n-k} \times  B_{R^*}(0) \setminus B_{R}(0)$ of $\R^{n-k}$-translation invariant hypersurfaces, homologous to $\R^{n-k} \times \p B_{R}(0)$, can be replaced by another one supported within a fixed positive distance of the two repelling borders and, thus, the limit $F^{k-1}_{C^k}$ has the asserted properties. \qed

\textbf{Localization} \, The support of each of the deformations we discussed in  the last section is the entire \emph{non-compact} set $\R^{n-k} \times B_{R^*}(0) \setminus B_{r^*}(0)$. Now we shrink the support to a compact subset so that it fits into a multi-pipe.  For this we add a cut-off process along the $\R^{n-k}$-axis, either along balls or $n-k$-dimensional rectangles
\[Q_L (0):= (-l_1,l_1) \times...\times (-l_{n-k},l_{n-k}), \mm{ for any vector } L:=  (l_1,...,l_{n-k}) \mm{ with } l_i >0.\]
Rectangles around other basepoints $x \in  \R^{n-k}$ are denoted $Q_L(x)$. Also, for $k \in \R^{>0}$,  we write $L>k$, when $l_i>k$, for all $i=1,..n-k$.\\

This localization is another instance where we need to sacrifice some of the $scal >0$. But different from the latter constructions the estimates for the lower bound of $scal$ weaken when we approach those singularities in $C$ perpendicular to  $\R^{n-k} \times \{0\}$. But for any $\rho>0$ and a sufficiently slow transition from $g^\Box$ to $g^+$ we can ensure such estimates for $\P(0,\rho) \subset C^k$.\\

This matches our $\bp$-multi-pipes setup: we use the singular threshold \ref{strs}  to see that the directions $\P^*(0,\rho)$ are iteratively surrounded by spokes until
we reach cones singular only in $0$.

\begin{lemma} \emph{\textbf{(Localization of $g^\Box$)}}\label{el2}\, For any cone $C^n= \R^{n-k} \times C^k$, $C^k \in {\cal{SC}}^{\R}_k$, and any quadruple $(R,r,\eta,L)$ with $R>r>0$, $\vartheta>0$ and $L>0$, there is a metric $g^\Delta_{C^n}$  with the following properties:
\begin{enumerate}
\item $g^\Delta_{C^n}$ is conformal to $g_{C^n}$.
\item $scal(g^\Delta_{C^n}) >0$ outside $Q_{L^*}(0) \setminus  Q_L(0) \times \P^*(0,\vartheta) \cap B_{R^*}(0) \setminus B_{r^*}(0)$
  \item $g^\Delta_{C^n}= g^{\bigodot}_{C^n}$ on $\P(0,\vartheta) \cap Q_L(0) \times B_R(0) \setminus B_r(0)$
  \item  $g^\Delta_{C^n}= g^+_{C^n}$ outside $\P(0,\vartheta/2) \cap Q_{L^*}(0) \times  B_{R^*}(0) \setminus B_{r^*}(0)$.
\end{enumerate}
for some triple $(R^*,r^*,L^*)$ with $R^* > R$,$r>r^*>0$ and $L^*>L$, that is, $l^*_i>l_i$,$i=1,..,n-k$, with $(R^*,r^*)$ depending only on $(R,r)$ and the dimension, and $L^*-L$, that is,  $l^*_i-l_i$  depending only on $\vartheta$ and on the dimension.
\end{lemma}

\textbf{Proof} \quad We set $h_{L,D}(x):= \phi((dist^*(x,Q_L(0))/D)$, where $dist^*$ is some smoothing of $dist$ we choose very close to $dist$ in $C^0$-topology as vanishing on $Q_L(0)$. As in D1 of \ref{el} we have
 \[( h_{L,D} \cdot \Psi_+ +(1-h_{L,D})\cdot \F_{R,r}(|y|) \cdot g_{C^n})^{\frac{n+2}{n-2}} \cdot \]
\[...scal\left(( h_{L,D} \cdot \Psi_+ +(1-h_{L,D})\cdot \F_{R,r}(|y|))^{4/(n-2)} \cdot g_{C^n}\right)\]
 \[= \lambda \cdot \bp^2 \cdot ( h_{L,D} \cdot \Psi_+ +(1-h_{L,D})\cdot \F_{R,r}(|y|))...\]
\[-\big(\Delta h_{L,D}  \cdot \Psi_+ -2 \cdot \langle \nabla h_{L,D}, \nabla \Psi_+ \rangle \big) - \big(\Delta  (1-h_{L,D})  \cdot \F_{R,r}(|y|)  + 2 \cdot \langle \nabla  (1-h_{L,D}), \nabla \F_{R,r}(|y|)  \rangle \big) \]
 \[= \lambda \cdot \bp^2 \cdot ( h_{L,D} \cdot \Psi_+ +(1-h_{L,D})\cdot \F_{R,r}(|y|) )...\]
\[ -  \big(\Delta h_{L,D}  \cdot (\Psi_+ - \F_{R,r}(|y|) ) + 2 \cdot \langle \nabla h_{L,D}, \nabla (\Psi_+ - \F_{R,r}(|y|) )  \rangle \big)\]

In the last expression $|\Psi_+ - \F_{R,r}(|y|)|$  and $|\nabla (\Psi_+ - \F_{R,r}(|y|))|$ remain bounded on $\P(0,\vartheta) \times \R^{n-k} \subset C^k \times  \R^{n-k}$.\\

For any $\ve >0$ there is some $D$ large enough, so that $|\Delta h_{L,D}|, |\nabla h_{L,D}| < \ve$, whereas $\lambda \cdot \bp^2 \cdot ( h_{L,D} \cdot \Psi_+ +(1-h_{L,D})\cdot \F_{R,r}(|y|) )$ remains positively bounded on $\R^{n-k} \times B_{R^*}(0) \setminus B_{r^*}(0)$.  From these observations one readily derives the assertions. Note that (iv) is a simple smooth cut-off without $scal$-control.
\qed

\subsubsection{Horizon Metrics and Inductive $\D$-Transfer} \label{lll}
\bigskip

Our plan is to transfer the deformations from $g^+_{C^n}$ to $g^\Delta_{C^n}$ to all pipes of multi-pipes $\sp$. To this end we notice that this deformation can also be accomplished adding the compactly supported and, as is readily checked, positive semi-definite tensor $g^\Delta_{C^n}-g^+_{C^n}$:

\begin{definition}\label{lhp} \, We call $\L_{C^n}:=g^\Delta_{C^n}-g^+_{C^n}$, a \textbf{local horizon pseudometric}.
\end{definition}

We will see that $g^\Delta$ has a localized repelling effect from the wedge of $g^{\bigodot}_{C^n}$ it contains. That is, it shields some portion of $\R^{n-k} \times \{0\}$ from
the access of area minimizers.\\

Towards the use for $\sp$, we notice that the spokes are not exactly Riemannian product but tapered product. Thus we also consider transfers using tapered versions of $\D$-maps.\\

\textbf{Tapered $\D$-Maps} \, When $C^n= \R^{n-k} \times C^k$ is an iterated tangent cone that locally approximates a spoke in the hub and spokes decomposition of Def.\ref{st2}. Then the $\D$-map is locally almost an isometry and we can use it to transfer $\L_{C^n}$. Globally, the spoke are tapered products but we can still use the $\D$-maps to globally transfer $\L_{C^n}$ to the tapered products using tapered versions of $\D$-maps, we explain now. This allows us to treat the shielding aspect of spokes as for proper Riemannian products.\\

We consider $k=n-1$, the cases $k <n-1$ follow similarly using some rather obvious, but notionally voluminous, induction we omit here.\\

 Thus let $C^n= \R^1 \times C^{n-1}$ be an (almost) tangent cone for some $p \in \p B_1(0) \cap C^*$, $C^* \in {\cal{SC}}_n$ and, therefore, for any point $t \cdot p$,  $t >0$.\\

 We modify the $\D$-map to a \textbf{tapered $\D$-map} $\D^\tau_{tap}$ defined along the ray $\R^{>0} \cdot p \subset C$: $C^{n-1}$ is a tangent cone of $S_{C^*}$ at $p \in \p B_1(0) \cap C^*$ and we choose the corresponding $\D$-map $\D_{C^{n-1} \ra \tau \cdot S_{C^*}}$, for some large $\tau \gg 1$:
\[\D^\tau_{C^{n-1}}: B_R(0)\setminus B_r (0) \cap  \P(0,\delta) \ra \tau \cdot  B_{R/ \tau}(p) \setminus B_{r/\tau} (p) \cap  \P_{\tau \cdot S_{C^*}}(p,\delta)\]

We extend this diffeomorphism radially to a diffeomorphic map $\D_{tap}$  from (subsets of) $\R^1 \times C^{n-1}$ to $C^*$
\begin{equation}\label{tap}
\D^\tau_{tap}: \R^{>0} \times B_R(0)\setminus B_r (0) \cap  \P(0,\delta) \ra C^*,\mm{ with }  \D^\tau_{tap}(t,z):= t \cdot \D^\tau_{C^{n-1}}(z)
\end{equation}

One readily checks that for increasing $\tau$ and suitable rescaling $\D^\tau_{tap}$ locally turns into an isometry:

 \begin{lemma} \emph{\textbf{(Nearly Product Shape of $\D^\tau_{tap}$)}}\label{tst}\, For any  $\delta >0$, $d>1$ and $\eta >0$  we can choose a large $\tau_{d,\eta} >1$ so that, for any $\tau \ge \max\{\tau_{d,\eta} ,10\}$ and any $a > 0$, after scaling by $1/a^2$
\begin{equation}\label{tap}
\D^\tau_{tap}: [a/d, d \cdot a] \times  B_R(0)\setminus B_r (0)  \cap  \P(0,\delta) \ra C^*
\end{equation}
is $\eta >0$-close in $C^3$-topology to an isometry onto its image.
\end{lemma}

The latter statements means that the $\D^\tau_{tap}$-pull-back of $g_{C^*}$ is $\eta >0$-close to the inherent product metric in $C^3$-topology.\\

This allows us to transfer $g^\Delta$ to $C^*$ with largely unchanged properties. This is also the elementary type of transfer of local horizons between nodes with increasing height in separator trees:
\begin{definition} \emph{\textbf{(Tapered Pseudometrics)}} \label{merg} \,  For any cone $C^*  \in {\cal{SC}}^{\R}_n$  and an (almost) tangent cone $C^n =  \R \times C^{n-1}$, we define a tapered version of the local horizon pseudometrics
 \begin{equation}\label{th}
 g_{\L,Y,C^*}:= (\D^\tau_{tap})_*(Y_*(\L_{\R^1 \times C^{n-1}}))
 \end{equation}
where $(\D^\tau_{tap})_*(\cdot)$ denotes the push-forward of the tensors and $Y$ is a transformation of $C$ defined by $Y(x,y):= \iota \cdot (x,y) + (\nu,0),$ for some $\iota>0$ and $\nu \in \R$.\\
\end{definition}

\subsubsection{Trees of Equidistributed Measures} \label{trr}
\bigskip

Now we define the main deformation of $g^+_H$ to a metric $g_{\T}$ on $H \setminus \Sigma$,  adding tapered pseudometrics to all elements of $\T$ and transferring them along the branches of $\T$ until we reach the root node. We also explain how the tree has to be chosen, to ensure the desired shielding and scalar curvature properties.\\

 For the sake of clarity we notice that the subsequent discussion involves the \emph{intermediate} use of area minimizers intersecting the singular set, cf.[AK], and the recent survey [H], for appropriated technical frameworks. However, the resulting global horizon will
 be entirely supported in the smooth manifold $H \setminus \Sigma$.\\

\textbf{Tightness of Open Ends} \, In our $\bp$-multi-pipe setup each pipe comes with open ends. The risk is that some area minimizers in $H$ could exploit this weakness to approach $\Sigma$ from such an apparently unprotected side of one of these pipes.\\

But when we choose trees $\T$ with small separators we find that even these open ends remain tight. We illustrate this simple recipe:\\

\emph{An area minimizing hypersurface spanned over a narrow hollow cylinder $B_1^{n-k}(0) \times \p B_r^k(0) \subset \R^n$, for   $r >0$ sufficiently small,  does not get deeply into the cylinder $B_1^{n-k}(0) \times B_r^k(0)$, but it remains closely to the boundary $\p (B_1^{n-k}(0) \times B_r^k(0))$.}\\

This Euclidean sample case can easily be checked applying elementary symmetrization arguments from the growth comparison  $Vol_{n-1}(\p B_1^{n-k}(0) \times B_r^k(0)) \sim r^k$, whereas we have $Vol_{n-1}(B_1^{n-k}(0) \times \p B_r^k(0)) \sim r^{k-1}$. \\

For the counterpart in our context, we first consider the proper product $C^n= \R^{n-k} \times C^k$ and use similar area comparison arguments to show, cf.Ch.\ref{geom} for details:

 \begin{lemma} \emph{\textbf{(Local Horizons)}}\label{haee}\, For any $R>0$, $L\ge 2$ and $\delta>0$, there is a radius $r(R,L,\delta)>0$, $R \gg r$, so that there is a    connected area minimizing $F \subset Q^{n-k}(L^*) \times B_{R^*}(0)$ relative $g^\Delta_{C^n}$ with
\begin{itemize}
  \item $F$ bounds an open set ${\cal{E}}_{C^n}$ with $(Q^{n-k}(L-\delta)  \times B_{\frac{2 \cdot R+ r}{3}}(0)) \subset {\cal{E}}_{C^n}$,
    \item $F \cap \p (Q^{n-k}(L^*) \times B_{R^*}(0)) \v$.
    \item For $L\ra \infty $, $F=F(L)$ converges compactly to a product shape hypersurface $F^{n-1} =\R^{n-k} \times F^{k-1}_{C^k}$ of \ref{lhi}.
    \end{itemize}
\end{lemma}

The actually tapered elements in $\T$ receive their  pseudometrics from an $\D_{tap}$-transfer to $C^*$. We choose the parameters $L$, for chosen $R$ and $r$,  in the localization construction  \ref{el2} so that $L/R^*$ is also very large. Since $L^*-L$ remains constant this also means we may assume $(L^*-L) \ll L$. Then we also get, in \ref{geom}:

\begin{lemma}  \emph{\textbf{(Tapered Local Horizons)}}  \label{l3} \, For any  radii $\rho_2 > \rho_1 >0$, $\zeta >0$ and for any $\delta \in (0,(\rho_2 - \rho_1)/100)$, there is a $\tau \gg 1$ and some $\zeta^*(\zeta) \ge 10 \cdot \zeta >0$ with $\zeta^*(\zeta) \ra 0$, for $\zeta \ra 0$, and some transformation $Y$ so that for $g_{\L,Y,C^*}+ g^+_{C^*}$
\[g_{\L,Y,C^*}= 0\mm{ outside } C(B_{\zeta^*}(p)) \cap B_{\rho_2 + \delta}(0) \setminus B_{\rho_1 - \delta}(0)\]
and there is a connected  connected area minimizing $F \subset C(B_{\zeta^*}(p)) \cap B_{\rho_2 + \delta}(0) \setminus B_{\rho_1 - \delta}(0)$ relative $g_{\L,Y,C^*}+ g^+_{C^*}$ so that
\begin{itemize}
  \item $F$ bounds an open set ${\cal{E}}^{tap}_{C^n}$ with $C(B_{\zeta}(p)) \cap B_{\rho_2 - \delta}(0) \setminus B_{\rho_1 + \delta}(0) \subset {\cal{E}}^{tap}_{C^n}$,
    \item $F \cap \p (C(B_{\zeta^*}(p)) \cap B_{\rho_2 + \delta}(0) \setminus B_{\rho_1 - \delta}(0) ) \v$.
    \item For $\tau\ra \infty $ and after rescaling by $\tau^2$, the $\D^\tau_{tap}$-pull-back of $F=F(\tau)$ converges compactly  to an area minimizer $F^{n-1} =\R^{n-k} \times F^{k-1}_{C^k}$ as in \ref{lhi}.\\
    \end{itemize}
\end{lemma}

\textbf{The Horizon Metrics $g_{\T}$} \, With these $g_{\L,Y,C^*}$ we can already write our central metric $g_{\T}$ on $H$, where we compose the iterated $\D^\tau_{tap}$-maps to a map $\D_{\CCC}$ for the upstream move that pushes each of the $\L_{\CCC}$ forward to its designated place on $H$, where we fill $\L_{\CCC}$ into one of the elements of the given multi-pipe. \\

This means, the support belongs to the element and, under the standard assumption of a very high accuracy and small separator, we may assume that the parameters for the localization $L$ and $L^*$ satisfy $(L^*-L)/L \ll 1$ and the ends of $\D_{tap}$-transfer of the rectangle $Q_{L^*}$ reach to ends of the elements they fill. Then \ref{l3} asserts the shielding effect
for any single element of our multi-pipe (before it is transferred to $H$).

\begin{definition} \emph{\textbf{(Horizon Metrics $g_{\T}$)}} \label{merg} \,  For any  separator tree $\T$ and  the almost tangent cones $\CC^n = C^k \times \R^{n-k} \in \T$, each contributing its local horizon pseudometric $\L_{\CCC}$, we define the \textbf{horizon metric}
\[g_{\T}:= \sum_{\CCC \in \T}(\D_{\CCC})_*(\L_{\CCC})+  g^+_H.\]
\end{definition}

For suitable separator trees $\T$, this becomes a raw version of the geometries on $H$ we described in our main theorems in the introduction. (We derive them in the next section Ch.\ref{smooth}.)

 \begin{proposition} \emph{\textbf{(Superposition of Local Horizons)}} \label{merg2} \, For any neighborhood $V$ of $\Sigma_H$, $\eta >0$ there  a separator tree $\T_{V,\eta}$ and a  horizon metric $g_{\T}$, so that there is a collection of hypersurfaces $F_{\CCC}$ with non-negative mean curvature relative $g_{\T}$  bounding open sets ${\cal{O}}_{\CCC}  \subset V$ with
\[\textstyle \Sigma_H \subset \bigcup_{\CCC \in \T} {\cal{O}}_{\CCC}\, \mm{ and }\, scal(g_{\T}) >0  \mm{ on } H \setminus \bigcup_{\CCC \in \T} {\cal{O}}_{\CCC}.\]
so that the flat norm distance of ${\cal{O}}_{\CCC}$ to $\D_{\CCC}({\cal{E}}^{tap}_{C^n})$ is less than $\eta$, where we only map the portion of ${\cal{E}}^{tap}_{C^n}$  in the domain of $\D_{\CCC}$.
 \end{proposition}

\textbf{Outline of Proof} \, We refer to Ch.\ref{wm} and \ref{smoo} for technical details. We start with area minimizers spanned over each of the individually deformed elements of $\sp$. That is, we consider $(\D_{\CCC})_*(\L_{\CCC})+  g^+_H$ for each $\CC \in \T$ separately. Here we use  \ref{l3}, for the model cone geometry $\CC$, and the skin uniformity of $H$ to show that a counterpart to \ref{l3} also holds, after $\D_{\CCC}$-transfer, within $(H,(\D_{\CCC})_*(\L_{\CCC})+  g^+_H)$.\\

 Then we argue that these individual area minimizers relative $(\D_{\CCC})_*(\L_{\CCC})+  g^+_H$ are not substantially damaged when we stepwise turn to $\sum_{\CCC \in \T}(\D_{\CCC})_*(\L_{\CCC})+  g^+_H$, following the branches of $\T$.\\

To this end, we choose, in \ref{bow3} and \ref{set}, a separator tree $\T$ for prescribed $\widehat{\Theta}^\bullet,\widehat{\ve}^\bullet,\widehat{s}^\bullet$, where we use very small self-similarity factors, accuracy values and separators but with  $\pi^\bullet_{n+1}[\kappa,1]/\pi^\bullet_{n+1}[\kappa+1,0]\gg 1$, $\kappa=1,...,c(n+1)$.
\begin{itemize}
  \item For some high accuracy $\widehat{\ve}^\bullet, \ll 1$ and for \emph{small} but \emph{broad} $\widehat{\Theta}^\bullet$, cf.\ref{sse},  the tapered products locally becomes almost isometric to products and narrow, but also we separate the sizes of pipes in the $c(n+1)$-subfamilies in one given layer. This makes the interference of the shielding effects arbitrarily small.
  \item An additionally small self-similarity factor $\widehat{s}^\bullet \ll 1$ is used to make the effect of deformations along different layers virtually independent.
\end{itemize}
For these choices, \ref{set} gives us the compatible set of parameters with the asserted properties of $g_{\T}$.\qed

\begin{remark}\textbf{(Equidistributed Measures)}\, We recall that $\Psi_{\bigodot}$ is the solution of $(L_{C^k,n})_\lambda \,w=0$ associated to the Lebesgue measure $\lambda_{\,\R^{n-k} \times \{0\}}$ supported on $\R^{n-k} \times \{0\} \subset \R^{n-k} \times C^k$, the canonical \emph{equidistributed measure} along the symmetry axes of $C^n$.\\

  The pseudometrics $\L_C$, for  the various (almost) tangent cones of the chosen separator tree $\T$, are localized versions of  the $g^{\bigodot}_{C^n}:=\Psi_{\bigodot}^{4/(n-2)}[n,k] \cdot (g_{\R^{n-k}} + g_{C^k})$, relative the underlying minimal metric $g^+_{C^n}$.  Thus the metric $g_{\T}$ imitates the bending effect of the metric we would get from an assembly of piecewise equidistributed measures.\qed
 \end{remark}

\subsubsection{Regular End Structures}\label{smooth}
\bigskip

Here we derive the main Theorems we formulated in the introduction from \ref{merg2} from some fine-tuning of the boundary structures. Thus, given a singular area minimizer $H^n \subset M^{n+1}$ within some smooth compact manifold $M$ with $scal \ge 0$, we start with Theorem 1:

\begin{proposition} \emph{\textbf{(Regular Open Ends)}} \label{merg3} \, For any  $\ve > 0$, there is a smoothly bounded neighborhood $W_{\ve}$ of the singular set $\Sigma \subset H$ within an $\ve$-neighborhood $U_{\ve}$ of
$\Sigma$ such that $H \setminus W_{\ve}$ admits a smooth metric $g^+_{\ve}$  so that
\begin{enumerate}
\item \quad $g^+_{\ve}$ is conformal to $g_H$, \item  \quad $scal (g^+_{\ve}) > 0$ on $H \setminus W_{\ve}$, \item  \quad $\p W_{\ve}$ has positive mean curvature.\\
\end{enumerate}
 \end{proposition}

\textbf{Proof} \,  We start from the $scal>0$-metric $g_H^+=\Phi(H)^{4/(n-2)} \cdot g_H$ on $H$.  $\Phi(H)(\cdot) =G(q_H,\cdot)>0$ is singular in  $q_H \in H \setminus \Sigma$, but this does not cause any issues:\\

 The growth of $\Phi(H)$ is of order $r^{-(n-2)}$, with $r=dist(x,q_H)$, for $x$ near $q_H$ and this means that $\Phi(H)^{4/(n-2)} \cdot g_H$ is a $scal >0$-metric with some additional singularity in $q_H$. More precisely, we recall that for $\rho>0$ sufficiently small there exists a
unique Green’s function $G_\rho(q_H,\cdot)>0$ on $B_\rho(q_H) \setminus \{q_H\}$ with
\[L_{H,\lambda} \, G_\rho=0, \, G_\rho|_{\p B_\rho(q_H)} =0\, \mm{ and }\,  r^{n-2} \cdot G_\rho(x) \ra 1, \mm{ for } x \ra q_H.\]
and so that compactly, outside $q_H$, $G_\rho \ra \Phi(H)$, for $\rho \ra diam(H,g_H)$.\\

Thus, $B_\rho(q_H) \setminus \{q_H\}$ is deformed into an asymptotic flat end, where $q_H$ becomes the point at infinity. For sufficiently small $\eta \in (0,\rho)$ we observe that $\p B_\eta(q_H)$ is a smooth positive mean curvature hypersurface in  $\Phi(H)^{4/(n-2)} \cdot g_H$. An elementary modification allows us to conformally close this end to a compact half sphere keeping $scal >0$. \\

Now we consider the deformation of $g_H^+$ to $g_{\T}$.  From \ref{merg2} we get a hypersurface  $F$ that bounds a small neighborhood $U(F)=\bigcup_{\CCC \in \T} {\cal{O}}_{\CCC}$ of $\Sigma$. We can choose it so that $U(F) \subset U_{\ve}$. Thus mean curvature of $F$ in regular points vanishes and $F$ is  \emph{quasiregular} in the sense of [G2],Ch.3.2. Thus, it can be smoothed. In [G2], Gromov explains his remarkably simple approach, cf. [G2],p.950. The outcome is that $F$ can be $\eta$-Hausdorff-approximated by some homologous \emph{smooth} compact hypersurface
$F_\eta \subset H \setminus \Sigma$ with \emph{strictly positive mean curvature}, for any $\eta >0$.\\

Formally, this smoothing result of [G2] applies to singular  hypersurfaces which already have positive mean curvature. To adapt this situation to this ours, we first $C^3$-perturb $g_{\T}$, so that each of the $\p {\cal{O}}_{\CCC}$ becomes a locally unique area minimizer. Then it is easy to check that $\p {\cal{O}}_{\CCC}$ can be Hausdorff-approximated by still singular hypersurfaces $\p {\cal{O}}^*_{\CCC}$ that bounds some open set ${\cal{O}}^*_{\CCC}$, so that $\p {\cal{O}}^*_{\CCC}$ has mean curvature $\ve >0$. It is an $\ve$-bubble, cf.\ref{muu} below. Then we apply the smoothing operation to the boundary of  $\bigcup_{\CCC \in \T} {\cal{O}}^*_{\CCC}$.  \qed

 \begin{remark} ($\mathbf{\mu}$\textbf{-bubbles}) \label{muu} \,  For any, for simplicity orientable, Riemannian manifold $(M^n,g)$, $\mu\in C^0(M^n, \R)$,
we define the $\mu$-volume $Vol_\mu(W)$, for any rectifiable hypersurface $W \subset M^n$ bounding an open set $W^+ \subset M^n$. Upon choosing some orientation, we think of $W^+$ as the inner of $W$. Then we consider the functional
\[Area_\mu(W) := Area(W) - Vol_\mu(W)= Area(W) - \int_{W^+} \mu(z) \: dVol\]
This functional can be extended to \emph{Caccioppoli sets} and, in general, \emph{rectifiable integer multiplicity currents}, cf.[L1], appendix.\\

We call local minimizers $H_\mu$ of $Area_\mu$, the $\mu$\emph{-bubbles}.  These hypersurfaces are \emph{almost minimizers}. Therefore, they share the regularity theory with area minimizers e.g. their singular sets are still of codimension $7$ and tangent cones in singular points exists and are area minimizing cones, cf.[L1], appendix.\\

  The geometric detail, we are interested in, is that for any regular point $ x \in H_\mu$ the \emph{mean
curvature} $ tr \: A_{H_\mu}(x) $ \emph{equals} $ \mu(x)$, cf.[G1],Ch.5 for additional information. \qed\end{remark}

Theorem 2 is a geometric interpretation and application of \ref{merg3}.

\begin{corollary}  \emph{\textbf{(Smooth Ambience I)}}  \quad For any  closed smooth hypersurface $L \subset H \setminus \Sigma$ there is some $\ve_L > 0$ such that for
any $\ve \in (0,\ve_L)$ we get an area minimizing hypersurface $L_{\ve} \subset (H \setminus W_{\ve},g^+_{\ve})$ so that $$L_{\ve} \mm{ is  homologous to } L \mm{ and }
L_{\ve}  \cap \p W_{\ve} = \emptyset.$$ Thus $L_{\ve}$ is an area minimizer within a smooth $scal >0$-curved ambient space.
\end{corollary}

\textbf{Proof} \quad Since $L$ is compact with $L \cap \Sigma \v$, we can choose $\ve_L > 0$ small enough so that $L \cap \overline{U_{\ve}} \v$. Then $H \setminus W_{\ve}$ is strictly mean-convex and for an area minimizing sequence of hypersurfaces $L_n \subset H \setminus W_{\ve}$ homologous to $L$, there is an $\eta >0$ and a second such sequence $L^*_n  \subset H \setminus W_{\ve}$ so that:
 \[ L^*_n \cap U_\eta(\p W_{\ve}) \v \, \mm{ and }\, Area(L^*_n) \le Area(L_n).\]
 But then there a subsequence converging to an area minimizing hypersurface limit $L_\ve \subset H \setminus W_{\ve}$ homologous to $L$ with
$L^*_\ve \cap U_\eta(\p W_{\ve}) \v$. \qed

Theorems 3 and 4 are extensions of these results to include cases where all representatives of a homology classes $H$ hit $\Sigma$. To derive them we deform $H$ near $\p W_\ve$  another time, but in a non-conformal way, to assemble smooth approximations of singular minimizers.

\begin{corollary}  \emph{\textbf{(Regular Closed Doublings)}} \quad For any  $\ve > 0$,  there is a smoothly bounded neighborhood $V_{\ve}$ of $\Sigma$ contained in an
$\ve$-neighborhood $U_{\ve}$ of $\Sigma$ such that the doubling $H^D_{\ve}:= H \setminus  V_{\ve}  \cup_\sim H\setminus  V_{\ve} $, where $\sim$ means gluing along $\p
V_{\ve}$ admits a smooth metric $g^*_{\ve}$  with $g^*_{\ve} \equiv g_{ini}$ on $H \setminus U_{\ve}$ so that
\begin{enumerate}
\item$scal (g^*_{\ve}) > 0$   on $H \setminus V_{\ve}$, \item$g^*_{\ve}$ is conformal to the induced metric on $H \setminus V_{\ve}\subset M$, except for a kind of
    generalized warped product deformation (arbitrarily) near $\p V_{\ve}$.
    \item $H^D_{\ve}$ is mirror symmetric along $\p V_{\ve}$ and, hence, $\p V_{\ve}$ is totally geodesic.\\
\end{enumerate}
\end{corollary}

\textbf{Proof} \quad We start with $V:=V_{\ve}$ and the metric $g_{\T}$ from \ref{merg2} so that $scal (g_{\T}) > 0$ on $H \setminus V_{\ve}$ and $\p V_{\ve}$ has positive mean curvature. Now the elementary non-conformal deformation of Almeida [A], Th.1.1, supported arbitrarily close to $\p V_{\ve}$, can be applied to transform this metric into another one $g^*_{\ve}$ on $H \setminus V_{\ve}$ so that
 \begin{equation}\label{ccc2}
 scal (g^*_{\ve}) > 0 \mm{ on } H \setminus V_{\ve} \mm{  and }   \p V_{\ve}  \mm{ is totally geodesic.}
 \end{equation}
 In fact $g^*_{\ve}$ can be chosen sufficiently regular along $\p V_{\ve}$ so that on the doubling the resulting metric $g^*_{\ve}$ on $H^D_{\ve} = H \setminus  V_{\ve}  \cup_\sim H\setminus  V_{\ve}$ is smooth. \qed

\begin{corollary} \emph{\textbf{(Smooth Ambience II)}} \quad Let $\beta \in H_{n-1}(H;\Z)$ be given. Then, for any $\ve
> 0$, we can find a hypersurface $N_{\beta,\ve} \subset H \setminus  V_{\ve}$, so that
\begin{enumerate}
\item  $N_{\beta,\ve}$ is an integral current and mass minimizing under perturbations with support in the interior of $H_{\ve}$
\item $N_{\beta,\ve}$ hits $\p V_{\ve}$ perpendicularly, $\p N_{\beta,\ve} \subset \p V_{\ve}$  is a minimal hypersurface and $\p N_{\beta,\ve} \subset N^D_{\beta,\ve}$ is
    totally geodesic.
\item  $N_{\beta,\ve}$  represents the image of $\beta$ under the projection map $\pi_*$ composed with the isomorphism ${I_{H,V_{\ve}}}$ between $H_*$ and $H^{IC}_*$ on
    the couple  $(H,V_{\ve})$: {\small
\[\begin{CD}
\beta \in H_{n-1}(H;\Z) @>{\pi_*}>>  H_{n-1}(H,V_{\ve};\Z) @>{I_{H,V_{\ve}}}>> [N_{\beta,\ve}] \in H^{IC}_{n-1}(H,V_{\ve};\Z).
\end{CD}\]}
\end{enumerate}
\end{corollary}

\textbf{Proof} \quad Given $\beta \in H_{n-1}(H;\Z)$, the image $\pi_*(\beta) \in  H_{n-1}(H,V_{\ve};\Z)$ can be represented by a hypersurface $L_{\beta,\ve}$ with smooth
boundary $\p L_{\beta,\ve} \subset \p V_{\ve}$. We take the doublings $L_{\beta,\ve}^D \subset H^D_{\ve}$ and can find a mirror-symmetric area minimizing hypersurface $L^\cs_{\beta,\ve}  \subset H^D_{\ve}$ homologous to $L_{\beta,\ve}^D$. Then $N_{\beta,\ve}:=L^\cs_{\beta,\ve} \cap H \setminus  V_{\ve}$ satisfies the asserted properties.
 \qed

\setcounter{section}{5}
\renewcommand{\thesubsection}{\thesection}
\subsection{Isoperimetric Excursions} \label{smo}
\bigskip

In this chapter we discuss the technical background needed to control area minimizing hypersurfaces with $H$ relative to the horizon metric $g_{\T}$. In essence, we combine different forms and applications of isoperimetric inequalities. An interesting aspect is the relation between these inequalities and the central asymptotic property of area minimizers we exploited in this series of papers, their skin uniformity. \\

This chapter is written as a supplement, contrasting the constructive chapters \ref{fm0} - \ref{smo0}. We restrict to indications of the main ideas and hints to references. We leave details, like the adaptation of well-documented techniques from the literature or inductive estimates of volumina, to the reader.

\subsubsection{Geometry of Minimal Metrics} \label{geom}
\bigskip

When we span area minimizers over local horizons $(\D_{\CCC})_*(\L_{\CCC})+ g^+_H$ we know different aspects of this geometry near the support of the inserted pseudometric
$(\D_{\CCC})_*(\L_{\CCC})$. But to understand its impact on area minimizers we want to span over the wedge the geometry $(\D_{\CCC})_*(\L_{\CCC})+ g^+_H$ has when compared to $g^+_H$, we also need to understand the minimal metric we think of a indifferent background metric.\\

 We start with a Plateau problem relative $g^+_H$, we use to derive \ref{haee} and \ref{l3}.

 \begin{lemma} \label{els}\, For  $C^n= \R^{n-k} \times C^k$, $C^k \in {\cal{SH}}^{\R}_k$, let $P^{n-1}(\zeta)$ be an area minimizing  solution of the following Plateau problem in $(C,g^+_{C^n})$:
\begin{equation}\label{plp}
P^{n-1} \subset C^n = \R^{n-k} \times C^k \,\mm{ with }\, \p P= \p B^{n-k}_1(0) \times \p B^k _\zeta(0)
\end{equation}
for some $\zeta >0$. Then, for any $\ve>0$, there is a $\zeta_0(\ve) >0$, independent of $C$, so that for any $\zeta \in (0,\zeta_0(\ve))$:
\begin{equation}\label{ppl}
P^{n-1} \cap  B_{1-\ve}(0) \times C^k \v, \mm{ and }
\end{equation}
\begin{equation}\label{ppl2}
P^{n-1} \cap  \R^{n-k} \times (C^k \setminus \overline{B^k _\zeta(0)}) \v.
\end{equation}
\end{lemma}
\medskip

\textbf{Proof} \quad We start with (\ref{ppl}) and  observe that for any cone $C^n= \R^{n-k} \times C^k$, $C^k \in {\cal{SH}}^{\R}_k$, the area $Vol_{n-1}(\p B^{n-k}_1(0) \times B^k_\zeta(0),g^+_C)$ is \emph{finite}. We start with the case of a cone $C^k$ singular only in $0$ and $n \ge k \ge 7$
\begin{equation}\label{v1}
 Vol_{n-1}(\p B^{n-k}_1(0) \times B^k_\zeta(0),g^+_C)\le
\end{equation}
 \[c_C\cdot  \int_{B^k_\zeta(0) } (r^{- (1- \sqrt{2/3}) \cdot \frac{k-2}{n-2}})^{n-k-1} \cdot (r^{- (1- \sqrt{2/3}) \cdot \frac{k-2}{n-2}})^{k-1} \cdot r^{k-1} dr d \theta_{k-1}\le \]
\[c^*_C\cdot  \zeta^{- (1- \sqrt{2/3}) \cdot (k-2)+k}\]
\begin{equation}\label{v2}
 Vol_{n-1}(B^{n-k}_1(0) \times \p B^k_\zeta(0),g^+_C)\ge \overline{c_C}\cdot   (\zeta^{- (1- \sqrt{2/3}) \cdot \frac{k-2}{n-2}})^{n-1}  \cdot  \zeta^{k-1}
 \end{equation}
 \[= \overline{c_C}\cdot   \zeta^{- (1- \sqrt{2/3}) \cdot \frac{(n-1)\cdot (k-2)}{n-2}+k-1}.\]
For  $c_C:=\sup_{\omega \in S_{C^k}\setminus \Sigma_{S_{C^k}}}  \psi_{C^k}(\omega)>0$ and $\overline{c_C}:=\inf_{\omega \in S_{C^k}\setminus \Sigma_{S_{C^k}}}  \psi_{C^k}(\omega)>0$, depending on $C$. Both  volumina shrink to zero, while $\zeta \ra 0$, since the exponents are $\ge k/2 >0$.\\

In dimension $7$ we readily see that there is common upper resp. lower bounds for the constant $c_7 \ge c_C$ and $\overline{c_7} \le \overline{c_C}$, for any $C \in \mathcal{SC}_7$ from the  compactness of $\mathcal{SC}_7$.
 In higher dimensions, where we expect singularities also outside $0$, we argue the same way for the pencils $B_1(0) \cap\P(0,\rho) \subset C$, for $\rho >0$.  Using the stable covers of \ref{bow0} we inductively get uniform area estimates for $B_1(0) \cap\P^*(0,\rho) \subset C$, for sufficiently small $\rho >0$, from the compactness of $\mathcal{SC}_k$.\\

For any $\delta >0$, there is a $\eta_\delta >0$ for any $C \in \mathcal{SC}_n$, so that for $\eta  \in (0,\eta_\delta]$:
\begin{equation}\label{rag}
 Vol_{n-1}(\p B^{n-k}_1(0) \times B^k_\eta(0),g^+_C)/ Vol_{n-1}(B^{n-k}_1(0) \times \p B^k_\eta(0),g^+_C) \le \ve
\end{equation}
From this estimate and the coarea formula e.g. [GMS], Vol1, Ch.2.1.5, Th.3, p.103, we infer  (\ref{ppl}). \\

Finally, (\ref{ppl2}) follows readily from the fact that the repelling effect of the distance tubes  towards $\{0\} \times C^k$, we experience from $tr\, A_{T_\varrho}(g^+_{C^n})  <0$ from \ref{etl}.  \qed

\textbf{Proof of \ref{haee} and \ref{l3}} \quad  We only need to consider a single cone $C^n = \R^{n-k} \times C^k$, $C^k \in {\cal{SC}}_k$, the usual compactness arguments for ${\cal{SC}}_k$ give us uniform bounds for all such cones.\\

We first estimate the area of any rectifiable hypersurface $S \subset (C^n,g^{\bigodot}_{C^n})$ that bounds an open set containing $B^{n-k}_1(0) \times B^k_b(0)\setminus B^k_a(0)$, for some $b>a>0$:
\[Vol_{n-1}(S) \ge Vol_{n-1}(\p B^{n-k}_1(0) \times B^k_b(0)\setminus B^k_a(0)) \ge \]
 \[k_C \cdot  \int_{B^k_b(0)\setminus B^k_a(0)} (r^{- (1+ \sqrt{2/3}) \cdot \frac{k-2}{n-2}})^{n-k-1} \cdot (r^{- (1+ \sqrt{2/3}) \cdot \frac{k-2}{n-2}})^{k-1} \cdot r^{k-1} dr d \theta_{k-1}\ge\]
\[k^*_C \cdot  \Big(a^{- (1+ \sqrt{2/3}) \cdot (k-2)+k} - b^{- (1+ \sqrt{2/3}) \cdot (k-2)+k}\Big), \mm{ where the exponent is }  \le - k/2.\]
For some constant $k_C>0$. Thus, for any given $b >0$ we can find a small $a \in (0,b)$ so that $Vol_{n-1}(S)$ exceeds any given lower bound.\\

The constant $k_C$ results from the global Harnack inequality \ref{skk}: recall that for $\Psi_-[n,k](\omega,r) = \psi_C(\omega) \cdot r^{\alpha_-}$:
$|\psi_C|_{L^1(S_C \setminus \Sigma_{S_C})} \le a_{n,k,\lambda} \cdot \inf_{\omega \in S_C \setminus \Sigma_{S_C}}  \psi_C(\omega)$. This Harnack estimate shows that $\Psi_-[n,k](\omega,r) \ge  r^{\alpha_-}/a_{n,k,\lambda}$, for $|\psi_C|_{L^1(S_C \setminus \Sigma_{S_C})} = 1$.\\

Thus, a comparison with \ref{els} (\ref{v1}) shows that, for $R/r$ large enough, there is a connected area minimizing $F \subset Q^{n-k}(L^*) \times B_{R^*}(0)$ relative $g^\Delta_{C^n}$ so that
\begin{itemize}
\item $F$ bounds an open set that contains  $Q^{n-k}(L-\delta/2)  \times B_{\frac{2 \cdot R+ r}{3}}(0) \setminus B_{\frac{R+ 2 \cdot r}{3}}(0)$
\end{itemize}

Now we use \ref{els} for the boundary $Q^{n-k}(L-\delta/2)  \times \p B_{\frac{R+ r}{2}}(0)$ to also see that $F$ surrounds $Q^{n-k}(L-\delta)  \times B_{\frac{R+ 2 \cdot r}{3}}(0)$.\\

 The remaining assertions follow from the repelling effect of the non-negative mean curvature of $\p (Q^{n-k}(L^*) \times B_{R^*}(0))$ relative $g^+_{C^n}$ and the compactness results for area minimizers. The argument then carries over to the tapered product \ref{l3}. \qed

\textbf{General Minimal Metrics} \, Now  we consider the shielding effects of the $\D$-transferred pseudometrics on the less symmetric approximated space.  Due to the almost isometric approximation via $\D$-maps the shielding argument carries over - up to a channel near the singularities, where $\D$-maps no longer apply. In our covers this channel appears as a subset of some singular pencil $\P^*(z,\rho)$, for some very small $\rho>0$. For a high accuracy $\widehat{\ve} \ll 1$,  this channel clearly becomes arbitrarily narrow relative to the size of the surrounding $\bp$-pipes.\\

 As before, we have to show that an area minimizing hypersurface does not enter such sufficiently narrow channels to any substantial extent.  The advantage of this two step argument is that we can now focus on area minimizing hypersurfaces within area minimizers $H \in {\cal{H}}^{\R}_n$, equipped with their minimal metric $(H, g^+_H)$. \\

The proofs of such results typically involve variants of isoperimetric inequalities. In the next sections, we derive an appropriate isoperimetric inequality for $(H, g^+_H)$. Here we start with some preliminary observations concerning the basic geometry of $(H, g_H)$ and $(H, g^+_H)$.

 \begin{lemma}  \emph{\textbf{(Diameter and Volume Bounds)}} \label{abb} \, We have  the following  bounds
 \begin{enumerate}
   \item \, For any $H \in {\cal{H}}^{\R}_n$, $p \in H^n$ we have
{\small \begin{equation}\label{sii}
\kappa_n(0) \le Vol_{n-1}(\p B_1(p),g^+_H) \le \kappa_n(1) \mm{ and } \kappa^*_n(0) \le Vol_{n}(B_1(p),g^+_H) \le \kappa^*_n(1)
\end{equation}}
for some $\kappa^*_n(i) , \kappa_n(i)>0$, $i=0,1$, independent of $H$.
 \item \, For any $H \in {\cal{H}}^c_n$, we have:
    \[diam (H \setminus B_\rho(q_H),\Phi(H)^{4/(n-2)} \cdot g_H) < \infty, \mm{ for any radius } \rho>0.\]
 \end{enumerate}
  \end{lemma}
For the proof one argues using the usual induction and compactness theorems, the freezing effect and as in  \ref{els} using the growth rate of $\Psi_+$ on cones estimated in \ref{ddd}(iii).  We leave the details to the reader.\qed

\begin{lemma} \label{confin} \quad  For $C \in \mathcal{C}_n$, the deformed space  $(C, g^+_C)$ is again a cone. That is, $(C,g^+_C)$ is invariant under scalings around $0$.
This equally applies to subcones of $C$.\\

Also, for the pencil $\P(0,\rho)$ we have some constant $c_n(\rho)\ge 1$, depending only on the dimension so  that for the respective  spherical measures 
\begin{equation}\label{aps}
c_n^{-1} \cdot \mu(\P(0,\rho)\cap S_C,g_C) \le \mu(\P(0,\rho)\cap S_C,g^+_C) \le c_n  \cdot \mu(\P(0,\rho)\cap S_C,g_C)
\end{equation}
We call this the \textbf{aperture stability} of the transition between $g^+_C$ and $g_C$.
\end{lemma}

In particular, we find that the tangent cone approximation of  $(H, g_H)$ by cones $C,g_C$ near singular points, including the freezing effects, translates to completely similar cone approximations of   $(H, g^+_H)$ by cones $C,g^+_C$.\\

{\bf Proof} \quad As before, we note that for $\Psi_+(\omega,r)= \psi_C(\omega) \cdot r^{\alpha_+}$, we have $\alpha_+ >- (1- \sqrt{2/3}) \cdot \frac{k-2}{2}$ and therefore $r^{2 \cdot\alpha_+/(n-2)}$ is integrable on $[0,1]$.\\

To see that $(C, g^+_C)$ is again a cone, we first consider a conformal deformation by $r^{2\cdot \beta} \cdot g$, $\beta > -1$. Namely, the distance from $0$ to a
point in $\p B_\rho(0)$ becomes $\int_0^{\rho} r^\beta = \frac{1}{\beta  +1 } \cdot \rho^{\beta  +1 }$ and the diameter of $\p B_\rho(0) \cap C$ becomes $\rho^{\beta + 1}\cdot
diam(\p B_1(0) \cap C,g) $. Choosing a reparameterization $\R^{\ge 0}$ by arc-length we see that $(C,r^{2\cdot \beta} \cdot g)$ is again a cone,  isometric to
 $(\R^{\ge 0} \times \p B_1(0) \cap C, g_{\R} + c_\beta \cdot r^2 \cdot g_{\p B_1(0) \cap C})$ for some $c_\beta > 0$.\\

Next we consider deformations depending only on $\omega$: $g' = c(\omega) \cdot g$ for a conical metric $g$. Let $f: C \to C$ be  the map $(r,\omega) \mapsto (\lambda \cdot r, \omega)$, for some $\lambda > 0$.  $g$ satisfies $ f^*g = \lambda^2 \cdot g$ and $g'$ shares this property:
{\small \[(f^* g')_{(r,\omega)} (X,Y) = c(f(r,\omega)) \cdot g_{(f(r,\omega))} (df_{(r,\omega)} X, df_{(r,\omega)} Y) = c(\omega) \cdot f^*g (X,Y) = \lambda^2 \; g'(X,Y).\]}
We combine the two arguments to see that $(C, g^+_C)$  is isometric to any copy scaled around $0$.\\

Finally, the asserted aperture stability is again an application of the usual compactness arguments on $\mathcal{C}_n$.
\qed

\textbf{Conclusion} \, The minimal $scal>0$-geometry $(H, g^+_H)$ is situated between the original geometry  $(H, g_H)$ and the Gromov hyperbolic skin geometry $(H, \bp^2 \cdot g_H)$. The space $(H, g^+_H)$ is compact. Towards singular points the space looks more and more like a $scal>0$-cone with a positively bounded aperture strictly smaller than that of tangent cones of $(H, g_H)$. In simple terms, $(H, g^+_H)$ is a $scal>0$-version of $(H, g_H)$, uniformly conformally expanded towards $\Sigma$.\\

\subsubsection{Isoperimetry and Skin Uniformity} \label{wm}
\bigskip

On area minimizing hypersurfaces, and it does not matter whether they are regular or singular,  in a smooth ambience we have isoperimetric inequalities. Specifically, there is some $\gamma>0$  depending only on $n$, so that for any
$H \in {\cal{H}}^{\R}_n$ and any open set $U \subset H$ with compact closure and rectifiable boundary $\p U$ we get the estimate:
\begin{equation}\label{ii}
 Vol_{n}(U)^{(n-1)/n} \le \gamma \cdot Vol_{n-1}(\p U)
\end{equation}
There are several classical proofs of (\ref{ii}) giving differently sharp estimates for $\gamma>0$.  The arguably simplest idea, that applies to area minimizers in arbitrary codimension,  is to take a cone over $\p U$. This gives us a coarse upper estimate for $Vol_{n}(U)^{(n-1)/n}$. To derive uniform bounds for $\gamma$, for any $U$,  one may use Federer's deformation theorem instead to find more efficient but basically similar extensions of $\p U$, cf.[F], 4.2.10.\\

Here we draw another proof that directly applies to $(H, g_H)$, but extends to $(H, g^+_H)$, from the skin uniformity of $(H, g_H)$, cf.Ch.\ref{mre}.C. This argument relies on more recent techniques due to Semmes [Se]. \\

\textbf{Thick Families of Curves} \, The upshot of [Se] is that a metric measure space $(X,d,\mu)$ with the property that for any pair of points $p,q \in X$  there is a  controllably thick spindle $\Gamma_{p,q}$ of rectifiable curves $\gamma: I_\gamma \ra X$, $I_\gamma \subset \R$, parameterized by arc-length joining them, with lengths comparable to the distance of the points, support  Sobolev and thus isoperimetric  inequalities. We recollect the needed \emph{Semmes conditions}:
\begin{enumerate}
 \item There exists a constant $C > 0$ so that any $\gamma \in \Gamma_{p,q}$: $l(\gamma|[s,t]) < C \cdot |s-t|$, for any two $s,t \in I_\gamma$.
 \item Each $\Gamma_{p,q}$ carries a probability measure $\sigma_{p,q}$.
 \item For any Borel set $A \subset X$, the assignment $\gamma \mapsto$  is $\sigma$-measurable with
 \begin{equation}\label{tcu}
 \int_{\Gamma_{p,q}} l(\gamma \cap A) \, d \sigma(\gamma) \le  C \cdot \int_{A_{C,p,q}} \left(\frac{d(p,z)}{\mu(B_{d(p,z)}(p)} + \frac{d(q,z)}{\mu(B_{d(q,z)}(q)}\right) d \mu(z)
 \end{equation}
   for $A_{C,p,q}:=(B_{C \cdot d(p,q)}(p) \cup B_{C \cdot d(p,q)}(q))\cap A$.
 \item $\mu$ is doubling: $\mu(B_{2 \cdot r}(q) \le C^*\cdot \mu(B_{2 \cdot r}(q)$, for any $q \in X$ and $r>0$.
 \item For balls $B^* \subset B \subset X$ and some $Q>0$, $C_Q>0$, we have
\[(diam(B^*)/diam(B))^Q \le C_Q \cdot \mu(B^*)/\mu(B).\]
 \end{enumerate}
At first sight, the conditions, in particular (iii), may appear somewhat cryptic. However, they are nothing but the plain components of the classical proof for the Poincar\'{e} inequality, viewed as an efficient doubled average. The last condition (v) is the one that allows us to derive the Sobolev inequality from the Poincar\'{e} inequality.
For further details cf.[Se] [Se], [He],Ch.4,p.29 and [HKST]. Ch.9.1 and 14.2. \\

Now, the point is that the skin uniformity of $(H, g_H)$, implies these conditions for thick families of curves relative $(H, g_H)$, and also relative $(H, g^+_H)$. (Once again, the weaker, while more common, uniformity would not be sufficient, cf. the discussions in [L1],Ch.1.1 and 4.2, and in [L2], after Th.1.)

\begin{proposition} \emph{\textbf{(Isoperimetric Inequalities)}}\label{isos} \,  There is some $\gamma, \gamma^*>0$  depending only on $n$, so that for any
$H \in {\cal{H}}^{\R}_n$ and any open set $U \subset H$ with compact closure and rectifiable boundary $\p U$ we get,  relative $(H, g_H)$ but also $(H, g^+_H)$:
\begin{equation}\label{iif}
 Vol_{n}(U)^{(n-1)/n} \le \gamma \cdot Vol_{n-1}(\p U)
\end{equation}
\begin{equation}\label{ii2}
\min \{Vol_{n}(B_{\rho} \cap U),Vol_{n}(B_{\rho} \setminus U)\}^{(n-1)/n} \le \gamma^* \cdot Vol_{n-1}(B_{\rho} \cap \p U)
\end{equation}
\end{proposition}

\begin{remark}\,  In [L1], we invested a (slightly weaker) version of (\ref{ii2}) for $(H, g_H)$ to derive its skin uniformity, namely [BG],Th.2. Thus skin uniformity may be regarded as another way to express this inequality, that also gives us means for improvements, for instance, in the extension to $(H, g^+_H)$.
 \end{remark}

\textbf{Proof} \quad The scaling invariance of the equation $(L_{C^k,n})_\lambda \phi =0,$ and  the uniform volume estimates of \ref{abb}  readily shows  (\ref{iif}) in the case of balls:  for any $H \in {\cal{H}}^{\R}_n$, radius $R>0$ and any point $p \in H$:
 \begin{equation}\label{iao}
 Vol_{n}((B_R(p), g^+_H))^{(n-1)/n} \le \gamma^+ \cdot Vol_{n-1}((\p B_R(p), g^+_H))
 \end{equation}
for some $\gamma^+>0$ depending only on $n$ and $\lambda$.  To derive the general form of (\ref{iif}) and (\ref{ii2}), we first note that for any given $R>0$ the geometries $(B_R(p), g^+_H)$ form a \emph{compact family} ${\cal{F}}_R$, that transforms into ${\cal{F}}_{R^*}$ under scaling by $R^*/R>0$.\\

Now we apply the skin uniformity of $(H,g_H)$. For any two $p,q \in  (B_{3 \cdot R/4}(p), g^+_H)$, with distance $R$ we have a skin uniform curve $\gamma: [a,b] \ra H$, for some $a <b$, with $\gamma(( a,b))\subset H \setminus \Sigma$, so that for some fixed $s_H>0$
  \begin{enumerate}
    \item \quad  $l(\gamma)  \le s_H \cdot d_{g_H}(p,q)$,
    \item \quad $l_{min}(\gamma(z)) \le s_H \cdot \delta_{\bp}(z)$, for any $z \in \gamma_{p,q}.$
    \end{enumerate}
We may improve the first condition to also holds for any two points on the curve and the segment between them. This is made to ensure Semmes condition (i). (For Poincar\'{e} inequalities on metric spaces, this plays an important role and is called called a quasi-convexity condition.)\\

In the canonical hyperbolic unfolding $(H,d_{\bp})$, [L2], Th.2, we choose $\gamma$ as a hyperbolic geodesic arc joining $p$ and $q$ in this complete geodesic space. These curves are again skin uniform in $(H,g_H)$ for some  $s^*_H\ge s_H>0$, cf. [L2], Lemma 2.14. That is,  we get $s^*_H$-skin uniform and quasi-convex curves linking any two given points.\\

Next we use the second skin uniformity condition, the Lipschitz continuity of $\delta_{\bp}$ and $|A| \le \bp$ to argue that we have a fine controls over the geometry along each skin uniform curve and in twisted double cone around it. This can be used to define carries the probability measure $\sigma_{p,q}$ of (ii) on $\Gamma_{p,q}$ we choose as a uniformly sized entourage of $2 \cdot s^*_H$-skin uniform curves for $g_H$ and from the aperture stability (\ref{aps}) also for $g^+_H$ and to establish (iii) using the arguments of [Se], Remark 6.10.\\

We leave this along with a verification of the condition (iv) and (v), for $Q=n$, and  to derive (\ref{iif}) and (\ref{ii2}), from [HKST]. Ch.9.1, to the reader. \qed

Our application is the proof of the following counterpart of \ref{els}.

\begin{corollary} \emph{\textbf{(Narrow Channels)}}\label{lue} \,  For any $\ve >0$, there is a $\delta >0$ depending only on $n$ and $\lambda$ so that for any
$H \in {\cal{H}}^{\R}_n$, there is an area minimizing  solution of the following Plateau problem in $(H,g^+_H)$:
\begin{equation}\label{nc0}
 P^{n-1} \subset H^n \,\mm{ with }\, \p P= \p B_1(z) \cap \P^*(z,\delta), \mm{ so that } P^{n-1} \subset  B_{1+\ve}(z) \setminus B_{1-\ve}(z)
\end{equation}
\end{corollary}

\textbf{Proof} \quad The volume bounds of \ref{abb} show that $Vol_n(B_1(z) \cap \P^*(z,\delta),g_H^+) \ra 0$, for $\delta \ra 0$, uniformly in $H \in {\cal{H}}^{\R}_n$ and $z \in H$.\\

From  (\ref{ii2}), we get as in [Gi], 5.14, some $\iota_n>0$ depending only on the dimension $n$ so that for any
$H \in {\cal{H}}^{\R}_n$ and any connected area minimizing hypersurface  $F^{n-1}\subset H^n$ bounding an open set $F^+\subset H^n$ the following estimate holds for any $r >0$
\begin{equation}\label{ii3}
Vol_n(F^+ \cap B_r(p)) \ge \iota_n \cdot  r^n.
\end{equation}
Thus we see that for sufficiently small $\delta>0$, the solution of the stated Plateau problem remains in $B_{1+\ve}(z) \setminus B_{1-\ve}(z)$. \qed

\subsubsection{Superposition Effects}\label{smoo}
\bigskip

In this concluding section,  we merge the properties and controls of our nanostructure $\sp$ with those we have obtained for individual elements above, to derive \ref{merg2}. \\

The basic idea is to show that a slightly perturbed union of those regions shielded by these elements is also shielded by $g_{\T}$. We henceforth we call ${\cal{E}}^{tap}_{C^n}$ the \emph{shielded region} of the element that belongs to $C^n$ of \ref{l3}.\\

 \textbf{Proof of \ref{merg2}} \quad We follow the stepwise assembly of multi-pipes which was part of their definition in Ch.\ref{grq} and recall that there are two types of superpositions of tapered horizon metrics alternatingly applied while the  induction runs over the height $\mathbf{h}$. We start from terminal nodes until we reach  $\sp \subset H \setminus \Sigma$. Also recall from \ref{tz} and \ref{tst} that also the almost tangent cone usually belong to lower-dimensional spheres in cones, we finally consider the, multiply taken, tapered products in the final $\sp \subset H$.\\

$\mathbf{h}=0$ \, The induction starts with multi-pipes of height $0$. This case is just \ref{l3} combined with \ref{lue} and directly applies to multi-pipelines of height $0$ since they constitute merely  deformations of disjoint multi-pipes, also after taking the tapered product and $\D$-map transfers,  they lead to disjoint  shielded regions.\\

$\mathbf{h}\Rightarrow \mathbf{h+1}$ \, In the induction step we assume the result is known for any multi-pipeline $\sp_\mathbf{h}$ of height $\mathbf{h}$. An also elementwise application of \ref{l3} gives the result for
$\st_{\mathbf{h}+1}(\sp_{\mathbf{h}})_{spokes}$ that is $C(\sp_\mathbf{h})  \cap B_{5 \cdot s_{m+1}^i}(0)\setminus  B_{\pi_{m+1}[j,1] \cdot s_{m+1}^i}(0)$
and for the hub  $\st_{\mathbf{h}+1}(\sp_{\mathbf{h}})_{hub}$.\\

Thus we need to understand the effect on the shielded regions under these two operations:\\

$\bullet$ \, We first consider the merging process of hubs and spokes to reach the  multi-pipe  $\st_{\mathbf{h}+1}(\sp_{\mathbf{h}})$ of height $\mathbf{h} +1$
\[\st_{\mathbf{h}+1}(\sp_{\mathbf{h}})_{spokes} \, \cup \, \st_{\mathbf{h}+1}(\sp_{\mathbf{h}})_{hub}\, \mm{ with } \,\st_{\mathbf{h}+1}(\sp_{\mathbf{h}})_{spokes} \, \cap \, \st_{\mathbf{h}+1}(\sp_{\mathbf{h}})_{hub}=...\]
\[...= C(\sp_\mathbf{h})  \cap B_{\pi_{m+1}[j,0]\cdot s_{m+1}^i}(0)\setminus   B_{\pi_{m+1}[j,1]\cdot s_{m+1}^i}(0).\]
Due to the symmetry of the metrics on the spokes, the area minimizer spanned over the hub, will remain to be a hypersurface of mean curvature zero  and vice versa. They intersect perpendicularly along their boundaries.\\

$\bullet$ \, In the second part of the induction step, we assemble the families of multi-pipes of height $\mathbf{h}$ to a multi-pipelines of height $\mathbf{h}$. This is the portion of the argument that follows not from an isolated look at the multi-pipes,  but on their interaction. In turn this depends on the parameters of separator tree $\T$.\\

We choose, in \ref{bow3} and \ref{set}, a separator tree $\T$ for prescribed $\widehat{\Theta}^\bullet,\widehat{\ve}^\bullet,\widehat{s}^\bullet$, where we use very small self-similarity factors $\widehat{s}^\bullet \ll 1$, accuracy values $\widehat{\ve}^\bullet \ll 1$  and \emph{small} but \emph{broad} separators $\widehat{\Theta}^\bullet$, that is,  with  $\pi^\bullet_{n+1}[\kappa,1]/\pi^\bullet_{n+1}[\kappa+1,0]\gg 1$, $\kappa=1,...,c(n+1)$.\\

We first consider a single layer of balls ${\cal{A}}^\cs[i]$, for one fixed $i$. It splits into  disjoint sublayers ${\cal{A}}^\cs[i,1],...,{\cal{A}}^\cs[i,c]$ with \[B_{10 \cdot s^i}(p)
     \cap B_{10 \cdot s^i}(q) \v, \mm{ for } B_{2 \cdot s^i}(p), B_{2 \cdot s^i}(q), \mm{ in the same sublayer }{\cal{A}}^\cs[i,j].\]
We use an induction over $j$ and  first place the multi-pipes belonging to balls in ${\cal{A}}^\cs[i,1]$. They and their shielded regions do not interfere. Thus we can proceed to the
essential induction step assuming we placed the multi-pipes belonging to balls in ${\cal{A}}^\cs[i,k]$, $k \le j$, so that the shielding hypersurfaces are surrounded by a new common area minimizer up to some small overlap $\eta_j >0$, with $\eta_j \ra 0$, for $\pi^\bullet_{n+1}[\kappa,1]/\pi^\bullet_{n+1}[\kappa+1,0]\ra \infty$.\\

For  $\pi^\bullet_{n+1}[\kappa,1]/\pi^\bullet_{n+1}[\kappa+1,0]\gg 1$ sufficiently large we observe that placing the multi-pipes belonging to balls in ${\cal{A}}^\cs[i,j+1]$
 the shielding hypersurfaces are now surrounded by a new common area minimizer up to some small overlap $\eta_{j+1} > \eta_j >0$. Since $j \le c(n)$, we can ensure that the resulting hypersurfaces surround $\bigcup_{\CCC \in \T} {\cal{O}}_{\CCC}$. \\

Now we consider the interference of multi-pipes in different layers: for sufficiently small $\widehat{s}^\bullet$, the multi-pipes layer ${\cal{A}}^\cs[i]$ and their shielding does not even recognize the presence of the larger ones, that is, ${\cal{A}}^\cs[j]$, $j<i$. In turn, the narrow channel argument of the last section shows that the effect of smaller sized layers ${\cal{A}}^\cs[j]$, $j>i$ does not push the area minimizer beyond their own support.\\

Finally, we invoke the singular threshold \ref{strs} and the surrounding relation \ref{s} of Ch.\ref{thr}, to see that, after another refinement of the tree parameters we also have $scal(g_{\T})>0$ outside  $\bigcup_{\CCC \in \T} {\cal{O}}_{\CCC}$. Here we need a minor adjustment. Inductively, over $\mathbf{h}$,  we slightly but uniformly extend the spokes in radial directions so that their shielded region surrounds the places of potential $scal <0$ of their associated hubs.\\

 For the terminal nodes this clearly does not apply, but it is also not needed since their are no singular direction where the localization deformation would cause any $scal <0$-contributions. \qed

\end{document}